\documentclass{article}
\usepackage[english]{babel}
\usepackage{amsthm}
\usepackage{amsmath}
\usepackage{mathtools}
\usepackage{parskip}
\usepackage{dsfont}
\usepackage{cleveref}
\usepackage{enumerate}
\usepackage{amsfonts}
\usepackage{graphicx} 
\usepackage{todonotes}
\usepackage{adjustbox}
\usepackage{fullpage}
\usepackage{bm}
\usepackage{subcaption}
\usepackage{enumitem}

\newtheorem{theorem}{Theorem}[section]
\newtheorem{claim}[theorem]{Claim}
\newtheorem{lemma}[theorem]{Lemma}

\newtheorem{remark}[theorem]{Remark}

\crefname{claim}{Claim}{Claims}

\catcode`\@=11
\@addtoreset{equation}{section}

\title{The evasion of tipping: pattern formation near a Turing-fold bifurcation}

\author{Dock Staal\footnotemark[1] {} and Arjen Doelman\footnotemark[1]}
 
\date{\today}

\begin{document}

\maketitle

\renewcommand{\l}{\left(}
\renewcommand{\r}{\right)}

\renewcommand{\thefootnote}{\fnsymbol{footnote}}
\footnotetext[1]{Mathematisch Instituut, Einsteinweg 55, 2333 CC Leiden,  Universiteit Leiden, the Netherlands}
\renewcommand{\thefootnote}{\arabic{footnote}}

\begin{abstract}
Model studies indicate that many climate subsystems, especially ecosystems, may be vulnerable to {\it tipping}: a {\it catastrophic process} in which a system, driven by gradually changing external factors, abruptly transitions (or {\it collapses}) from a preferred state to a less desirable one. In ecosystems, the emergence of spatial patterns has traditionally been interpreted as a possible {\it early warning signal} for tipping. More recently, however, pattern formation has been proposed to serve a fundamentally different role: as a mechanism through which an (eco)system may {\it evade tipping} by forming stable patterns that persist beyond the tipping point. 

Mathematically, tipping is typically associated with a saddle-node bifurcation, while pattern formation is normally driven by a Turing bifurcation. Therefore, we study the co-dimension 2 Turing-fold bifurcation and investigate the question: {\it When can patterns initiated by the Turing bifurcation enable a system to evade tipping?} 

We develop our approach for a class of phase-field models and subsequently apply it to $n$-component reaction-diffusion systems -- a class of PDEs often used in ecosystem modeling. We demonstrate that a two-component system of modulation equations governs pattern formation near a Turing-fold bifurcation, and that tipping will be evaded when a critical parameter, $\beta$, is positive. We derive explicit expressions for $\beta$, allowing one to determine whether a given system may evade tipping.  Moreover, we show numerically that this system exhibits rich behavior, ranging from stable, stationary, spatially quasi-periodic patterns to irregular, spatio-temporal, chaos-like dynamics.
\end{abstract}


\paragraph{Highlights}
\begin{itemize}[label=--]
    \item Investigated the interactions of tipping and pattern formation by studying co-dimension 2 Turing-fold bifurcations.
    \item Derived a new coupled system of modulation equations -- the AB-system -- that captures the associated dynamics.
    \item Analyzed the existence and stability of periodic solutions in the AB-system -- i.e., at the tip of the Busse balloon.
    \item Greatly reduced the complexity involved in calculating the sign of the Landau coefficient.
    \item Showed that tipping is evaded if the Turing-fold bifurcation -- equivalently, the Turing bifurcation -- is supercritical.
    \item Demonstrated the possibility of complex, chaos-like dynamics near the Turing-fold bifurcation.
\end{itemize}

\section{Introduction}
\label{sec:introduction}

In studies on the impact of climate change, the phenomenon of tipping plays an increasingly central role: various climate subsystems have been designated as {\it tipping elements} \cite{lenton2008tipping}. These elements are (typically large-scale) subcomponents of Earth's climate system, which have been indicated by (certain) models  to be at a risk of pass through a tipping point -- i.e., a critical parameter combination at which seemingly environmental changes may push the subsystem from a favorable state into a much less desirable one. Since this tipping process is expected to take place on a time scale relatively fast compared to that of the driving environmental changes, tipping is often associated with {\it catastrophic collapse}.  These potential {\it catastrophic shifts}, and especially their impact on ecosystem resilience, have long been a key topic in ecology \cite{scheffer2001catastrophic}. Naturally, identifying (observable) indicators of proximity to {\it critical transition} for complex system has been a central (sub)theme within the literature on tipping \cite{scheffer2009early}. In ecosystems, the appearance of {\it self-organized patchiness} has been proposed to be such an {\it early warning signal} \cite{rietkerk2004self}. More recently, however, the (potential) role of vegetation patterns in ecosystems and the associated {\it multi-stability} has been interpreted differently. Instead of being a signal of an impending collapse, pattern formation may increase the resilience of an ecosystem, allowing an ecosystem to {\it evade tipping} by gradually transforming into a patterned state \cite{rietkerk2021evasion}.  

From the mathematical point of view, tipping is most commonly identified with a stable configuration in a governing model -- the previously mentioned favorable state -- approaching a saddle-node (or fold) bifurcation due to one or more slowly varying parameters. When these parameters pass through their critical saddle-node value, the stable state vanishes, and the model drives the system towards another stable state -- the less desirable state. (In the literature, various types of tipping are distinguished; we focus in this paper on {\it bifurcational tipping} \cite{ashwin2012tipping}.) This picture is at the same time appealing and overly simplified. Predictions of tipping in ecosystems are typically based on non-spatial (ODE) models \cite{rietkerk2021evasion}. However, recently combined ecological-mathematical research has shown that in the (conceptual) Klausmeier/Gray-Scott model for dryland ecosystems, which incorporates spatial effects \cite{klausmeier1999regular,van2013rise}, a Turing bifurcation precedes the fold bifurcation \cite{siteur2014beyond}. In this setting, the Turing bifurcation initiates the formation of a family of vegetation patterns that introduces the multi-stability, which provides the ecosystem with alternative pathways of response to the changing environmental conditions. Rather than undergoing spatially homogeneous collapse, the spatial system may instead transition through a variety of distinct, stable patterned states, thereby adapting more gradually to changing external circumstances \cite{bastiaansen2020effect}. These states are all embedded in the family of stable solutions known as the {\it Busse balloon} \cite{busse1978non}, which can be viewed as the global continuation of the family of stable states that appear locally at a Turing bifurcation \cite{doelman2019pattern}. In other words, a Busse balloon defines a region of stable spatially periodic patterns in $(\text{parameter, wavenumber})$-space through which an ecosystem -- or any spatially extended complex system -- may circumvent critical collapse. Observations of vegetation patterns in Somalia indicate that such Busse balloons indeed exist in dryland ecosystems \cite{bastiaansen2018multistability}. Moreover, Turing bifurcations are widespread in spatial ecological (and biological) models \cite{meron2015nonlinear,murray2007mathematicalII}, and Busse balloons have recently been shown to play a central role in coral reef and seagrass ecosystem models \cite{detmer2025busse, moreno2025spatiotemporal}. Thus, the mathematical question whether (or not!) a Turing bifurcation can increase the resilience of an ecosystem by enabling it to follow a pathway around a tipping point, thereby evading collapse, is of central importance to our understanding of ecosystem dynamics driven by a changing climate.    

\begin{figure}[t]
    \centering
    \begin{subfigure}[t]{0.32\textwidth}
        \centering
        \includegraphics[width=\textwidth]{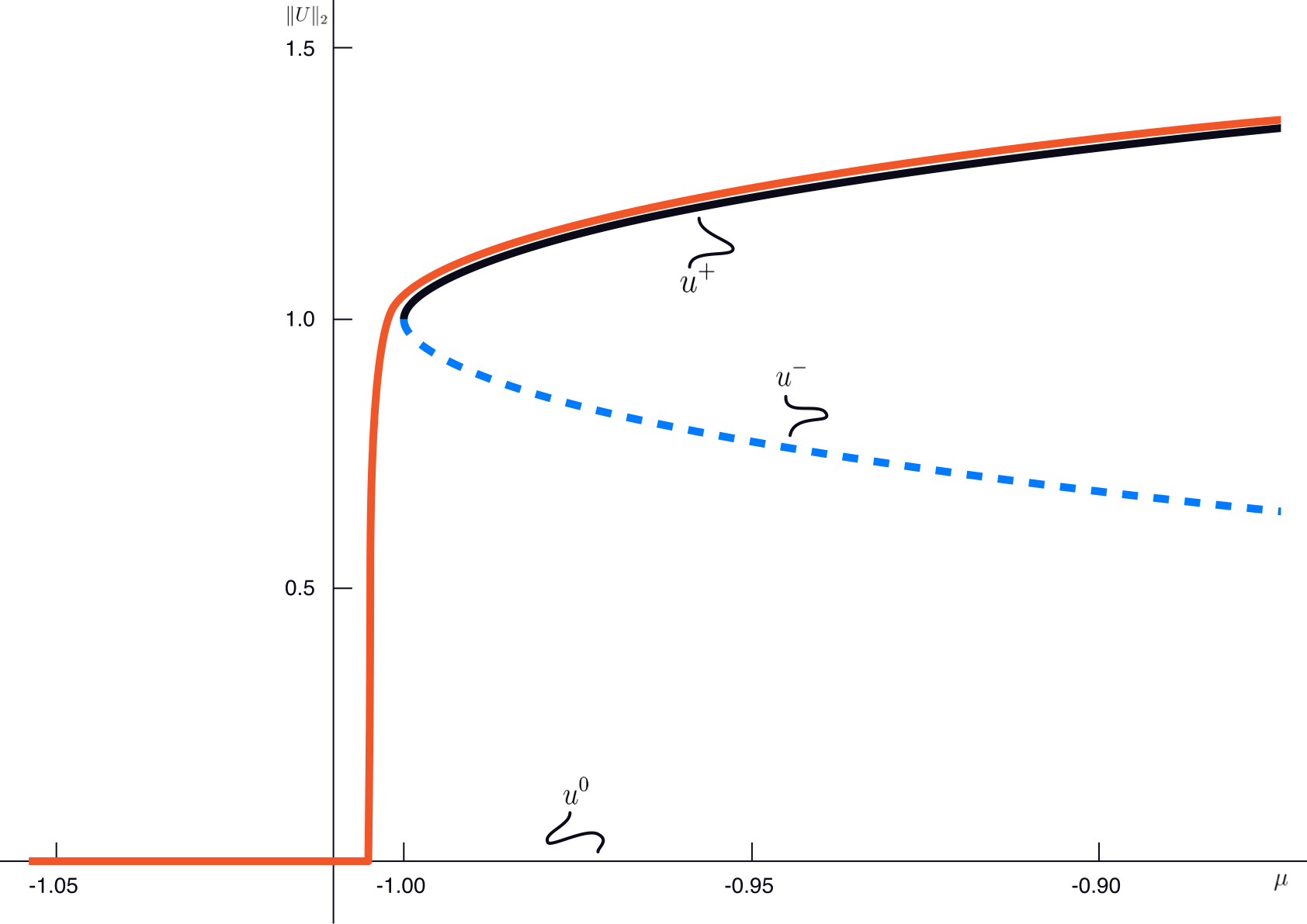}
        \caption{\underline{Spatially homogeneous tipping}: as $\mu(t)$ decreases through the saddle-node bifurcation at $\mu^* =-1$, the system collapses to the trivial state $U(x,t) \equiv 0$.}
        \label{fig:Intro-ODE}
    \end{subfigure}
    \hfill
    \begin{subfigure}[t]{0.32\textwidth}
        \centering
        \includegraphics[width=\textwidth]{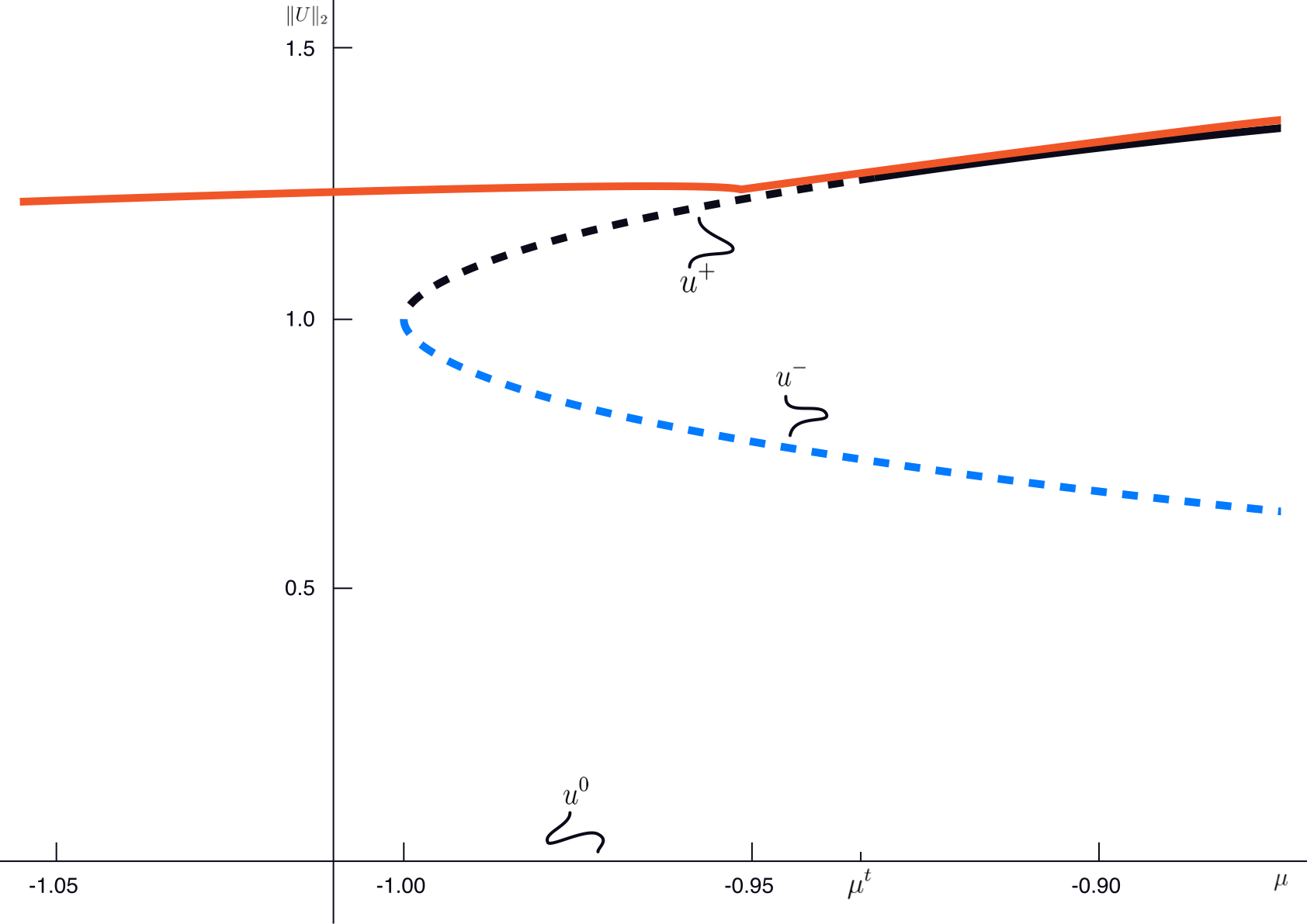}
        \caption{$\eta=2$. The stable patterns emerging from the Turing bifurcation persist as $\mu(t)$ decreases through the tipping point: \underline{the evasion of tipping}.}
        \label{fig:Intro-Turing}
    \end{subfigure}
    \hfill
    \begin{subfigure}[t]{0.32\textwidth}
        \centering
        \includegraphics[width=\textwidth]{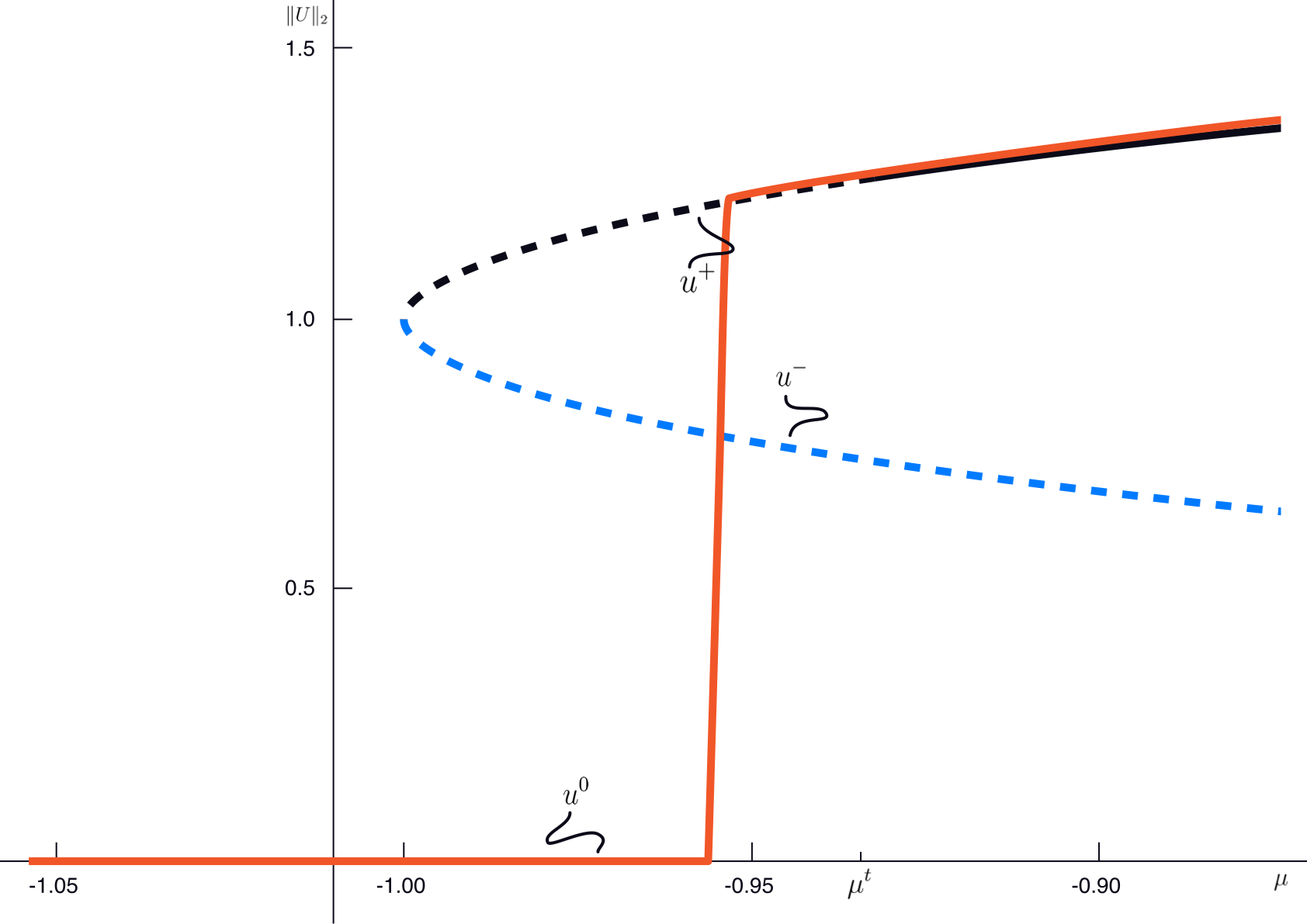}
        \caption{$\eta=-1$. The Turing bifurcation initiates the collapse to the trivial state $U(x,t) \equiv 0$ before $\mu(t)$ reaches the tipping point: \underline{Turing triggers tipping}.}
        \label{fig:Intro-TTT}
    \end{subfigure}
    \caption{The (averaged) $L^2$-norm $\|U\|_2/L$ of solutions $U(x,t)$ of the scalar sixth-order model \eqref{eq:scalar example} as function of the slowly decreasing parameter $\mu(t)$, obtained by simulations on the interval $x \in (0,L)$ with $L = 4 \pi$ (and periodic boundary conditions) for $t \in (0,T)$ with $T = 2000$ starting near the stable homogeneous state $u^+(\mu_0)$; $\mu(t) = \mu_0 - \epsilon t$ with $\mu_0 \approx -0.874 > \mu_t(\nu) \approx 0.929$ , with $\epsilon = \delta^2/4000$, $\nu = 1 - \delta$ and $\delta = 0.6$ (see \Cref{subsubsec:basic analysis:a model} and Remark \ref{rem:backgroundFig1}).}
    \label{fig:Intro-tipping}
\end{figure}

In this paper, we study this question under the assumption that the underlying (spatially extended, complex) system exhibits both a fold and a Turing bifurcation for parameter combinations that are {\it sufficiently close} to each other. Specifically, we consider the nature of pattern formation near a co-dimension 2 point at which a saddle-node and Turing bifurcation coincide. Our analysis is carried out in two settings: systems of $n$-component reaction-diffusion equations,
\begin{equation}
\label{eq:reaction-diffusion system:rds-Intro}
    \partial_t U = F(U;\mu,\nu) + D \partial_x^2 U,
\end{equation} where $U(x,t): \mathbb{R}\times \mathbb{R}^+ \to \mathbb{R}^n$, $D$ a diagonal diffusion matrix, $F:\mathbb{R}^{n+2} \to \mathbb{R}^n$, and parameters $\mu, \nu \in \mathbb{R}$ (see \Cref{sec:reaction diffusion system}); and higher-order scalar PDEs, 
\begin{equation}
\label{eq:illustration pde:rtab-Intro}
    \partial_t U = F(U;\mu) + \Psi(\partial_x^2 U, \dots, \partial_x^{2m} U, U;\nu),
\end{equation}
with $U(x,t): \mathbb{R}\times \mathbb{R}^+ \to \mathbb{R}$, $F: \mathbb{R}^2 \to \mathbb{R}$, $\Psi:\mathbb{R}^{m+2} \to \mathbb{R}$, and parameters $\mu, \nu \in \mathbb{R}$ (see \Cref{sec:scalar}). Since the co-dimension 2 Turing-fold bifurcation degenerates for $n=2$ in \eqref{eq:reaction-diffusion system:rds-Intro} and $m=2$ in \eqref{eq:illustration pde:rtab-Intro} (see \Cref{rem:k*=0}), we assume throughout that $n,m > 2$. Our focus on pattern formation near a Turing-fold bifurcation in \eqref{eq:reaction-diffusion system:rds-Intro} is natural, since ecosystem models are typically given in terms of systems of reaction-diffusion equations (see \cite{meron2015nonlinear, murray2007mathematicalII, rietkerk2021evasion} and the references therein). However, for transparency of presentation, we begin our investigations in the context of higher-order scalar phase-field models of the form \eqref{eq:illustration pde:rtab-Intro} \cite{caginalp1986higher}, using the sixth-order equation
\begin{equation}
\label{eq:scalar example}
    \partial_t U = \mu U + 2U^2 - U^3 + \nu \partial_x^2 U + 2 \partial_x^4 U + \partial_x^6 U + \eta (\partial_x^2 U)^2
\end{equation}
as the main motivational example. While this equation is of the form of the sixth-order phase-field models studied in \cite{gardner1990traveling} (for $\eta = 0$), we use it purely as the simplest model for which our approach can be demonstrated, without reference to a physical interpretation. 

In \Cref{fig:Intro-tipping}, we present three simulations of model \eqref{eq:scalar example} in which one of the main parameters $\mu$ decreases slowly in time (indicating {\it worsening environmental conditions}). \Cref{fig:Intro-ODE} shows the typical critical transition associated with tipping in a spatially homogeneous ODE: as $\mu(t)$ decreases through saddle-node bifurcation at $\mu^* =-1$ (\Cref{subsec:a model}), the system collapses on a fast time scale to the trivial state $U(x,t) \equiv 0$. By contrast, in \Cref{fig:Intro-Turing,fig:Intro-TTT} spatial variation is allowed, and we observe two possible scenarios. In \Cref{fig:Intro-Turing}, the system undergoes a Turing bifurcation as $\mu$ decreases through $\mu_t > \mu^*$, leading to the formation of stable spatial patterns that persist beyond the tipping point. In this way, the dynamics effectively evades the spatially homogeneous collapse observed in \Cref{fig:Intro-ODE} through pattern formation. However, \Cref{fig:Intro-TTT} shows that a change in the parameter $\eta$ -- which does not affect the locations (in parameter space) of the saddle-node or Turing bifurcation -- has the opposite effect: the Turing bifurcation at $\mu_t > \mu^*$ triggers collapse before the tipping point is reached -- a scenario more in line with the original idea of pattern formation as early warning signal for tipping (see however \Cref{rem:backgroundFig1}).

We begin our analysis by assuming the existence of a critical combination $(\mu^*,\nu^*)$ of parameters $(\mu,\nu)$ in \eqref{eq:reaction-diffusion system:rds-Intro}/\eqref{eq:illustration pde:rtab-Intro} at which a background state $U(x,t) \equiv u^*$ of \eqref{eq:reaction-diffusion system:rds-Intro}/\eqref{eq:illustration pde:rtab-Intro} simultaneously undergoes a saddle-node and a Turing bifurcation. Additionally, we introduce a small parameter $0 < \delta \ll 1$ that measures the distance of $(\mu,\nu)$ from the co-dimension 2 point $(\mu^*,\nu^*)$. Based on the classical Ginzburg-Landau approach near a `pure' Turing destabilization \cite{schneider2017nonlinear} -- which, by assumption, is not of Turing-Hopf type \cite{rademacher2007instabilities, scheel2003radially} -- together with a set of non-degeneracy conditions (see \Cref{subsubsec:non-degeneracy conditions:gs,subsec:non-degeneracy conditions:rd}), we introduce the following Ansatz that models the dynamics of small amplitude solutions $U(x,t)$ of \eqref{eq:reaction-diffusion system:rds-Intro}/\eqref{eq:illustration pde:rtab-Intro} near the co-dimension 2 point:
\begin{equation}
\label{eq:Ansatz:abd:rd-Intro}
U_{AB}(x, t) = u^* + \delta \left[\left(A(\xi,\tau)  v^{t}_* e^{i k^* x}  + \; {\rm c.c.}\right) + \left(\tilde{u}^s + B(\xi,\tau) v^{s}_* \right) \right] \; +  \mathcal{O}(\delta^{3/2}),
\end{equation}
where $A(\xi,\tau): \mathbb{R} \times \mathbb{R}^+ \to \mathbb{C}$ and $B(\xi,\tau): \mathbb{R} \times \mathbb{R}^+ \to \mathbb{R}$ are unknown functions corresponding to the Turing and saddle-node bifurcations, respectively; time and space are rescaled as $\tau = \delta t$ and $\xi = \sqrt{\delta} x$; $k^* > 0$ (Remark \ref{rem:k*=0}) is the critical wavenumber associated to the Turing bifurcation; the vectors $v^t_* \in \mathbb{R}^n$ and $v^s_* \in \mathbb{R}^n$, respectively, span the (one-dimensional) kernels of the $n \times n$ matrix
\[
T(k; \mu^*, \nu^*) = F_u(u^*;\mu^*,\nu^*) - k^2 D,
\]
at $k=k^*$ and $k=0$, respectively, where $F_u(U;\mu,\nu)$ is the Jacobian of the vectorfield $F(U;\mu,\nu)$; and $\tilde{u}^s = (u^s - u^*)/\delta + \mathcal{O}(\delta)$, with $u^s \coloneqq u^s(\mu,\nu)$ denoting the fold point value of $U$ away from $(\mu^*,\nu^*)$ (i.e., $u^* = u^s(\mu^*,\nu^*)$). In the scalar models considered in \Cref{sec:scalar}, we have $v^t_* = v^s_* = 1$ and $\tilde{u}^s = 0$. Note that Ansatz \eqref{eq:Ansatz:abd:rd-Intro} differs in a few essential aspects from those used in the standard Ginzburg-Landau approach. In \Cref{subsec:a model}, we present an extensive motivation of the structure of this Ansatz in the context of the example model \eqref{eq:scalar example}, by considering the standard Ginzburg-Landau approach near the Turing bifurcation and showing how it degenerates as the Turing point approaches the fold point. 

By substituting an extended version of Ansatz \eqref{eq:Ansatz:abd:rd-Intro} -- both as asymptotic series in $\sqrt{\delta}$ and as Fourier series in $e^{i k^* x}$ -- into the governing equations \eqref{eq:reaction-diffusion system:rds-Intro}/\eqref{eq:illustration pde:rtab-Intro}, and following the standard procedure for deriving modulation equations (see \cite{doelman2019pattern} and the references therein), we obtain a coupled system of modulation equations for the (rescaled) amplitudes $(A(\xi,\tau), B(\xi,\tau)): \mathbb{R} \times \mathbb{R}^+ \to \mathbb{C} \times \mathbb{R}$ that reads, in `canonical' form, 
\begin{equation}
\label{eq:canonical form:csd:gs}
\begin{cases}
     \phantom{\frac{1}{\alpha}}A_\tau = A_{\xi\xi} + A - AB \\
     \frac{1}{\alpha}B_\tau = d B_{\xi\xi} + 1 - R - B^2 + \beta|A|^2.
\end{cases}
\end{equation}
Here, the parameter $R \in \mathbb{R}$ unfolds the Turing-fold bifurcation and thus has taken over the role of original parameters $(\mu,\nu)$: in \eqref{eq:canonical form:csd:gs}, the Turing bifurcation takes place at $R=0$, and the fold bifurcation follows as $R$ increases through $1$. One of the key results of the present work is that it provides explicit (relatively concise) expressions for the coefficients $\alpha >0$, $d > 0$, and $\beta \in \mathbb{R}$ in terms of underlying systems \eqref{eq:reaction-diffusion system:rds-Intro} and \eqref{eq:illustration pde:rtab-Intro}; see \Cref{sec:reaction diffusion system,sec:scalar}. In particular, our methods yield explicit information on the sign of the $\beta$ coefficient in the modulation system \eqref{eq:canonical form:csd:gs} for any model of type  \eqref{eq:reaction-diffusion system:rds-Intro}/\eqref{eq:illustration pde:rtab-Intro}. This coefficient plays a role analogous to the Landau coefficient in the standard Ginzburg-Landau approximation that governs the Turing bifurcation (see \Cref{thm:close enough:t:cs,thm:supercrit:t:cs}).

Naturally, the next question to be answered is: {\it Under what conditions  does AB-system \eqref{eq:canonical form:csd:gs} admit spatially periodic stationary {\it plane waves}, i.e., solutions of the form}
\begin{equation}
\label{eq:periodic solutions:cs-Intro}
A(\xi,\tau) = \bar{A} e^{iK\xi} =  \bar{A}(K,R) e^{iK\xi},  \; \;  B(\xi,\tau) = \bar{B} = \bar{B}(K,R)
\end{equation}
{\it with $K \in \mathbb{R}$ -- that appear at the Turing bifurcation ($R=0$) and persist as stable patterns beyond $R=1$?} In other words: {\it Under which conditions can the tipping point at $R=1$ be evaded through the formation of spatial patterns?} (as in \Cref{fig:Intro-Turing}). The answer follows from the existence and stability analysis presented in \Cref{sec:coupled system}:  stable plane wave patterns exist \eqref{eq:canonical form:csd:gs} for $R > 0$ if $\beta > 0$, whereas no such stable solutions exist when $\beta < 0$. 

\begin{claim}
\label{thm:close enough:t:cs}
Consider \eqref{eq:reaction-diffusion system:rds-Intro}/\eqref{eq:illustration pde:rtab-Intro} and assume that all the non-degeneracy conditions are satisfied (cf. Sections \ref{subsubsec:non-degeneracy conditions:gs} and \ref{subsec:non-degeneracy conditions:rd}). Then, for sufficiently small $0 < \delta \ll 1$ -- i.e., when $(\mu,\nu)$ is sufficiently close to the co-dimension 2 point $(\mu^*,\nu^*)$ -- the following holds. \underline{If $\beta > 0$} in the associated AB-system \eqref{eq:canonical form:csd:gs}, the Turing-fold bifurcation is supercritical and generates stable (small amplitude) spatially periodic patterns in \eqref{eq:reaction-diffusion system:rds-Intro}/\eqref{eq:illustration pde:rtab-Intro} that persist beyond the fold point (and correspond through \eqref{eq:Ansatz:abd:rd-Intro} to plane waves solutions \eqref{eq:periodic solutions:cs-Intro} of \eqref{eq:canonical form:csd:gs}). Conversely, if $\beta < 0$, the Turing-fold bifurcation is subcritical and there are no stable plane waves near the Turing-fold point. Thus, tipping is evaded by the formation of (small amplitude) patterns \underline{only for $\beta > 0$}.
\end{claim}

In this paper, we do not pursue a rigorous validity proof of this claim (see \Cref{rem:validity}), the main goals of the present work are to develop a systematic and asymptotically accurate derivation of coupled AB-system of modulation equations \ref{eq:canonical form:csd:gs}, to study its most simple plane wave patterns \eqref{eq:periodic solutions:cs-Intro} and to test its validity by various kinds of numerical simulations (see \Cref{sec:simulations}). \Cref{thm:close enough:t:cs} is formulated based on the insights obtained by this approach. 

In the context of the AB-system derivation procedure, we also examine in detail how the classical Ginzburg-Landau approximation is embedded in the $AB$-system near the Turing point. Specifically, we `zoom in' on Turing bifurcation by taking our parameters $\mathcal{O}(\varepsilon^2)$ close to a Turing point, with $0 < \varepsilon \ll \delta \ll 1$, so that the fold point remains `far away' from the perspective of the associated Ginzburg-Landau approach (see \Cref{subsubsec:ginzburg-landau derivation:a model,subsubsec:ginzburg-landau derivation:gs,subsubsec:Reducing the AB-system:gs:se} for the scalar model and \Cref {subsec:ginzburg-landau derivation:rd,subsec:reducing the AB-system back to the ginzburg landau system:rd} for system \eqref{eq:reaction-diffusion system:rds-Intro}). Through this, we uncover the relation between the AB-system and its embedded Ginzburg–Landau equation, as stated in  \Cref{thm:supercrit:t:cs}.

\begin{claim}
\label{thm:supercrit:t:cs}
Consider \eqref{eq:reaction-diffusion system:rds-Intro}/\eqref{eq:illustration pde:rtab-Intro}, and assume that all non-degeneracy conditions are satisfied and that $(\mu,\nu)$ is sufficiently close to $(\mu^*,\nu^*)$. Then the Turing-fold bifurcation is supercritical (i.e., $\beta > 0$), respectively subcritical ($\beta < 0$), if and only if the Turing bifurcation in the `zoomed in' Ginzburg-Landau equation is likewise super- or subcritical.
\end{claim}

It is remarkable that such an explicit relation can be shown, as determining the nature of a Turing bifurcation is typically an overwhelming analytic task. Even for two-component reaction–diffusion systems, there are only a few models for which this analysis has been carried out explicitly (see \cite{doelman2019pattern, kuramoto1984chemical} and the references therein). However, in \Cref{sec:scalar,sec:reaction diffusion system}, we show that the Landau coefficient $L$ of the `zoomed in' Ginzburg-Landau equation -- which determines the nature of the Turing bifurcation -- is of the form $L = L^*/\delta + \mathcal{O}(1)$, where $0 < \delta \ll 1$ measures the distance to the fold point. Thus, the sign of $L^*$ determines whether the Turing bifurcation is sub- or supercritical. Our method enables us to derive explicit expressions for $L^*$ -- even in the fully general setting of $n$-component reaction-diffusion systems \eqref{eq:reaction-diffusion system:rds-Intro} -- and to show that $L^*$ can be written as $L^* = -(\text{a positive number}) \times \beta$, from which \Cref{thm:supercrit:t:cs} follows. Naturally, a proof of this claim must be based on the same rigorous underpinning as that of \Cref{thm:close enough:t:cs}.

\begin{figure}[t]
    \centering
    \begin{subfigure}[t]{0.35\textwidth}
        \centering
        \includegraphics[width=\textwidth]{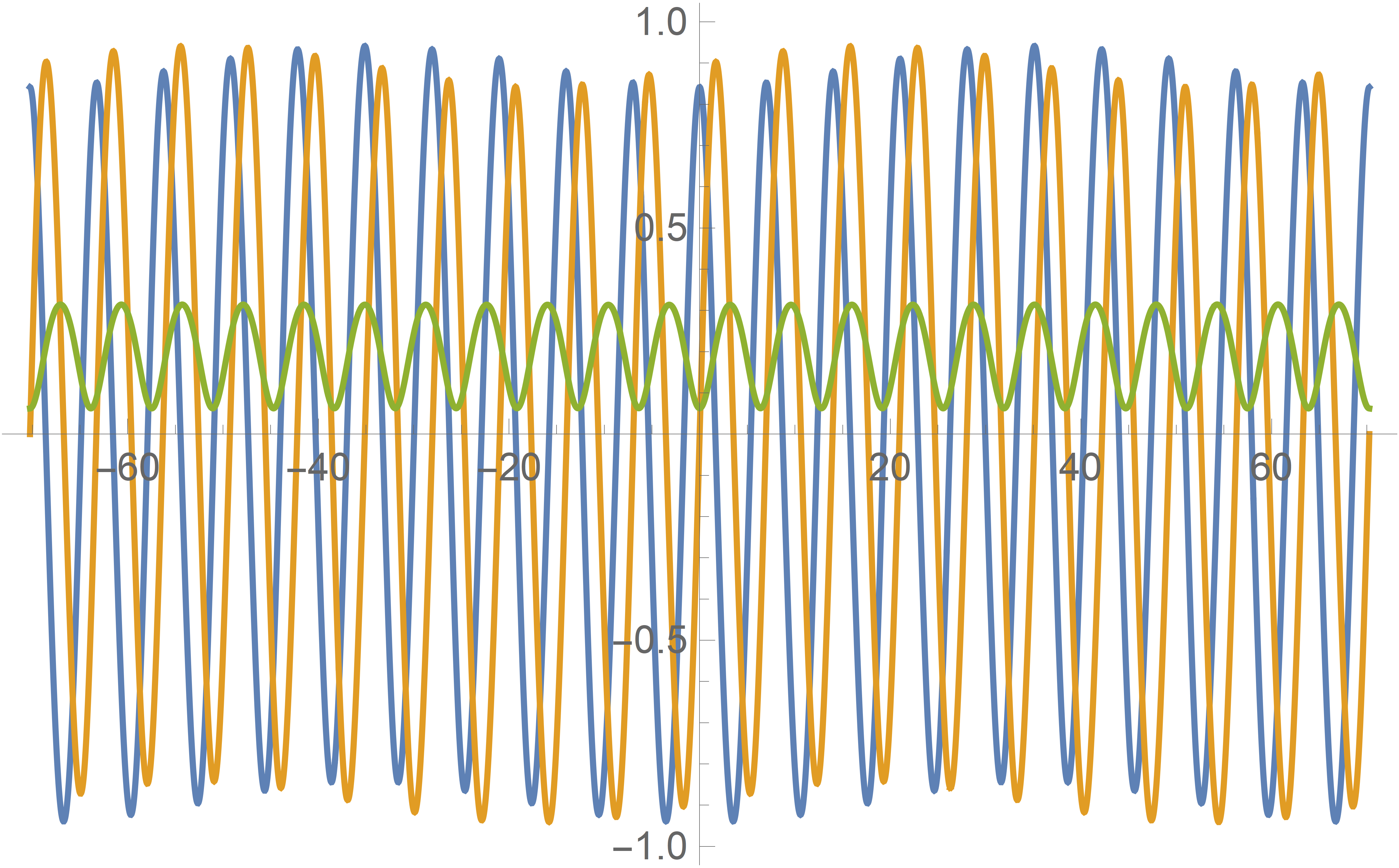}
        \caption{Re$\,(A(\xi,\tau))$ (blue), Im$\,A((\xi,\tau))$ (orange), $B(\xi,\tau)$ (green) at $\tau=2000$.}
        \label{fig:Intro-qp-a}
    \end{subfigure}
    \hfill
    \begin{subfigure}[t]{0.35\textwidth}
        \centering
        \includegraphics[width=\textwidth]{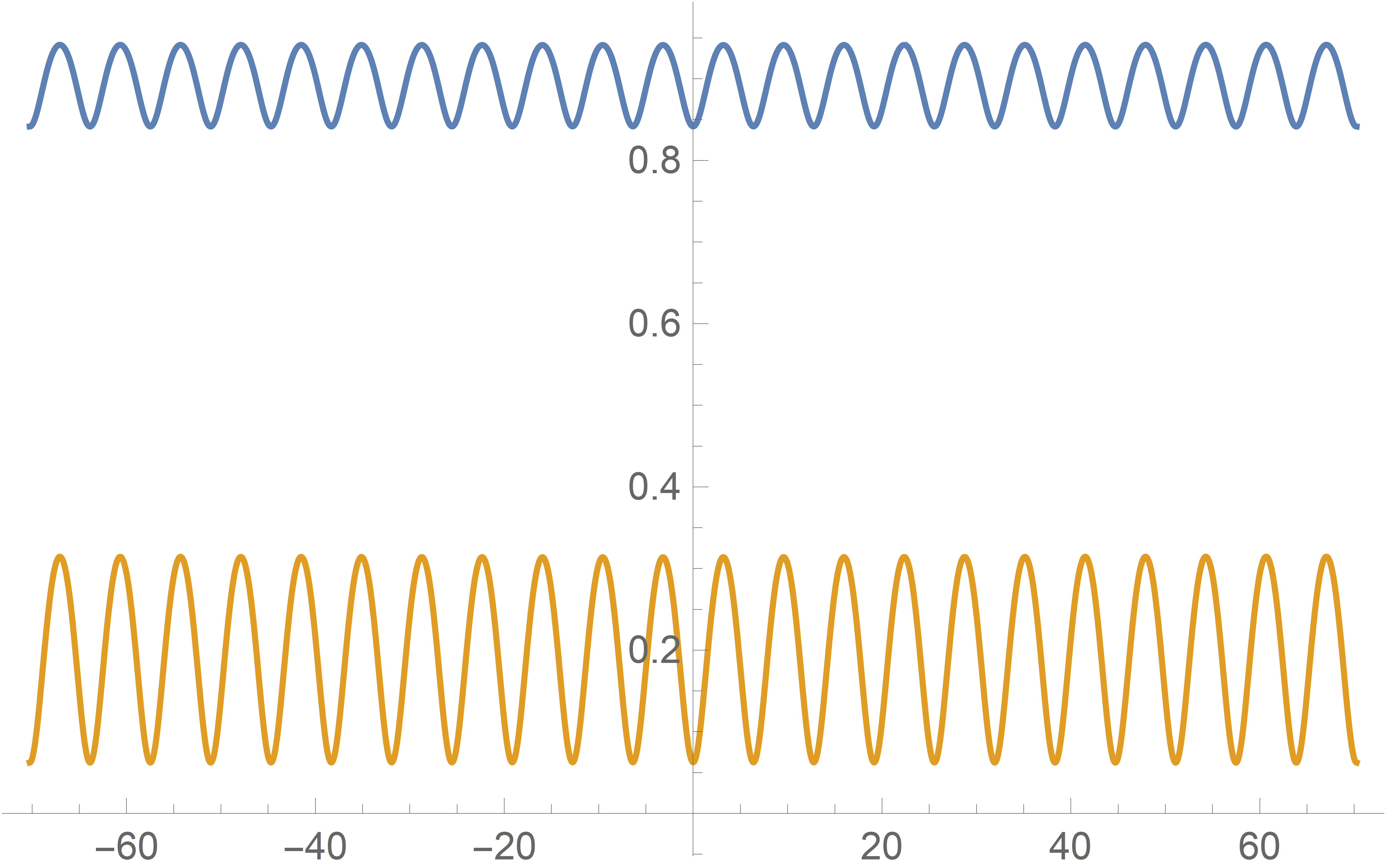}
        \caption{$|A(\xi,\tau)|$ (blue) and $B(\xi,\tau)$ (orange) at $\tau=2000$.}
        \label{fig:Intro-qp-b}
    \end{subfigure}
    \hfill
    \begin{subfigure}[t]{0.28\textwidth}
        \centering
        \includegraphics[width=\textwidth]{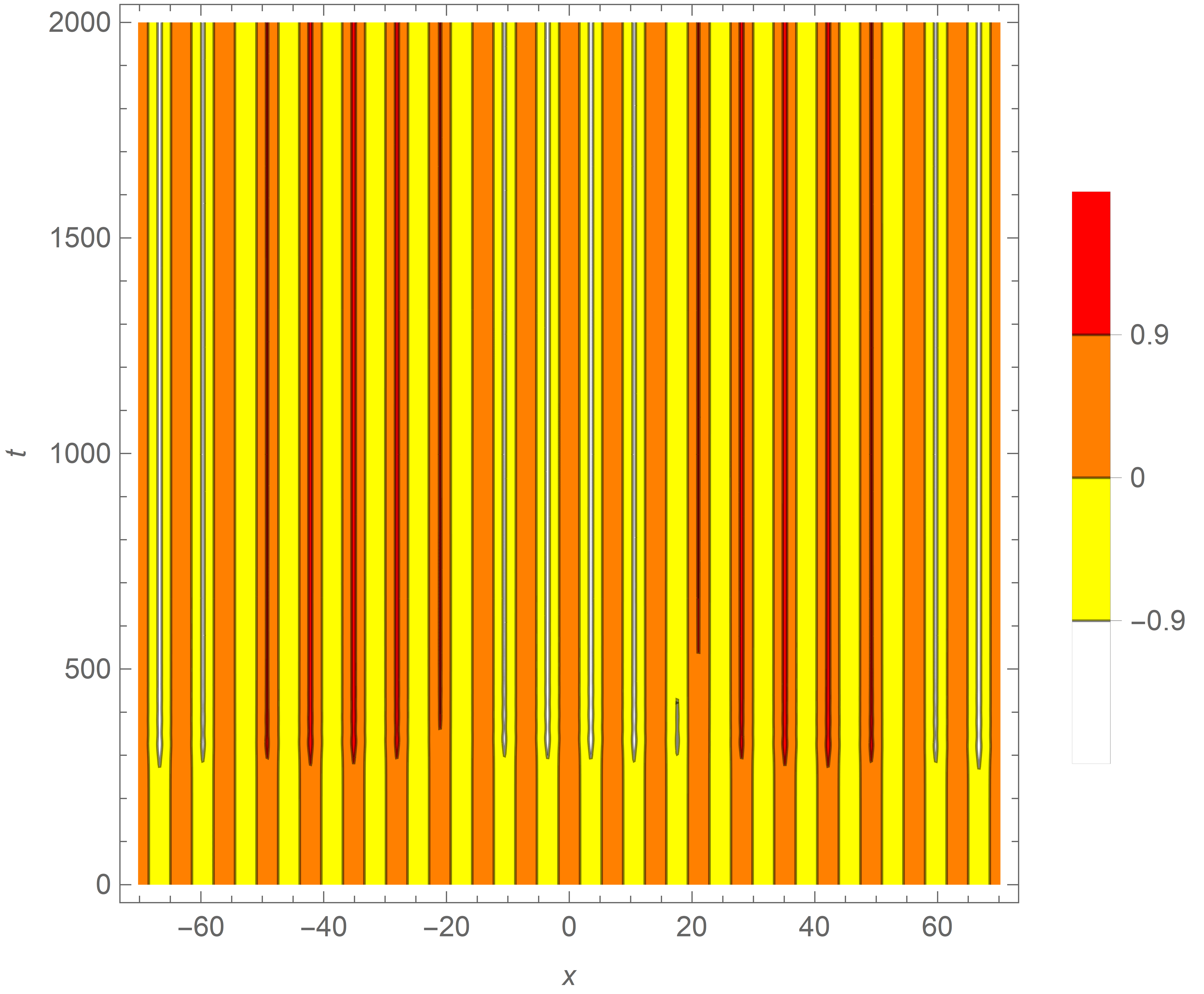}
        \caption{contour plot of Re$\,(A(\xi,\tau))$ for $\tau \in (0,2000)$.}
        \label{fig:Intro-qp-c}
    \end{subfigure}
    \caption{A numerically stable (stationary) quasi-periodic pattern in AB-system \eqref{eq:canonical form:csd:gs} with $\alpha = 0.8$, $d=1/3$, $\beta =1$. The simulation is done on the interval $\xi \in (-20 \pi \sqrt{4/5}, 20 \pi \sqrt{4/5})$) with homogeneous Neumann boundary conditions and the plane wave \eqref{eq:periodic solutions:cs-Intro} with $K = \sqrt{4/5}$ as initial condition. Parameter $R$ is chosen such that this wave is weakly unstable with respect to a Turing bifurcation, i.e., close to the bifurcation value: $R = R_t(\sqrt{4/5}) - 0.01 = 53/30-0.01$ (see \Cref{sec:coupled system}). }
    \label{fig:Intro-qp}
\end{figure}

As with the classical (co-dimension 1) Ginzburg-Landau equation, it is natural to expect that stable patterns -- or, in more general terms, attractors -- of AB-system \eqref{eq:canonical form:csd:gs} have counterparts. In other words, likely there are leading order approximations of similar patterns and attractors in the underlying system \cite{aranson2002world,busse1978non,cross1993pattern,mielke2002ginzburg,mielke1995attractors,schneider1998nonlinear}. In \Cref{subsec:corroborating validity:s}, this expectation is numerically confirmed for the plane waves \eqref{eq:periodic solutions:cs-Intro}. Although the attractors exhibited by the real Ginzburg-Landau equation -- i.e., the Ginzburg-Landau equation with real coefficients that governs Turing bifurcations -- do not significantly go beyond (modulations of) plane waves \cite{aranson2002world,cross1993pattern}, it is natural to investigate whether the AB-system may generate a richer kind of dynamics. This indeed is the case, as we find by numerical simulations. 

In \Cref{fig:Intro-qp} we show numerically that AB-system \eqref{eq:canonical form:csd:gs} exhibits stable, stationary patterns that are quasi-periodic in $\xi$. In \Cref{fig:Intro-qp-a}, Re$\,(A(\xi,\tau))$, Im$\,A((\xi,\tau))$, and $B(\xi,\tau)$ are plotted as function of $\xi$ (for $\tau = 2000$), while \Cref{fig:Intro-qp-b} indicates that both $|A(\xi,2000)|$ and $B(\xi,2000)$ are periodic functions of $\xi$. The contour plot of \Cref{fig:Intro-qp-c} shows the evolution towards the stationary quasi-periodic pattern. It should be noted that the real Ginzburg-Landau equation also has a family of quasi-periodic patterns of this type -- i.e., solutions for which $|A|$ varies periodically in $\xi$ -- but that these patterns are unstable \cite{doelman1995instability}. In fact, it is to be expected that these patterns can also be recovered in the AB-system by considering $R$ asymptotically close to $0$, and that these patterns will inherit the instability of the Ginzburg-Landau limit. However, the patterns shown in \Cref{fig:Intro-qp} are not of this type: they appear from parameter settings (and initial conditions) close to a Turing instability of the plane waves (see \Cref{sec:coupled system,sec:simulations,sec:discussion}).          

\begin{figure}[t]
    \centering
    \begin{subfigure}[t]{0.4\textwidth}
       \centering
        \includegraphics[width=\textwidth]{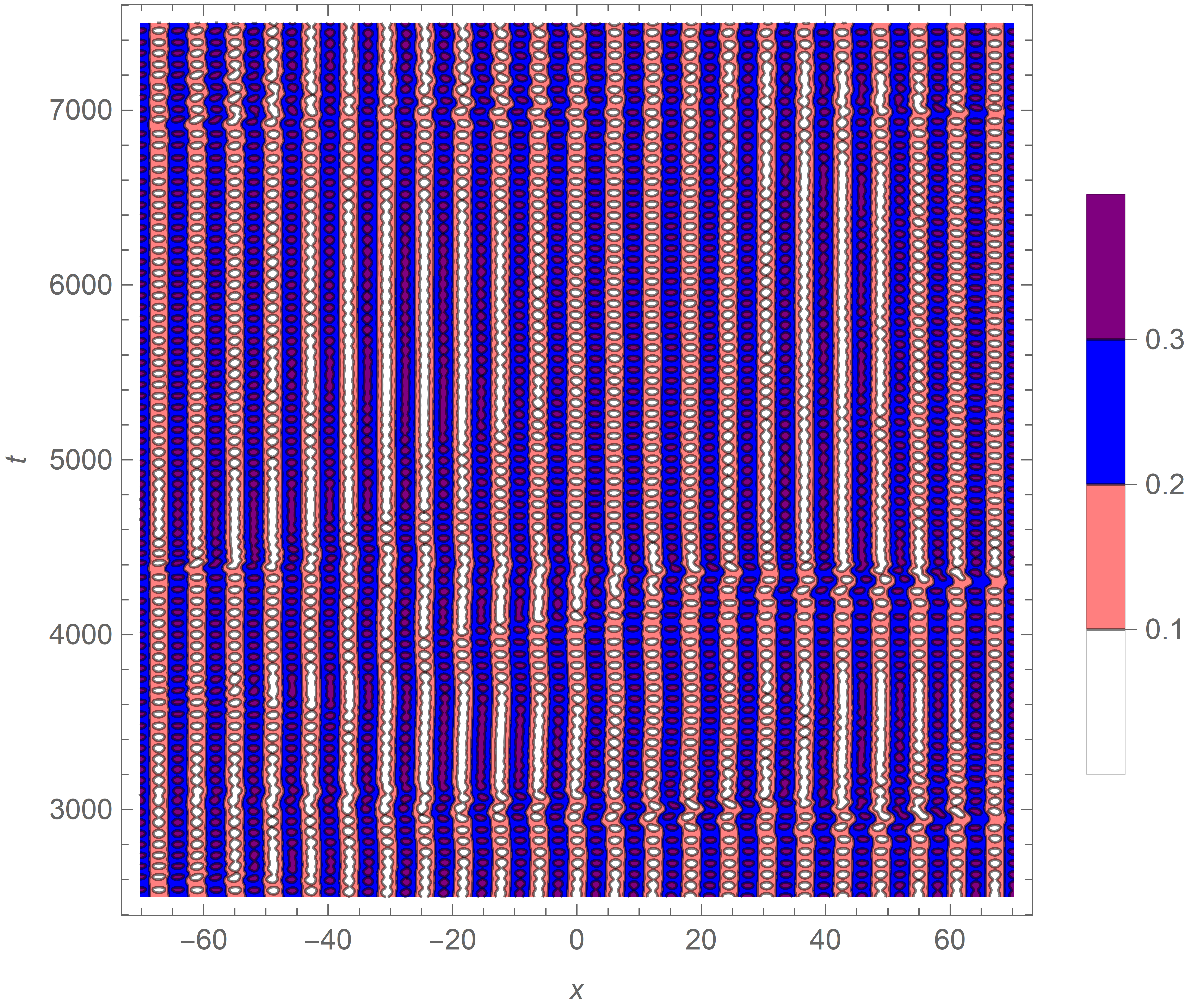}
        \caption{contour plot of $B(\xi,\tau)$ for $\tau \in (2500,7500)$}
        \label{fig:Intro-chaos-c}
    \end{subfigure}
    \hfill
    \begin{subfigure}[t]{0.58\textwidth}
        \centering
        \includegraphics[width=\textwidth]{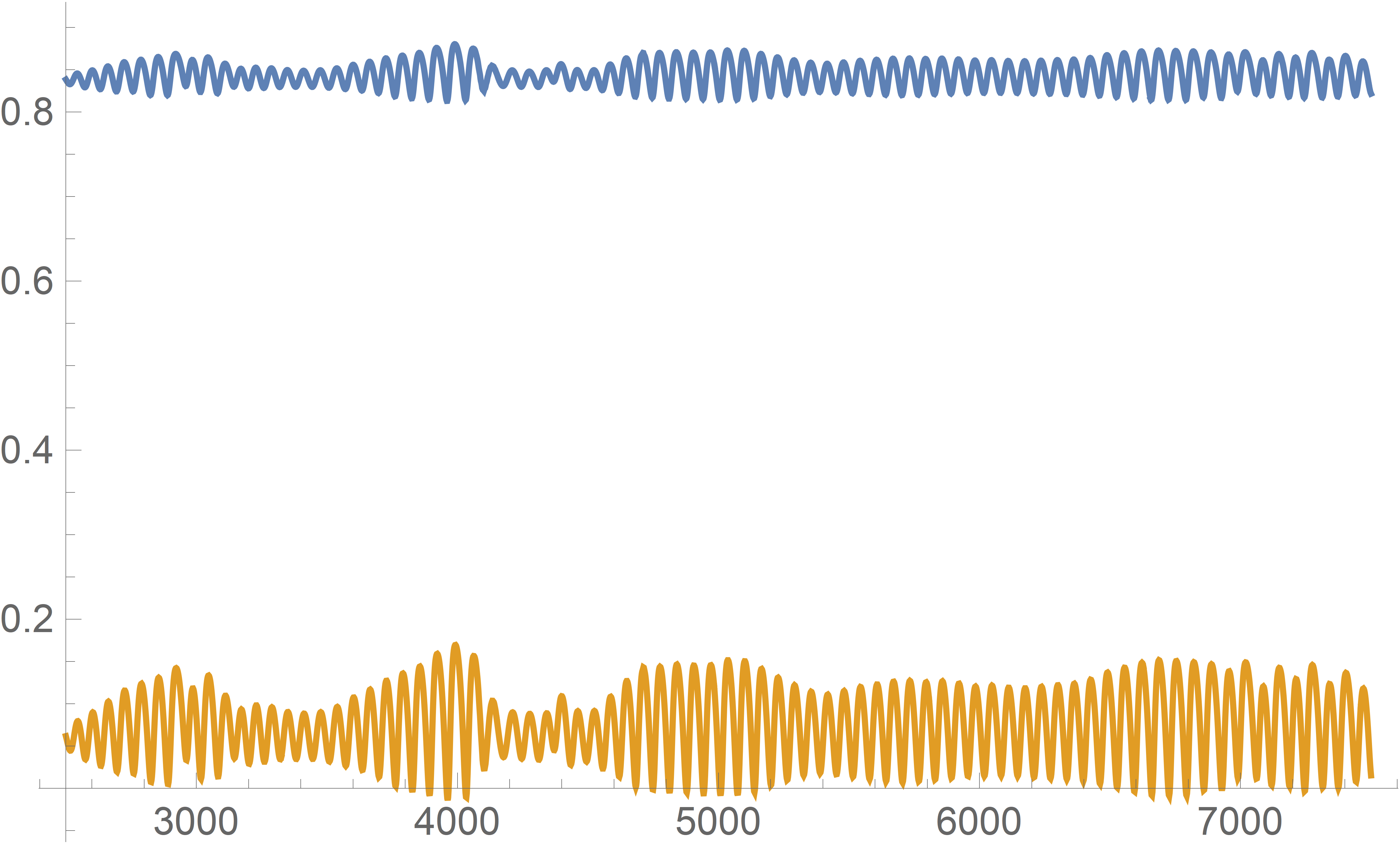}
        \caption{$|A(0,\tau)|$ (blue), $B(0,\tau)$ (orange) for $\tau \in (2500,7500)$}
        \label{fig:Intro-chaos-b}
    \end{subfigure}
    \hfill
    \caption{Irregular patterns suggesting spatio-temporal chaotic dynamics: a simulation of AB-system \eqref{eq:canonical form:csd:gs} on the interval $\xi \in (-20 \sqrt{4/5} \pi, -20 \sqrt{4/5} \pi)$ for $t \in (0,10000)$ with homogeneous Neumann boundary conditions. Initial conditions and parameters are chosen as in \Cref{fig:Intro-qp}, except for $\alpha$: $\alpha = 0.75715$.}
    \label{fig:Intro-chaos}
\end{figure}

The Turing instability of the plane waves -- i.e., a destabilization initiated by perturbation with a {\it most critical wave number} $k_c > 0$ -- does not occur in the real Ginzburg-Landau equation (though it does arise in the complex Ginzburg-Landau equation \cite{matkowsky1993stability}). This suggests that the presence of such an instability may give rise to richer dynamics than those exhibited by the real Ginzburg–Landau equation. This is indeed the case: in \Cref{fig:Intro-chaos-c}, we present a contour plot of $B(\xi,\tau)$ that displays spatio-temporal chaos-like behavior. \Cref{fig:Intro-chaos-b} shows the corresponding plots of $|A(\xi,\tau)|$ and $B(\xi,\tau)$ at the fixed position $\xi = 0$ over the interval $\tau \in (2500, 7500)$. The times series $|A(0,\tau)|$ and $B(0,\tau)$ indeed seems to lack regularity. 

In section \ref{subsec:interesting behavior in the AB-system:s}, we show more examples of complex dynamics exhibited by the AB-system. Moreover, we show (again numerically) that the same type of (complex) dynamics is inherited by the underlying system. This seems to be the case for all types of dynamics we have encountered in our (not overly extensive) numerical simulations. We, for instance, uncover the same kind of chaotic behavior as shown \Cref{fig:Intro-chaos} in a sixth-order scalar equation that is very similar to the example model \eqref{eq:scalar example}. Thus, the AB-system \eqref{eq:canonical form:csd:gs} we derived in the present work seems to open up a promising `tool' by which the rich and, a priori, unexpected dynamics associated with a co-dimension 2 Turing-fold bifurcation may be appreciated and understood.  

The set-up of this article is as follows. In \Cref{subsec:a model}, we examine the Turing-fold bifurcation in the simplest setting of the example model problem \eqref{eq:scalar example} and motivate in detail the approach we developed for the derivation of the coupled AB-system of modulation equations. Next, we consider a more general class of scalar models in \Cref{subsec:general setting:se} and show that the results of \Cref{subsec:a model} are generic given some natural non-degeneracy conditions. In \Cref{sec:reaction diffusion system}, we develop the same method -- and obtain the same AB-system \eqref{eq:canonical form:csd:gs} -- in the setting of the general $n$-component system of reaction-diffusion equations \eqref{eq:reaction-diffusion system:rds-Intro} (near a Turing-fold bifurcation). In \Cref{sec:coupled system}, we analyze the existence and stability of the plane wave solutions \eqref{eq:periodic solutions:cs-Intro} of the AB-system. We (numerically) corroborate the validity of the AB-system and the (rich) dynamics it generates in \Cref{sec:simulations}. We conclude the paper with a Discussion section.

\begin{remark}
\label{rem:k*=0}
\rm  
Throughout this paper, we assume that the Turing bifurcation described by the Ginzburg-Landau equation occurs at a critical wave number $k^c$ that remains bounded away from zero as the Turing bifurcation merges with the fold bifurcation. That is, we assume that $k^c \to k^* > 0$ as the governing pair of parameters $(\mu,\nu)$ approaches the co-dimension 2 point $(\mu^*,\nu^*)$. While this is a generic and natural assumption, it entails that the scalar higher-order phase-field models we consider must be at least sixth order, and that the $n$-component reaction-diffusion models must have at least three components. A straightforward linear stability analysis shows that $k^c \to k^* = 0$ in the case of a Turing-fold bifurcation for both scalar fourth-order (Swift-Hohenberg type) equations and two-component reaction-diffusion systems. In such cases, the Ginzburg-Landau approximation degenerates as the Turing-fold point is approached: it transforms into an extended Fisher-Kolmogorov/Swift-Hohenberg-type equation \cite{rottschafer1998transition}. Accordingly, the associated bifurcation obtains a co-dimension 3 character, which may be studied by combining the present approach with that of \cite{rottschafer1998transition}.  
\end{remark}

\begin{remark}
\label{rem:backgroundFig1}
\rm
The $\eta = -1$ simulation of {\it Turing-triggers-tipping} shown in \Cref{fig:Intro-TTT} is chosen for illustrative `dramatic' reasons. Keeping all other parameter values equal, there exist values of $\eta \in (-1, 1)$ (with $\eta \geq 1$ corresponding to $\beta \geq 0$ in the associated AB-system; see \Cref{subsec:a model}) where the system does not collapse to the trivial $U(x,t) \equiv 0$ state, but instead transitions to another spatially periodic state. Thus, also for these values of $\eta$, pattern formation is not an early warning signal for tipping, but a mechanism by which collapse (to the trivial state) is evaded. Both types of behavior -- direct collapse to a trivial state and transition to alternative spatial patterns -- also arise in a recent model for tree–grass dynamics  \cite{van2025vegetation}. In particular, transitions to alternative spatially periodic states when the Turing bifurcation itself is subcritical (thus there are no stable small-amplitude spatially-periodic states), such as in \eqref{eq:scalar example} for $\eta < 1$, are not uncommon; see, for instance, \cite{champneys2021bistability,elvin2009transient,guisoni2022transient, knobloch2015spatial} and the references therein. Moreover, for similar illustrative reasons, we have chosen the magnitude of $\delta$ beyond the asymptotically small range on which the derivation of the system of modulation equations \eqref{eq:canonical form:csd:gs} is based. Specifically, it cannot be assumed that the choice $\delta = 0.6$ satisfies the condition `$\delta$ is sufficiently small' that underlies the asymptotic analysis of \Cref{sec:scalar}. This, combined with our choice of domain length $L = 4 \pi$, may be the reason that the Turing bifurcations are somewhat delayed in \Cref{fig:Intro-Turing,fig:Intro-TTT}. Namely, the choice $L = 4 \pi$ allows for two full periods of wavelength 1 Turing patterns. However, in practice, the critical wave number $k^c$ is not equal to $1$, but to $k^c = 1 +\frac14 \delta + \mathcal{O}(\delta)^2$ (see \Cref{subsubsec:basic analysis:a model}); thus, for $\delta =0.6$, the actual value of $k^c$ may deviate significantly from $1$, so that a two-period pattern can only be expected to bifurcate for values of $\mu$ strictly smaller than $\mu_t$. Additionally, it is known that a Turing bifurcation driven by a slowly varying parameter will also be delayed under more perfect circumstances \cite{avitabile2020local,chen2015patterned}.
\end{remark}

\begin{remark}
\label{rem:validity}
\rm    
In this paper, we do not address the important issue of rigorously justifying the validity of the AB-system of modulation equations as a {\it normal form} describing the dynamics of small amplitude solutions near a co-dimension 2 Turing-fold point (for systems of evolutionary PDEs on cylindrical domains). In \Cref{sec:simulations}, we numerically verify that the dynamics generated by the AB-system persist in the underlying PDE -- even when the patterns exhibited are more complex than the standard plane waves. Preliminary analysis strongly suggests that the procedure by which the validity of the AB-system can be established through an approach analogous to that used for the  (`classical, co-dimension 1') Ginzburg-Landau equation, , both for general patterns on finite time intervals and for long-time dynamics and attractors as $t \to \infty$ (see \cite{mielke2002ginzburg,mielke1995attractors,schneider1998nonlinear}). This is the subject of ongoing work. 
\end{remark}

\section{Scalar equations}\label{sec:scalar}

In this section we consider the Turing-fold bifurcation in (one-dimensional) higher-order scalar models. In \Cref{subsec:a model} we illustrate the essential properties of this co-dimension 2 bifurcation via our basic scalar example \eqref{eq:scalar example}. In \Cref{subsec:general setting:se} we study general higher-order scalar systems and demonstrate why the properties highlighted in \Cref{subsec:a model} are characteristics of the Turing-fold bifurcation. Furthermore, the analysis in \Cref{sec:scalar} elucidates our rationale for selecting this example: unless the scalar equation is at least of sixth-order, the critical wave number associated with the co-dimension 2 Turing-fold bifurcation must be zero, which introduces further degeneracies -- see \Cref{rem:k*=0,rem:sixth-order}. Additionally, we will show that the Turing bifurcation near the Turing-fold point can only be supercritical (in the class of scalar models considered here -- see \Cref{subsubsec:ginzburg-landau derivation:gs}) if the nonlinear terms contain spatial derivatives (see \Cref{subsubsec:coupled system derivation:a model}) -- which is why we introduced the $\eta (\partial_x^2 U)^2$ term.

\subsection{The Turing-fold bifurcation in model \eqref{eq:scalar example}}\label{subsec:a model}

In \Cref{subsubsec:basic analysis:a model}, the basic (linear) bifurcation analysis for \eqref{eq:scalar example} -- as a function of bifurcation parameters $\mu$ and $\nu$ -- is presented. To study the dynamics near the co-dimension 2 Turing-fold bifurcation point $(\mu,\nu) = (\mu^*,\nu^*)$ ($=(-1,1)$ as we shall find), we introduce the artificial,  sufficiently small, parameter $0 < \delta \ll 1$ and set $\nu = \nu^* - \delta$. We show that this leads to the relation $|\mu^t - \mu^*| = \mathcal{O}(\delta^2)$ -- where $\mu^t = \mu^t(\nu)(> \mu^*)$ is the value at which one of the background states of \eqref{eq:scalar example} resulting from the saddle-node bifurcation undergoes a Turing bifurcation.  To focus on the interactions between the Turing and fold bifurcations, we introduce a scaled bifurcation parameter $r$ by setting $\mu = \mu^t - r\delta^2$. This order-$\delta^2$ scaling allows us to look beyond the saddle-node point (for $r$ sufficiently large), while still capturing the effects of the Turing bifurcation by considering $\mu$ sufficiently close to $\mu^t$. In \Cref{subsubsec:ginzburg-landau derivation:a model}, assume that $\mu$ is sufficiently close to $\mu^t$ to derive (which we subsequently do) the Ginzburg-Landau equation that governs the (weakly) nonlinear dynamics of \eqref{eq:scalar example} near the Turing point. Moreover, we show that sufficiently close must mean $|\mu - \mu^t| \ll \delta^2$ since the standard Ginzburg-Landau approach degenerates as one attempts to derive the governing modulation equations for $\mu = \mu^t - r \delta^2$. In \Cref{subsubsec:coupled system derivation:a model}, we adapt the Ansatz and derive a coupled system of modulation equations -- which we refer to as the AB-system -- that takes over the role of the Ginzburg-Landau equation: describing the dynamics near and beyond the fold-point in the Turing-fold setting.

\subsubsection{Linear stability analysis}\label{subsubsec:basic analysis:a model}
\begin{figure}[t]
    \centering
    \includegraphics[width=0.6\linewidth]{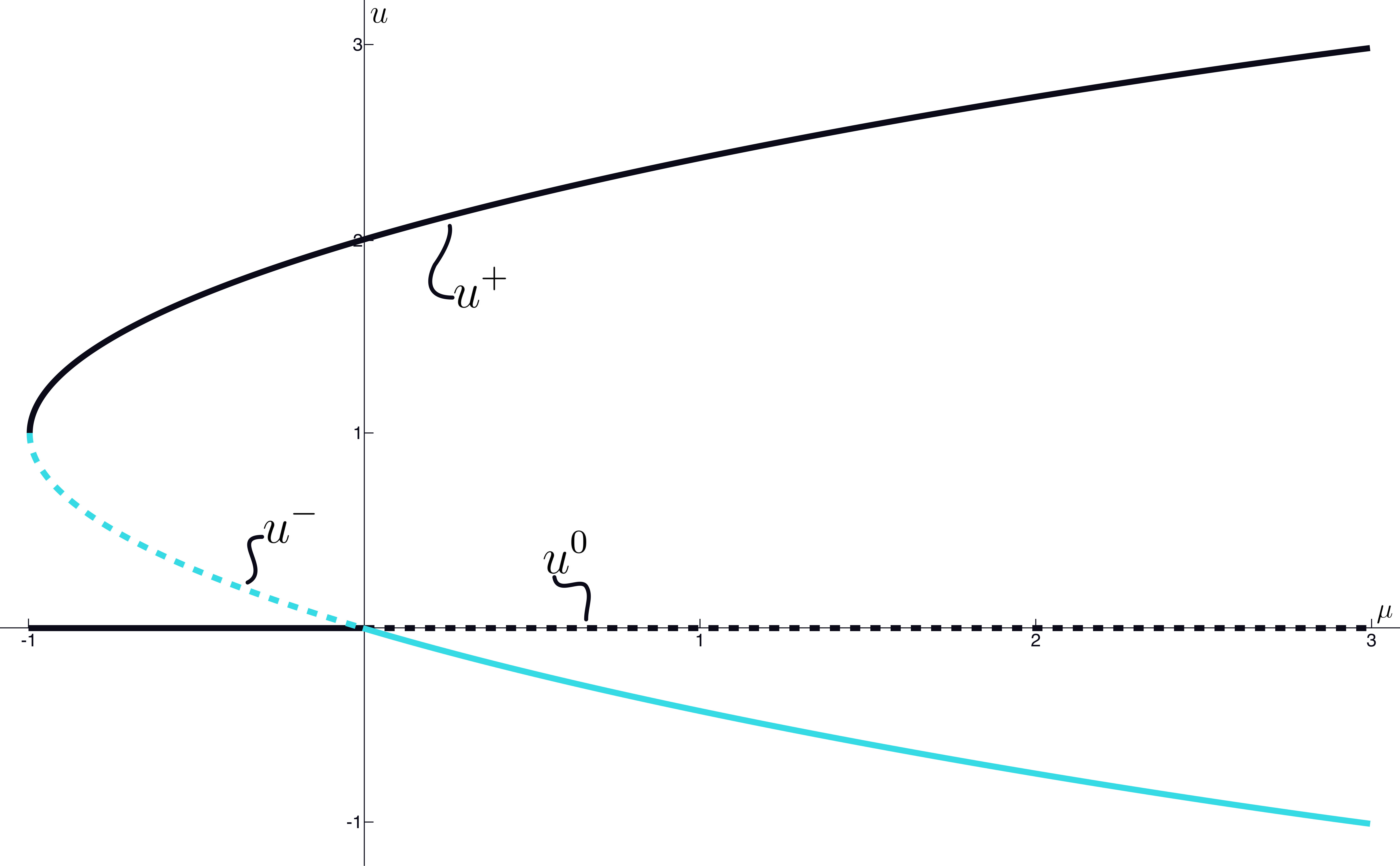}
    \caption{The bifurcation diagram associated to \eqref{eq:ode problem associated with scalar example}.}
    \label{fig:bifurcation diagram ode:a model}
\end{figure}

For $\mu > -1$, equation \eqref{eq:scalar example} has 3 spatially homogeneous states:
\begin{equation}
\label{model-homstates}
    u^0 = 0, \quad \text{and} \quad u^\pm(\mu) = 1 \pm \sqrt{1 + \mu}.
\end{equation}
The pair $u^\pm(\mu)$, with $u^+ > u^-$, results from a saddle-node bifurcation as $\mu$ passes through $\mu^* \coloneqq -1$. Naturally, these states are critical points of the spatially homogeneous ODE associated to \eqref{eq:scalar example},
\begin{equation}\label{eq:ode problem associated with scalar example}
    \dot{u} = F(u;\mu) \coloneqq \mu u + 2u^2 - u^3.
\end{equation}
Since,
\begin{equation}
\label{Fuu0upm}
F_u\left(u^0; \mu\right) = \mu, \quad \text{ and } \quad F_u\left(u^\pm(\mu);\mu\right)\coloneqq F_u\left(u^\pm; \mu\right) = 2 \left(-1 - \mu \mp \sqrt{1 + \mu}\right),
\end{equation}
we note that $u^+$ is stable as critical point of \eqref{eq:ode problem associated with scalar example} (thus, as a solution of $\eqref{eq:scalar example}$, $u^+$ is stable with respect to spatially homogeneous perturbations) and that $u^0$ and $u^-$ exchange stability in a transcitical bifurcation as $\mu$ passes through $\mu = 0$ -- see the bifurcation diagram associated to \eqref{eq:ode problem associated with scalar example} in \Cref{fig:bifurcation diagram ode:a model}.

To investigate the stability of $u^+$ against nonhomogeneous perturbations, we set
\begin{equation}\label{eq:perturbation:a model}
    U = u^+ + \left(\bar{u} e^{ikx + \omega(k; \mu, \nu) t} + \text{c.c.}\right),
\end{equation}
and linearize,
\begin{equation}\label{eq:spectral stability equation:ba:am:se}
    \omega(k; \mu, \nu) = \mu + 4 u^+ - 3 \left(u^+\right)^2 - \nu k^2 + 2 k^4 - k^6 \coloneqq F_u\left(u^+; \mu\right) + G(k; \nu).
\end{equation}
Suppose a Turing bifurcation occurs at $\l \mu,k\r  = \l \mu^t,k^c \r = \l \mu^t\l \nu\r ,k^c\l \nu\r \r $, then necessarily $\omega\l k^c; \mu^t, \nu\r  = 0$ and $\omega_k\l k^c; \mu^t, \nu\r  = 0$. The latter,  
\begin{equation}\label{eq:kc highlight non-degeneracy:a model}
    \omega_k\l k^c; \mu^t,\nu\r  = G_k\l k^c;\nu\r  = - 2 k^c\l \nu - 4\l k^c\r ^2 + 3 \l k^c\r ^4\r  = 0,
\end{equation}
has solutions:
\begin{equation}
\label{kc0kcpm}
    k^c_0 = 0, \quad \left(k^c_-\right)^2 = \frac{1}{6}\left(4 - \sqrt{16 - 12 \nu} \right), \quad \text{ and } \quad \left(k^c_+\right)^2 = \frac{1}{6}\left(4 + \sqrt{16 - 12 \nu} \right),
\end{equation}
which determine $k^c(\nu)$. Now suppose that $\nu > 0$ -- which we will do from now on -- then it follows that $k^c_0 = 0$ is a maximum of the sixth-order polynomial $G(k; \nu)$ \eqref{eq:spectral stability equation:ba:am:se}. Hence, $k^c_-$ is a minimum and $k^c_+$ is also a maximum. Alternatively, if $\nu < 0$, then $k^c = 0$ is a minimum, $k^c_-$ has disappeared and $k^c_+$ is still a maximum. Consequently, since $G(0;\nu) = 0$ for all $\nu \in \mathbb{R}$, and $F_u(u^+;\mu) < 0$ for all $\mu > \mu^*$ \eqref{Fuu0upm}, \eqref{eq:spectral stability equation:ba:am:se} implies that the critical wave number is $k^c_+$ if $u^+$ destabilizes via a Turing bifurcation at $\mu = \mu^t > \mu^*$. Henceforth, we denote $k^c_+$ simply by $k^c$.

Next, we observe that,
\begin{equation}
\label{Gkcnu}
G(k^c(\nu);\nu) = \frac{2}{27}\left(\sqrt{4-3\nu} + 2\right)^2\left(\sqrt{4-3\nu}-1\right),
\end{equation}
which implies that $G(k^c;\nu) > 0$ if and only if $\nu < 1$. Hence, it follows from \eqref{Fuu0upm} and \eqref{eq:spectral stability equation:ba:am:se} that for all $\nu < 1$, there exists a unique $\mu^t = \mu^t(\nu) > \mu^*$ such that
$\omega(k^c(\nu); \mu^t(\nu), \nu) = 0$. Additionally,
\begin{equation}\label{eq:omega mu:a model}
\omega_\mu\l k^c;\mu^t,\nu\r  =  F_{uu}\l u^+\l \mu^t\r ;\mu^t\r u^+_\mu\l \mu^t\r  + F_{u\mu}\l u^+\l \mu^t\r ;\mu^t\r = -2 - \frac{1}{\sqrt{1+\mu^t}} < 0.
\end{equation}
from which we conclude that $u^+$ destabilizes via a Turing bifurcation with critical wave number $k^c$ as $\mu$ decreases through $\mu^t$. Since $G(k^c(\nu);\nu) \downarrow 0$ as $\nu \uparrow 1$ \eqref{Gkcnu} and 
$F_u(u^+(\mu);\mu) \uparrow 0$ as $\mu \downarrow -1$ \eqref{Fuu0upm}, we indeed find that the Turing and fold bifurcations coincide at the co-dimension 2 Turing-fold bifurcation point $(\mu,\nu) = (\mu^*,\nu^*) = (-1,1)$. Moreover, this (degenerate) Turing bifurcation has critical wave number $k^c(\nu^*) \coloneqq k^* = 1$ \eqref{kc0kcpm}.

\begin{figure}[t]
\centering
\includegraphics[width=1.0\linewidth]{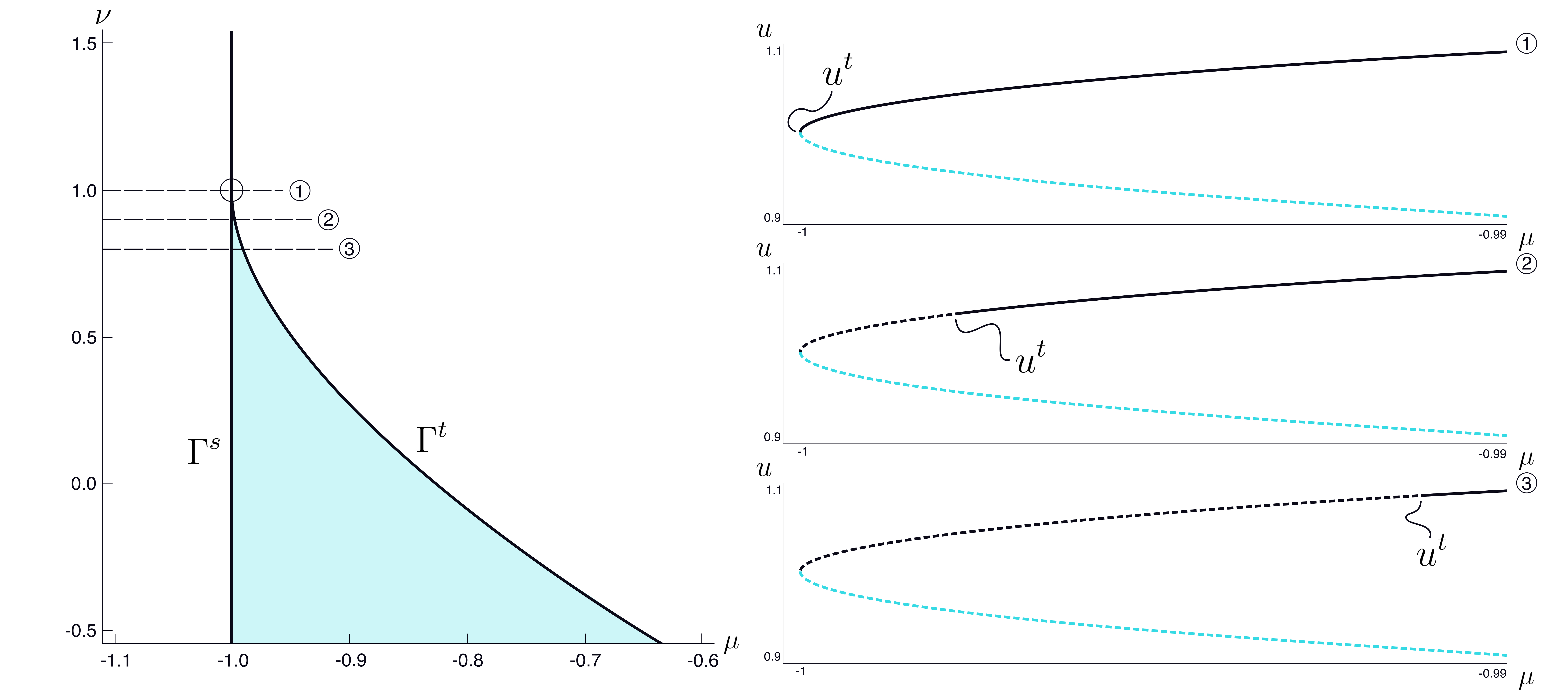}
\caption{
In the left image, the curves $\Gamma^s$ of saddle-node points and $\Gamma^t$ of Turing points are plotted in the $(\mu, \nu)$-plane. These curves meet tangentially at the co-dimension 2 point $(\mu^*, \nu^*) = (-1, 1)$. The blue area between the two curves represents the region in parameter space where the Turing bifurcation has occurred, but the saddle-node bifurcation has not. On the right, three bifurcation diagrams (as function of $\mu$) are plotted for three different values of $\nu$, with $\nu_1 = \nu^* = 1$ and $\nu_1 < \nu_2 < \nu_3$; black segments represent $u^+$ and the blue segments represent $u^-$.}
\label{fig:pde bifurcation diagram scalar example:a model}
\end{figure}

Naturally, we are interested in the dynamics of \eqref{eq:scalar example} as the Turing bifurcation approaches the fold bifurcation. Therefore, we introduce a small parameter $0 < \delta \ll 1$, and set $\nu = \nu^* - \delta$ so that $k^c = 1 + \frac{\delta}{4} + \mathcal{O}(\delta^2)$ \eqref{kc0kcpm}. Since $\omega(k; \mu, \nu) = 0$ at the Turing bifurcation, it follows by \eqref{Fuu0upm} and \eqref{eq:spectral stability equation:ba:am:se} that
\[
F_u(u^+;\mu^t) = -2 \left(1 + \mu^t +\sqrt{1 + \mu^t}\right) = -\delta + \mathcal{O}(\delta^2) = - G(k^c; 1-\delta).
\]
Introducing the notation $u^+(\mu^t) = u^t$, we thus have,
\begin{equation}\label{eq: scalar example PDE}
    \mu^t = -1 + \frac{\delta^2}{4} + \mathcal{O}\l \delta^3\r  \quad \text{ and } \quad u^t = 1 + \frac{\delta}{2} + \mathcal{O}\l \delta^2\r .
\end{equation}
Additionally, we note that $F_u\l u^t;\mu^t\r$ is asymptotically small,
\begin{equation}
\label{Fuu+mut}
F_u\l \mu^t;\mu^t\r = -\delta + \mathcal{O}\l\delta^2\r,    
\end{equation}
which is a key ingredient driving the upcoming analysis. Moreover, we conclude from \eqref{eq:omega mu:a model} and \eqref{eq: scalar example PDE} that $\omega_\mu\l\mu^t\r = \mathcal{O}_s\l \delta^{-1}\r$ (see Remark \ref{rem:Os}), more precisely,
\begin{equation}
\label{wmutO1delta}
\omega_\mu\l\mu^t\r = \frac{-2 + \mathcal{O}(\delta)}{\delta} = \mathcal{O}_s\l\delta^{-1}\r
\end{equation}
In conclusion, we have found that, 
\begin{equation}
\label{kcuutmut}
    k^c = k^* + \tilde{k}\delta + \mathcal{O}\l \delta^2\r, \quad
    u^t = u^* + \tilde{u}\delta + \mathcal{O}\l \delta^2\r, \quad
    \mu^t = \mu^* + \hat{\mu}\delta^2 + \mathcal{O}\l \delta^3\r, \quad
    \omega^t_\mu = \frac{\tilde{\omega} + \mathcal{O}(\delta)}{\delta},
\end{equation}
with $\tilde{u} \neq 0$, $\hat{\mu} \neq 0$, $k^* \neq 0$, and $\tilde{\omega} \neq 0$: these are the four essential properties of the Turing-fold bifurcation that will be used in \Cref{subsubsec:coupled system derivation:a model} to justify the adapted Ansatz that allows us to derive the AB-system. Beyond the obvious characteristics of the Turing-fold bifurcation, namely the presence of both Turing and fold bifurcations, the non-degeneracy conditions in \Cref{subsubsec:non-degeneracy conditions:gs,subsec:non-degeneracy conditions:rd} ensure that anything we call a Turing-fold bifurcation has these four additional properties. Note that one of the consequences of these properties is that $\Gamma^s$, where $\Gamma^s$ is defined as the curve in $(\mu, \nu)$-space of saddle-node points (i.e., $\Gamma^s = \{(-1, \nu): \nu \in \mathbb{R}\}$), and $\Gamma^t$, where $\Gamma^t$ denotes the curve of Turing points (i.e., $\Gamma^t = \{(\mu^t(\nu), \nu), \nu < 1 \}$), merge with a quadratic tangency at $(\mu^*, \nu^*)$ -- see \Cref{fig:pde bifurcation diagram scalar example:a model}.

To simplify notation in the upcoming derivations, we mark $F$, $\omega$, and their derivatives by a superscript $t$ to indicate evaluation at $(u^t; \mu^t)$ and $(k^c;\mu^t, \nu^* - \delta)$, respectively. For example, we write:
\[
    F^t \coloneqq F\l u^t;\mu^t\r,  \quad F_{uu}^t \coloneqq F_{uu}\l u^t;\mu^t\r,  \quad \omega^t \coloneqq \omega\l k^c,\mu^t,\nu^*-\delta\r, \quad \text{ and } \quad \omega^t_k = \omega_k\l k^c, \mu^t, \nu^* - \delta\r.
\]

\begin{remark}
\label{rem:Os}  
\rm
We use the notation $\mathcal{O}_s(\tilde{\delta}^j)$, $j=-1,0,1$, where $\tilde{\delta}$ is an asymptotically small quantity, to indicate that a specific expression is exactly of order $\tilde{\delta}^j$. More precisely, the function $h(\tilde{\delta}) = \mathcal{O}_s(\tilde{\delta}^j)$ ($j=-1,0,1,2$) if there exist constants $0 < C_1 < C_2$ such that $C_1 \tilde{\delta}^j < |h(\tilde{\delta})| < C_2 \tilde{\delta}^j$ as $\tilde{\delta} \to 0$ (cf. \cite{eckhaus1979asymptotic}).
\end{remark}

\subsubsection{Weakly nonlinear analysis: the Ginzburg-Landau approximation}\label{subsubsec:ginzburg-landau derivation:a model}

To approximate the dynamics near a Turing bifurcation, one usually derives the Ginzburg-Landau equation as the governing modulation equations near this point. However, this approach is unsuitable for understanding the dynamics in a radius large enough to include the saddle-node bifurcation, since the fold and Turing bifurcations interact. To illustrate the precise issues we encounter when applying the standard Ginzburg-Landau approach, we will derive the Ginzburg-Landau equation for \eqref{eq:scalar example} for $\mu$ sufficiently close to the Turing bifurcation.

We set $\mu = \mu^t - r\varepsilon^2$ ($r \in \mathbb{R}$) with a priori $0 < \varepsilon \ll \delta \ll 1$ (naturally we keep $\nu = 1 - \delta$) and impose the standard Ginzburg-Landau approximation \cite{schneider2017nonlinear}, 
\begin{equation*}
    U_{GL}(x,t) = u^t + \left(e^{i k^c x}A\left(\varepsilon x, \varepsilon^2 t\right) + \text{c.c.}\right)
\end{equation*}
(at leading order), for solutions $U$ of \eqref{eq:scalar example} $\mathcal{O}(\varepsilon^2)$ close to the marginally stable/unstable state $u^t$ -- where $A(\xi,\tau): \mathbb{R} \times \mathbb{R}^+ \to \mathbb{C}$ is an a priori unknown amplitude that varies slowly in space and time: $\xi = \varepsilon x$ and $\tau = \varepsilon^2 t$. The (weakly) nonlinear evolution is captured by the standard Ginzburg-Landau Anzatz:
\begin{equation}\label{eq:scalar ginzburg-landau ansatz:ginzburg-landau derivation}
    \begin{split}
        E^0 \big[&\varepsilon^2 X_{02} + \varepsilon^3 X_{03} + \mathcal{O}\l\varepsilon^4\r \big] \\
        U_{\text{GL}}(x,t) = u^t+E \big[\varepsilon A + &\varepsilon^2 X_{12} + \varepsilon^3 X_{13} + \mathcal{O}\l\varepsilon^4\r \big] + \text{c.c.} \\
        E^2 \big[&\varepsilon^2 X_{22} + \varepsilon^3 X_{23} + \mathcal{O}\l\varepsilon^4\r \big] + \text{c.c.} \\
         & \phantom{iiiiii}\, E^3 \big[\varepsilon^3 X_{33} + \mathcal{O}\l\varepsilon^4\r \big] + \text{c.c.},
    \end{split}
\end{equation}
where we introduced $E = e^{ik^cx}$, $A(\xi,\tau)$, and the unknown functions $X_{0j}(\xi,\tau): \mathbb{R} \times \mathbb{R}^+ \to \mathbb{R}$ and $X_{ij}(\xi,\tau): \mathbb{R} \times \mathbb{R}^+ \to \mathbb{C}$, $i \geq 1$, which will be determined in terms of $A(\xi,\tau)$ by substituting \eqref{eq:scalar ginzburg-landau ansatz:ginzburg-landau derivation} into \eqref{eq:scalar example} and collecting the terms at the $\varepsilon^jE^i$-levels (see for instance \cite{doelman2019pattern} and the references therein).

The terms at the $\varepsilon^0 E^0$-, $\varepsilon E$ and $\varepsilon^2 E$-levels cancel by construction. At the $\varepsilon^2 E^0$- and $\varepsilon^2 E^2$-levels we find:
\[
- F_{u}^t X_{02} = \left(F_{uu}^t + 2 \eta (k^c)^4 \right)|A|^2
\quad \text{ and } \quad
-\left(F_u^t + G\l2k^c; \nu\r\right)X_{22} = \frac{1}{2}\left(F_{uu}^t + \eta \l k^c\r^4\right)A^2.
\]
The leading order Ginzburg-Landau equation (in $\varepsilon$) is retrieved at the $\varepsilon^3E$-level:
\begin{equation}\label{eq:scalar example gl:ginzburg-landau derivation}
    A_\tau = - \frac{1}{2}\omega_{kk}^t A_{\xi\xi} - r \omega_{\mu}^t A + L |A|^2A,
\end{equation}
with $\omega_{\mu}^t$ as in \eqref{wmutO1delta},
\[
\omega_{kk}^t = G_{kk}\l k^c,\nu\r = -2 \nu + 24 \l k^c\r^2 - 30 \l k^c\r^4 = -8 + \mathcal{O}\l \delta\r, 
\]
\eqref{kc0kcpm} and the (Landau) coefficient
\[
L = -3-\frac{\l F_{uu}^t + 8 \eta \l k^c\r ^4\r \l F_{uu}^t + \eta \l k^c\r ^2\r }{2\l F_u^t + G\l 2 k^c; \nu\r \r }-\frac{F_{uu}^t\l F_{uu}^t + 2 \eta \l k^c\r ^4\r }{F_u^t} = \frac{4\l 1 - \eta \r }{\delta} + \mathcal{O}\l 1\r;
\]
thus, \eqref{eq:scalar example gl:ginzburg-landau derivation} reduces to:
\begin{equation}\label{eq:ginzburg-landau:a model}
A_\tau = \left(4 + \mathcal{O}(\delta)\right)A_{\xi \xi} + \frac{r}{\delta} \left(2 + \mathcal{O}(\delta) \right) A + \frac{4}{\delta} \left(1-\eta + \mathcal{O}(\delta) \right)|A|^2 A,
\end{equation}
at leading order in $\varepsilon$.

The procedure above provides a valid approximation of the dynamics of small-amplitude patterns of \eqref{eq:scalar example} for $(U,\mu)$ $\mathcal{O}(\varepsilon^2)$ close to $(u^t,\mu^t)$ as long as $\varepsilon \ll \delta$ ($\ll 1$) -- see \Cref{rem:validity}. However, since $|\mu^* - \mu^t| = \mathcal{O}(\delta^2)$ \eqref{kcuutmut} we need to increase $\varepsilon$ to $\varepsilon = \mathcal{O}(\delta)$ to describe the dynamics near and beyond the saddle-node bifurcation. It is straightforward to check that the equation for $X_{0(2n+2)}$ at the $\varepsilon^{2n+2} E^0$-level has the structure:
\[
    -F_u^t X_{0(2n+2)} = F_{uu}^t X_{0(2n)} X_{02} + \dots
\]
(for all $n \geq 1$). Hence, we conclude that we conclude that $X_{0(2n)} = \mathcal{O}(\delta^{-2n + 1})$, since $F_u^t = \mathcal{O}_s(\delta)$ \eqref{Fuu+mut}. Consequently, higher-order terms as $|A|^{2n}A$ -- that appear at order $\varepsilon^{2n-2}$ in the standard setting and thus are neglected in the leading order Ginzburg-Landau equation \eqref{eq:scalar example gl:ginzburg-landau derivation}/\eqref{eq:ginzburg-landau:a model} -- are, in fact, terms of order $\delta^{-2n + 1} \varepsilon^{2n-2}$ in the (extended) Ginzburg-Landau equation. As a result, all $|A|^{2n}A$ terms grow to $\mathcal{O}(\delta^{-1})$ as $\varepsilon$ is increased to $\mathcal{\delta}$: these higher-order effects thus all increase to the same order as leading order term $|A|^2A$. This completely invalidates the standard approximation procedure. 

\begin{remark} {\rm
    The Turing bifurcation is subcritical (i.e., no stable spatially periodic patterns exist in \eqref{eq:ginzburg-landau:a model}) if $\eta < 1$ and supercritical if $\eta > 1$. This is why we introduced the $\eta(\partial_x^2 U)^2$ term to \eqref{eq:scalar example}. As we will show in \Cref{subsec:a model}, a supercritical Turing bifurcation cannot occur in this Turing-fold setting if the scalar PDE includes only linear spatial derivative terms.}
\end{remark}

\subsubsection{Weakly nonlinear analysis: the coupled AB-system}\label{subsubsec:coupled system derivation:a model}
Keeping $\nu = 1 - \delta$, we now set $\mu = \mu^t - r \delta^2$, which implies (by \eqref{eq: scalar example PDE}) that $\mu < \mu^*$ for $r > 1/4$: the bifurcation parameter $r$ connects the Turing bifurcation (at $r = 0$) via the saddle-node bifurcation (at $r = 1/4$) to the `tipped' situation in which the (Turing unstable) background state $u^+$ \eqref{model-homstates} has ceased to exist (for $r > 1/4$). 

However, as we saw in \Cref{subsubsec:ginzburg-landau derivation:a model}, we must adapt the standard Ginzburg-Landau procedure to describe the dynamics of patterns in \eqref{eq:scalar example} near the co-dimension 2 point $(\mu^*,\nu^*)$. First, we note that it is necessary to introduce a second (now real) amplitude $B(\xi,\tau)$ to capture the impact of the saddle-node bifurcation. The standard Ansatz \eqref{eq:scalar ginzburg-landau ansatz:ginzburg-landau derivation} centers around $u^t$, the present co-dimension 2 setting, however, asks for replacing $u^t$ by $u^*$, because $u^t$ depends on $\delta$. The difference between $u^t$ and $u^*$ -- which is of $\mathcal{O}_s(\delta)$ -- thus needs to be incorporated in the saddle-node mode. Therefore, we conclude that the effect of the saddle-node bifurcation should be governed by the term $\delta B(\xi,\tau)$ (at leading order in $\delta$).

Second, we also need to adapt the scaling of $\xi$, and $\tau$. To see this, we (re)examine the linear destabilization process, which is, by definition, driven by $\omega(k; \mu^t - r \delta^2, 1-\delta)$ \eqref{eq:perturbation:a model}. Specifically, set $k = k^c + K$ with $K$ small, so that (at leading order),
\begin{equation*}
\omega( k^c + K; \mu^t - r \delta^2, 1-\delta) = \frac{1}{2} \omega_{kk}^t K^2 - \omega_{\mu}^t r \delta^2 = -4 K^2 + 2r \delta,
\end{equation*}
\eqref{eq:spectral stability equation:ba:am:se}, \eqref{wmutO1delta} and observe that there exist (linearly) unstable spatial periodic patterns for $K \in (-\sqrt{r/2} + \mathcal{O}(\delta), \sqrt{r/2} + \mathcal{O}(\delta))$, i.e., for $K = \mathcal{O}_s(\sqrt{\delta})$. Additionally, note that the linear growth rate of these patterns, i.e., the magnitude of $\omega(k^c + K, \mu^t - r \delta^2, 1-\delta)$, is of $\mathcal{O}_s(\delta)$. Therefore, we conclude that in the co-dimension 2 Turing-fold setting, the natural spatial and temporal scales are $\xi = \sqrt{\delta} x$ and $\tau = \delta t$  -- in contrast to the Ginzburg-Landau setting where the standard choices for these scales are $\mathcal{O}(\delta)$ and $\mathcal{O}(\delta^2)$ (for $\mu$ at an $\mathcal{O}(\delta^2)$ of its bifurcational value) -- see section \ref{subsubsec:ginzburg-landau derivation:a model}. 

Finally, we need to (re)consider the scaling of $A(\xi,\tau)$. However, there is a priori, in contrast to a posteriori (\Cref{rem:degenerate}), no reason to adapt this scaling. Therefore, the $\mathcal{O}(\varepsilon) = \mathcal{O}(\sqrt{|\mu^t - \mu|})$ scaling of the classical Ginzburg-Landau Ansatz \eqref{eq:scalar ginzburg-landau ansatz:ginzburg-landau derivation} is retained here: $\sqrt{|\mu^* - \mu|} = \mathcal{O}(\delta)$. 

Together, this yields the following adapted Ansatz to capture the weakly nonlinear behavior of patterns in the Turing-fold bifurcation:
\begin{equation}\label{eq: Coupled system Ansatz}
    \begin{split}
    U_{AB}(x, t) = 1 + E^0 \big[\delta B +\delta^{3/2}X_{02} + &\delta^{2} X_{03} +\delta^{5/2} X_{04}+\mathcal{O}\l \delta^3\r\big] \\
    + E \big[\delta A + \delta^{3/2} X_{12} + &\delta^{2} X_{13} + \delta^{5/2} X_{14} +\mathcal{O}\l \delta^3\r \big] + \text{c.c.} \\
    +E^2 \big[&\delta^2 X_{23} + \delta^{5/2} X_{24} + \mathcal{O}\l \delta^3\r \big] + \text{c.c.} \\
     &\phantom{iiiii00000000000iiiiii+E_c^3}\,  + \text{h.o.t.},
\end{split}
\end{equation}
where $E = e^{i k^* x} = e^{i x}$ (Note that $k^* \neq 0$ is indeed crucial: the decomposition \eqref{eq: Coupled system Ansatz} breaks down if $k^* = 0$ -- see \Cref{rem:k*=0}.) The unknown functions $X_{0j}(\xi,\tau): \mathbb{R} \times \mathbb{R}^+ \to \mathbb{R}$ and $X_{ij}(\xi,\tau): \mathbb{R} \times \mathbb{R}^+ \to \mathbb{C}$, $i \geq 1$ will be determined in terms of amplitudes $A(\xi,\tau): \mathbb{R} \times \mathbb{R}^+ \to \mathbb{C}$ and $B(\xi,\tau): \mathbb{R} \times \mathbb{R}^+ \to \mathbb{R}$, again by substitution and collecting terms at the $\delta^j E^i$-levels. 

By construction, we (again) find that the $\delta E^i$ and $\delta^{3/2} E^i$-levels reduce to $0 = 0$. However, at the $\delta^2 E^0$-level, we obtain:
\begin{equation}
    B_{\tau} = B_{\xi \xi} + \frac{1}{4} - r - B^2 + 2(\eta - 1)|A|^2.
\end{equation}
Additionally, the $\delta^2 E$-level yields:
\begin{equation}
    A_{\tau} = 4 A_{\xi \xi} + A - 2 A B.
\end{equation}
Consequently, we have found a new system of modulation equations that describes the leading order dynamics of small-amplitude solutions of the model \eqref{eq:scalar example} near the co-dimension 2 Turing-fold point $(\mu^*,\nu^*) = (-1,1)$:
\begin{equation}\label{eq:coupled system:a model}
    \begin{cases}
        A_{\tau} = 4 A_{\xi \xi} + A - 2 A B, \\
        B_{\tau} = B_{\xi \xi} + \frac{1}{4} - r - B^2 + 2(\eta - 1)|A|^2.
    \end{cases}
\end{equation}
As we will show in \Cref{subsec:general setting:se,sec:reaction diffusion system} and discuss in \Cref{sec:discussion}, this system of modulation equations -- that we call the AB-system -- represents the modulation system that is universal within the Turing-fold setting (the nature of the underlying model only changes the specific values of the coefficients of \eqref{eq:coupled system:a model}). Additionally, in \Cref{subsubsec:Reducing the AB-system:gs:se} we show that the classic Ginzburg-Landau equation is, in the general setting, embedded in the AB-system by reintroducing $\varepsilon$. Finally, in \Cref{sec:coupled system}, we will analyze the coupled system in its canonical form. Among other insights, we will demonstrate that there is a supercritical bifurcation of periodic solutions if the coefficient of the $|A|^2$ term is positive, i.e., if $\eta > 1$ in \eqref{eq:coupled system:a model}, and a subcritical bifurcation of periodic patterns if this coefficient is negative. Moreover, we will (formally) validate the statement of Claim \ref{thm:close enough:t:cs} and show that the stable periodic patterns of the supercritical case persist for values of $r$ beyond the saddle-node bifurcation at $r = \frac{1}{4}$: the Turing patterns enable the system to evade tipping -- see also \Cref{fig:Intro-Turing} in which $\eta > 1$ and \Cref{fig:Intro-TTT} with $\eta < 1$.

\begin{remark}
\label{rem:degenerate}
\rm    
In the analysis of \Cref{sec:coupled system}, we show that the (sub/supercritical) nature of the Turing bifurcation in \eqref{eq:coupled system:a model} is determined by the sign of the coefficient of the $|A|^2$ term in the $B$ equation. The transition from a sub- to a supercritical bifurcation can be studied by following the approach of \cite{doelman1991periodic, shepeleva1997validity} in which this transition has been studied in the setting of the Ginzburg-Landau equation. Therefore, we introduce $\gamma$ and $\tilde{\eta}$ by writing $\eta = 1 + \tilde{\eta} \delta^\gamma$, where $\gamma \geq 0$ measures the relative magnitude of the coefficient of the $|A|^2$ term in terms of $\delta$ -- the magnitude of the distance between the saddle-node and the Turing bifurcation. (Note that we are thus zooming in into a co-dimension 3 bifurcation (at $(\mu^*,\nu^*, \eta^*) = (-1, 1, 1)$) -- unlike the co-dimension 2 Ginzburg-Landau case considered in \cite{doelman1991periodic,eckhaus1989strong,shepeleva1997validity}.) As a consequence, the magnitude of the coefficient of the $|A|^2$ term in the AB-system \eqref{eq:coupled system:a model} becomes of $\mathcal{O}(\delta^{\gamma})$, which implies that the magnitude of the plane wave solutions \eqref{eq:periodic solutions:cs} $(A(\xi,\tau),B(\xi,\tau)) = (\bar{A} e^{iK\xi},\bar{B})$ becomes of $\mathcal{O}(\delta^{-\gamma/2})$ -- see \Cref{sec:coupled system}. Therefore, it is natural to rescale $A(\xi,\tau)$ and introduce $\tilde{A}(\xi,\tau)$ by $\tilde{A}(\xi,\tau) = A(\xi,\tau) \delta^{\gamma/2}$. For $\gamma = \frac12$ -- a significant degeneration in the terminology of \cite{eckhaus1979asymptotic} -- this yields as (leading order) $\tilde{A}$B-system: 
\[
\begin{cases}
\tilde{A}_\tau = 4 \tilde{A}_{\xi \xi} + \tilde{A} - 2 \tilde{A} B \\
B_{\tau} = B_{\xi \xi} + \frac{1}{4} - r - B^2 + 2 \tilde{\eta} |\tilde{A}|^2 + 4 i(\tilde{A} \bar{\tilde{A}}_{\xi} - \bar{\tilde{A}} \tilde{A}_{\xi}),
\end{cases}
\]
where the new $i(\tilde{A} \bar{\tilde{A}}_{\xi} - \bar{\tilde{A}} \tilde{A}_{\xi})$ term has appeared from the originally higher-order corrections to the leading order AB-system \eqref{eq:coupled system:a model}. However, when the magnitude of the $|A|^2$ term decreases further -- at a next significant degeneration $\gamma = 1$ -- the leading system will change into:
\[
\begin{cases}
\tilde{A}_\tau = 4 \tilde{A}_{\xi \xi} + \tilde{A} - 2 \tilde{A} B - 3 \tilde{A}|\tilde{A}|^2 \\
B_{\tau} = B_{\xi \xi} + \frac{1}{4} - r - B^2 + 2 \tilde{\eta} |\tilde{A}|^2 +  \frac{4 i}{\sqrt{\delta}}(\tilde{A} \bar{\tilde{A}}_{\xi} - \bar{\tilde{A}} \tilde{A}_{\xi}) - 6B|\tilde{A}|^2 + 8 |\tilde{A}_\xi|^2 - 2(\tilde{A} \bar{\tilde{A}}_{\xi \xi} + \bar{\tilde{A}} \tilde{A}_{\xi \xi}).
\end{cases}
\]
Note that this situation is similar to that encountered in \cite{eckhaus1989strong}, in which the transition between sub- and supercritical (Turing-type) bifurcations in the complex Ginzburg-Landau equation was studied: one `special' term, $(\tilde{A} \bar{\tilde{A}}_{\xi} - \bar{\tilde{A}} \tilde{A}_{\xi})$, is (asymptotically) large with respect to the `usual' terms and the system thus evolves on two timescales. Note that it may be natural to rescale $\xi$ and introduce $\tilde{\xi} = \sqrt{\delta} \xi$ so that the $i(\tilde{A} \bar{\tilde{A}}_{\xi} - \bar{\tilde{A}} \tilde{A}_{\xi})$ term again becomes of $\mathcal{O}(1)$ magnitude. However, this rescaling of $\xi$ will also introduce an extra factor $\delta$ in front of the diffusion terms $\tilde{A}_{\tilde{\xi} \tilde{\xi}}$ and $B_{\tilde{\xi} \tilde{\xi}}$ (and of the terms $|\tilde{A}_{\tilde{\xi}}|^2$ and $(\tilde{A} \bar{\tilde{A}}_{\tilde{\xi} \tilde{\xi}} + \bar{\tilde{A}} \tilde{A}_{\tilde{\xi} \tilde{\xi}})$), so that the leading order $\tilde{A}$B-system will no longer be of reaction-diffusion type. Understanding the mechanisms driving this (complex and relevant) bifurcation will be the subject of future research. 
\end{remark}

\subsection{A general scalar higher-order model }\label{subsec:general setting:se}
To serve as a bridge between the (relatively) straightforward analysis of \Cref{subsec:a model} and the technically more involved analysis presented in \Cref{sec:reaction diffusion system}, we set up a Turing-fold analysis in a general class of scalar, higher-order models. Specifically, we consider the following scalar $2m$-th order equation -- that is a simplified version of the more general model \eqref{eq:illustration pde:rtab-Intro} (see \Cref{subsec:rescaling the ab system:gs:s}),
\begin{equation}\label{eq:general scalar eq:gs}
\partial_t U = F(U;\mu) + \sum_{j=1}^m (a_j - \nu \tilde{a}j) \partial_x^{2j} U + \sum_{j=1}^{\lfloor \frac{m-1}{2} \rfloor} \sum_{l=j}^{\lfloor \frac{m-1}{2} \rfloor} b_{jl} \partial_x^{2j} U \partial_x^{2l} U,
\end{equation}
where $U(x,t): \mathbb{R} \times \mathbb{R}^+ \to \mathbb{R}$, $\mu, \nu \in \mathbb{R}$ are the bifurcation parameters, $2 < m \in \mathbb{N}$, and $F: \mathbb{R}^2 \to \mathbb{R}$ is sufficiently smooth. The coefficients $a_j$, $\tilde{a}_j$, and $b_{jl}$ are real-valued. We assume that the problem is well-posed and therefore restrict our attention to values of $\nu \in I \subset \mathbb{R}$, where $I$ is an open interval, such that $a_m - \nu \tilde{a}_m > 0$ for odd $m$, and $a_m - \nu \tilde{a}_m < 0$ for even $m$.

Note that scalar $2m$-th order equations of this type -- excluding the nonlinear terms involving spatial derivatives such as $\partial_x^{2j} U$ (i.e., with $b_{jl} \equiv 0$) -- appear in (higher-order) phase-field models \cite{caginalp1986higher,gardner1990traveling}. Additionally, observe that, for simplicity, we have assumed that these coefficients do not vary with the main parameters $\mu$ and $\nu$ -- we will see (and explain) in \Cref{subsec:rescaling the ab system:gs:s} that this affects the nature of the coefficients in the governing AB-systems.

Similar to example \eqref{eq:scalar example}, the nonlinear term $F$ depends solely on $\mu$, while the second bifurcation parameter $\nu$ appears only in the linear spatial derivative terms. This choice has been made for simplicity; more general scenarios can be treated in a similar, albeit notationally more complex, manner. The same applies to our choice of letting \eqref{eq:general scalar eq:gs} depend only linearly on $\nu$.

The spatially homogeneous ODE associated to \eqref{eq:general scalar eq:gs} is given by
\begin{equation}
\label{genscalarODE}
    \dot{u} = F(u; \mu).
\end{equation}
We assume that a saddle-node bifurcation occurs as $\mu$ passes through a value $\mu^* \in \mathbb{R}$, resulting in the creation/annihilation of two equilibrium solutions, $u^+(\mu)$ and $u^-(\mu)$, with $u^+(\mu) > u^* > u^-(\mu)$, where $u^* \coloneqq u^\pm(\mu^*)$. This (completely standard) saddle-node bifurcation (see \cite{homburg2024bifurcation}) is characterized by the conditions
\begin{equation}
\label{scalar-SNcondition}
F(u^*; \mu^*) = 0 \quad \text{and} \quad F_{u}(u^*; \mu^*) = 0.
\end{equation}
Additionally, we assume that the critical points $u^\pm(\mu)$ emerge as $\mu$ increases through $\mu^*$, and that $u^+$ is stable as critical point of \eqref{genscalarODE} (i.e., that $F_{u}(u^+; \mu) < 0$ for $\mu > \mu^*$). The (signed) non-degeneracy conditions for this specific (standard) saddle-node bifurcation are given by \cite{homburg2024bifurcation}:
\begin{equation}\label{eq:non-degeneracy_conditions:gs}
    F_{uu}(u^*;\mu^*) < 0 \quad \text{ and } \quad F_{\mu}(u^*;\mu^*) > 0.
\end{equation}

The spectral stability of the homogeneous solution $u^+$, considered as a solution of \eqref{eq:general scalar eq:gs}, is determined by
\begin{equation} \label{eq:eigenvalue_problem_scalar_system}
    \omega(k; \mu, \nu) = F_u(u^+; \mu) + \sum_{j=1}^m (a_j-\nu\tilde{a}_j)(-k^2)^j := F_u(u^+; \mu) + G(k; \nu),
\end{equation}
where we note that $\omega(k; \mu, \nu) \in \mathbb{R}$, since \eqref{eq:general scalar eq:gs} includes only even spatial derivatives. Next, we assume the existence of $\nu^* \in I$ and $\bar{\nu} \in I$, with $\bar{\nu} < \nu^*$, such that for all $\nu \in (\bar{\nu}, \nu^*)$, the state $u^+(\mu)$ undergoes a (non-degenerate) Turing bifurcation with critical wave number $k^c(\nu) > 0$, as $\mu$ decreases through a critical value $\mu^t(\nu)$ with $\mu^t(\nu) > \mu^*$. As usual, the pair $(k^c(\nu), \mu^t(\nu))$ is determined by the conditions
\begin{equation}
\label{genscalarTuring}    
\omega(k^c;\mu^t,\nu) = 0, \quad \omega_k(k^c;\mu^t,\nu) = 0, \quad \text{ and } \quad \omega(k; \mu^t(\nu), \nu) < 0 \quad \forall k \in \mathbb{R}\backslash\{-k^c,k^c\}.
\end{equation}

The non-degeneracy conditions for this specific Turing bifurcation are 
\begin{equation}\label{eq:non-deg conditions turing:gs:se}
    \omega_{kk}(k^c;\mu^t,\nu) < 0, \quad \omega_{\mu}(k^c;\mu^t,\nu) < 0, \quad \text{ and } \quad  \omega_{\nu}(k^c;\mu^t,\nu) < 0,
\end{equation}
where we emphasize that the sign conditions on $\omega_\mu(k^c;\mu^t,\nu^t)$ and $\omega_{\nu}(k^c;\mu^t,\nu)$ reflect our choice to mimic the situation of example \eqref{eq:scalar example}: the Turing bifurcation occurs for $\nu < \nu^*$ (hence $\omega_\nu(k^c;\mu^t,\nu^t) < 0$) as $\mu$ decreases through $\mu^t > \mu^*$ (so $\omega_\mu(k^c;\mu^t,\nu)< 0$). 

Finally, we assume that a Turing-fold bifurcation occurs for $\nu = \nu^*$, i.e., when $\mu^t(\nu^*) = \mu^*$. At this point, the Turing and saddle-node bifurcations coincide, and the non-degeneracy conditions specified in \Cref{subsubsec:non-degeneracy conditions:gs} are satisfied.

Again, to simplify notation, we mark the functions $F$, $G$, $\omega$, and their derivatives, with a superscript $*$ to indicate evaluation at $(u^*, \mu^t, \nu^*, k^*)$, and with a superscript $t$ to indicate evaluation at $(u^t, \mu^t, \nu, k^c)$, where $u^t \coloneqq u^+(\mu^t)$. For example, we write:
\[
    F^* \coloneqq F(u^*;\mu^*), \quad \omega_{\mu}^t \coloneqq \omega_{\mu}(k^c;\mu^t,\nu), \quad \text{ and } \quad G^t_\nu \coloneqq G_\nu(k^c;\nu).
\]

\begin{remark}
\rm  
Note that $\omega_{\mu}^t < 0$ follows from our assumption that the saddle-node bifurcation occurs as $\mu$ increases through $\mu^t$, implying that $F_{\mu}^* F_{uu}^* < 0$. Specifically, if the saddle-node bifurcation occurs as $\mu$ increases, then on the ODE-stable branch (in our case $u^+$ since $F_{uu}^* < 0$), we have

\[
\frac{d}{d\mu} F_u(u^+; \mu) = \omega_{\mu}(k; \mu, \nu) < 0,
\]

because $F_{u}^* = 0$ and $F_u(u^+; \mu) < 0$ for $0 < \mu - \mu^* \ll 1$. Therefore, $\omega_\mu^t < 0$ if the saddle-node bifurcation occurs as $\mu$ increases through $\mu^*$. Conversely, $\omega_\mu^t > 0$ if the bifurcation occurs as $\mu$ decreases through $\mu^*$.
\end{remark}

\subsubsection{Non-degeneracy conditions for the Turing-fold bifurcation}\label{subsubsec:non-degeneracy conditions:gs}
\begin{figure}[t]
    \centering
    \includegraphics[width=1.0\linewidth]{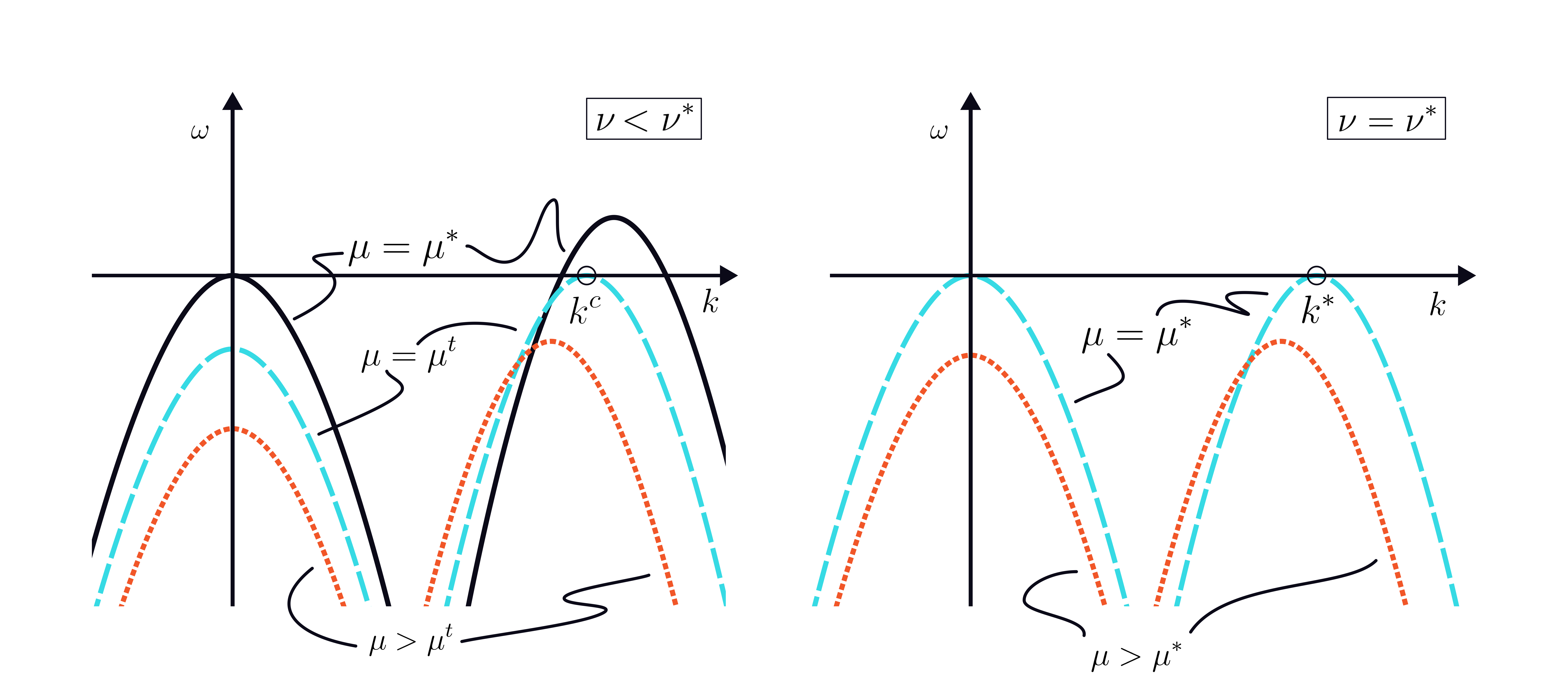}
    \caption{Illustrations of the $\omega(k; \mu, \nu)$ curves \eqref{eq:eigenvalue_problem_scalar_system} as assumed by non-degeneracy conditions \eqref{genscalar-allnondegconds} for various choices of $\nu$ and $\mu$ near the co-dimension 2 point $(\mu^*,\nu^*)$. The left panel represents the case $\nu < \nu^*$ in which the Turing bifurcation occurs just before the fold bifurcation: the $\omega$-curve crosses the $k$-axis at a critical value $k = k^c$ for $\mu < \mu^t$ while $\omega < 0$ near $k = 0$. The right panel depicts the case $\nu = \nu^*$. At the co-dimension 2 point, the $\omega(k; \mu^*, \nu^*$ curve touches the $k$-axis simultanuously at $k = 0$ and $k=k^* > 0$.} 
    \label{fig:dispersion relation:ndc:gs}
\end{figure}

Naturally, it is necessary to formulate and impose non-degeneracy conditions on the co-dimension 2 Turing-fold bifurcation. We begin by stating the conditions required for the occurrence of this bifurcation:
\begin{equation}
\label{genscalar-TuringFold}    
\begin{array}{rrrr}
{\rm (Fold)}_{\rm bif}:  &  F^* = 0, & F_{u}^* = 0; & \\
{\rm (Turing)}_{\rm bif}:  &  \omega^* = F_u^* + G^* = G^* = 0, & \omega_k^* = G_k^* = 0, & \omega(k; \mu^*, \nu^*) < 0 \;\, \forall k \in \mathbb{R}\backslash\{-k^*,0,k^*\}
\end{array}
\end{equation}
for some $k^* > 0$, as per the eigenvalue problem \eqref{eq:eigenvalue_problem_scalar_system}. This bifurcation is non-degenerate if the following conditions hold:
\begin{equation}
\label{genscalar-allnondegconds}    
\begin{array}{rrrr}
{\rm (Fold)}_{\rm n-d}:  &  F_{uu}^* < 0, & F_{\mu}^* > 0, & \rho^*_{kk} := \omega_{kk}(0; \mu^*, \nu^*) = G_{kk}(0; \nu^*) < 0;\\
{\rm (Turing)}_{\rm n-d}:  &  \omega_{kk}^* = G_{kk}^* < 0, & \omega_{\nu}^* = G_{\nu}^* < 0. &
\end{array}
\end{equation} 
We note (again) that the specific sign conditions on $F_{uu}^*$, $F_{\mu}^*$, and $\omega_{\nu}^*$ are chosen to match the setup of example equation \eqref{eq:scalar example}. \Cref{fig:dispersion relation:ndc:gs} shows sketches of the eigenvalue curve $\omega(k; \mu, \nu)$ for several values of $\mu$ and $\nu$ near $(\mu^*,\nu^*)$. 

As in the case of the example problem, conditions \eqref{genscalar-TuringFold}  and \eqref{genscalar-allnondegconds} together imply the existence of a curve of saddle-node points, $\Gamma^s = \{(\mu^*, \nu): \nu \in I\}$, and a curve of Turing points, $\Gamma^t = \{(\mu^t(\nu), \nu), \nu \in (\bar{\nu},\nu^*)\}$, that meet tangentially at the co-dimension 2 point $(\mu^*, \nu^*) \in \Gamma^s \cap \Gamma^t$; see \Cref{fig:pde bifurcation diagram scalar example:a model}. In fact, for $\nu = \nu^* - \delta$, we have
\begin{equation}\label{eq:mut, kc, ut form:ndc:gs:se}
    \mu^t = \mu^* + \hat{\mu} \delta^2 + \mathcal{O}(\delta^3), \quad u^t = u^* + \tilde{u} \delta + \mathcal{O}(\delta^2), \quad k^c = k^* + \tilde{k} \delta + \mathcal{O}(\delta^2), \quad \omega_{\mu}^t = \frac{\tilde{\omega} + \mathcal{O}(\delta)}{\delta},
\end{equation}
which were highlighted in \Cref{subsubsec:basic analysis:a model} as the crucial properties of the Turing and fold bifurcations near a non-degenerated co-dimension 2 Turing-fold point (cf. \eqref{eq: scalar example PDE}, \eqref{wmutO1delta}, and \eqref{kcuutmut}). More specifically, the non-degeneracy conditions \eqref{genscalar-allnondegconds} imply
\begin{equation}\label{eq:tilde mu u and omega:ndc:gs}
    \hat{\mu} = -\frac{(G_\nu^*)^2}{2 F_{uu}^* F_{\mu}^*} > 0, \quad   
    \tilde{u} = \frac{G_\nu^*}{F_{uu}^*} > 0, \quad
    \tilde{k} = \frac{G_{k\nu}^*}{G_{kk}^*}, \quad \tilde{\omega} = - \frac{F_{uu}^* F_{\mu}^*}{G_\nu^*} < 0,
\end{equation}
which follow from direct expansion and the conditions imposed in \eqref{genscalar-TuringFold}, and \eqref{genscalar-allnondegconds}. Namely, from
\[
\begin{array}{rcl}
0 & = & F(u^t; \mu^t) \\
& = & F(u^* + (u^t-u^*); \mu^* + (\mu^t-\mu^*)) \\
& = & F^* + F^*_u(u^t-u^*) + \frac12 F^*_{uu}(u^t-u^*)^2 + F^*_\mu(\mu^t-\mu^*) + \text{h.o.t.}     
\\
& = & \frac12 F^*_{uu}(u^t-u^*)^2 + F^*_\mu(\mu^t-\mu^*) + \text{h.o.t.,}      
\end{array}
\]
we conclude that $\mu^t-\mu^* = \mathcal{O}((u^t-u^*)^2)$. Moreover, this expansion, together with
\[
\begin{array}{rcl}
0 & = & \omega(k^c; \mu^t, \nu^* - \delta) \\
& = & F_u(u^* + (u^t-u^*); \mu^* + (\mu^t-\mu^*)) + G(k^* + (k^c-k^*); \nu^* - \delta) \\
& = & F_u^* + F^*_{uu}(u^t-u^*) + G^* + G_k^*(k^c-k^*) - G^*_\nu \delta + \text{h.o.t.}     
\\
& = & F^*_{uu}(u^t-u^*) - G^*_\nu \delta + \text{h.o.t.}, 
\end{array}
\]
yield:
\[
u^t = u^* + \frac{G_\nu^*}{F_{uu}^*} \delta + \mathcal{O}(\delta^2), \quad \text{ and } \quad \mu^t = \mu^* -\frac{(G_\nu^*)^2}{2 F^*_\mu F_{uu}^*} \delta^2 + \mathcal{O}(\delta^3).
\]
Next, we note that $\tilde{k}$ is determined from
\[
\begin{array}{rcl}
0 & = & \omega_k(k^c; \mu^t, \nu^* - \delta) \\
& = & G_k(k^* + (k^c-k^*); \nu^* - \delta) \\
& = & G_k^* + G_{kk}^*(k^c-k^*) - G_{k \nu}^* \delta + \text{h.o.t.} \\
& = & G_{kk}^*(k^c-k^*) - G_{k \nu}^* \delta + \text{h.o.t.,}   
\end{array}
\]
which leads directly to $\tilde{k} = G_{k\nu}^*/G_{kk}^*$.
Finally, using \eqref{genscalar-TuringFold}, \eqref{eq:mut, kc, ut form:ndc:gs:se}, and \eqref{eq:tilde mu u and omega:ndc:gs}, we obtain
\begin{equation}
\label{Fut}    
F_u(u^t; \mu^t) = F_u^* + F_{uu}^*(u^t - u^*) + \text{h.o.t.} = G^*_\nu \delta  + \mathcal{O}(\delta^2),
\end{equation}
(cf. \eqref{Fuu+mut}), which in combination with
\[
\omega_{\mu}^t = F_{uu}^t u_{\mu}^t + F_{u\mu}^t,
\quad {\rm and} \quad
0 = \frac{d}{d \mu} F(u^+(\mu); \mu)|_{\mu = \mu^t} = F_u^t u^t_\mu + F_\mu^t,
\]
yields
\[
\omega_{\mu}^t = - F_{uu}^t \frac{F_\mu^t}{F_u^t} + F_{u\mu}^t.
\]
Therefore, using \eqref{Fut}, we recover the leading order approximation $\tilde{\omega} = - (F_{uu}^* F_{\mu}^*)/(G_\nu^*)$, where $\tilde{\omega}$ is defined in \eqref{eq:mut, kc, ut form:ndc:gs:se}.
 
We must verify that the bifurcation and non-degeneracy conditions \eqref{genscalar-TuringFold} and \eqref{genscalar-allnondegconds} for the co-dimension 2 Turing-fold bifurcation are consistent with those for the saddle-node bifurcation \eqref{scalar-SNcondition} and \eqref{eq:non-degeneracy_conditions:gs}, as well as the Turing bifurcation conditions \eqref{genscalarTuring} and \eqref{eq:non-deg conditions turing:gs:se}. Apart from the global conditions on $\omega(k; \mu, \nu)$ in \eqref{eq:non-deg conditions turing:gs:se} and \eqref{genscalar-allnondegconds}, this consistency follows immediately from the assumed smoothness of the involved expressions. Moreover, together with the local control over $\omega(k)$, which we obtained above, near its (non-degenerate) maximums at $k^c$ and $k^*$, the same smoothness argument can be applied to conclude that these global conditions match up on compact $k$-intervals. Since \eqref{eq:general scalar eq:gs} is well-posed, we know that $\omega(k)$ decays as $-|a_m - \nu \tilde{a}_m| k^{2m}$ for large $|k|$ (with $|a_m - \nu \tilde{a}_m| > 0$). Therefore we conclude that, near $(\mu^*,\nu^*)$, all existence and non-degeneracy conditions of the 2 co-dimension 1 bifurcations are in agreement with those of the co-dimension 2 bifurcation.

Finally, we highlight that the non-degeneracy condition $\rho^*_{kk} = \omega_{kk}(0; \mu^*, \nu^*) < 0$ in \eqref{genscalar-allnondegconds} is a new condition, which we did not impose on the underlying saddle-node bifurcation. This condition does not appear in \eqref{eq:non-degeneracy_conditions:gs} since there the saddle-node bifurcation was introduced as a bifurcation in the spatially homogeneous ODE counterpart \eqref{genscalarODE} of PDE \eqref{eq:general scalar eq:gs}. The condition $\rho^*_{kk} < 0$ clearly pertains to the spatial (linear, long-wavelength) behavior of solutions to the PDE \eqref{eq:general scalar eq:gs} at the saddle-node bifurcation. 

From the PDE point of view, it would have been natural to add the condition $\rho^t_{kk} = \omega_{kk}(0; \mu^t, \nu) < 0$ to \eqref{eq:non-degeneracy_conditions:gs}, aligning it with $\rho^*_{kk} < 0$ of \eqref{genscalar-allnondegconds}, to indicating that $\omega(k; \mu^t, \nu)$ has a maximum at $k=0$. Due to the $k \to -k$ symmetry of \eqref{eq:eigenvalue_problem_scalar_system}, we have $\omega_k(0; \mu, \nu) = 0$ for all $\mu, \nu$, implying that $\omega(k;\mu^t, \nu)$ must have an extreme at $k=0$. Therefore, imposing $\rho^t_{kk} < 0$ guarantees that $k = 0$ is a non-degenerate maximum rather than a minimum, in which case the state $u^+(\mu)$ would have been unstable with respect to long wavelength perturbations. This is precisely the desired setting for a (non-degenerate) saddle-node bifurcation (of homogeneous states) in a PDE. 

Finally, since two extrema of a parametrized family of smooth curves can only merge with a vanishing second derivative, we note that the natural assumption that $\omega(k;\mu^t(\nu), \nu)$ has a non-degenerate maximum at $k=0$, which persists in the co-dimension 2 limit $\nu \to \nu^*$, implies that $k^c(\nu)$ cannot approach $0$ as $\nu \to \nu^*$: $k^c(\nu) \to k^*$ with necessarily $k^* > 0$. In other words, the non-degeneracy condition $\rho^*_{kk} = \omega_{kk}(0; \mu^*, \nu^*) < 0$ forces $k^*$ to be strictly positive, eliminating the need to assume it as an independent condition (which we did in \eqref{genscalar-TuringFold} and now proves to be redundant).  

\begin{remark}\label{rem:sixth-order}
\rm  
    The fact that both $k = \pm k^* \neq 0$ and $k = 0$ are absolute maxima of the polynomial $\omega(k; \mu^*, \nu^*)$, implies that a co-dimension 2 Turing-fold bifurcation can only occur in a scalar PDE of the form \eqref{eq:general scalar eq:gs} if the equation is at least sixth order.
\end{remark}

\subsubsection{The Ginzburg-Landau approximation}\label{subsubsec:ginzburg-landau derivation:gs}

Analogous to the approximation procedure carried out in \Cref{subsubsec:ginzburg-landau derivation:a model} for the example model \eqref{eq:scalar example}, we begin the weakly nonlinear analysis of small-amplitude patterns by deriving the classical Ginzburg-Landau approximation for $\mu$ sufficiently close to $\mu^t(\nu)$ with $\nu = \nu^* - \delta = -\delta$ (assuming, for simplicity and without loss of generality, that $\nu^* = 0$). Accordingly, we set $\mu = \mu^t(-\delta) - r \varepsilon^2$ with $r \in \mathbb{R}$ and $0 < \varepsilon \ll \delta \ll 1$. Additionally, we define
\begin{equation*}
    P(k_1, k_2) \coloneqq 2 \sum_{j=1}^{\lfloor\frac{m-1}{2} \rfloor} \sum_{l=j}^{\lfloor \frac{m-1}{2} \rfloor} b_{jl} (-k_1^2)^j(-k_2^2)^l,
\end{equation*}
and introduce the shorthand notation
\begin{equation}
\label{def:GjP}
    P^* \coloneqq P(k^*, k^*) \quad \text{ and } \quad P^t \coloneqq P(k^c, k^c).
\end{equation}
Following the standard procedure, we the substitute Ginzburg-Landau Ansatz \eqref{eq:scalar ginzburg-landau ansatz:ginzburg-landau derivation} into \eqref{eq:general scalar eq:gs} and expand:
\begin{equation}
    \begin{split}
        &F(U_{GL};\mu^t - \varepsilon^2 r) + \sum_{j=1}^m (a_j + \delta \tilde{a}_j) \partial_x^{2j} U_{GL} + \sum_{j=1}^{\lfloor \frac{m-1}{2} \rfloor} \sum_{l=j}^{\lfloor \frac{m-1}{2}\rfloor} b_{jl} \left(\partial_x^{2j} U_{GL}\right) \left(\partial_x^{2l} U_{GL}\right) = \\
        &F^t + F^t_u U_{GL} + \frac{1}{2}F^t_{uu} U_{GL}^2+ \frac{1}{6} F^t_{uuu}U_{GL}^3 - r \varepsilon^2 F_{\mu}^t - r \varepsilon^2 U_{GL} F_{u\mu}^t +
        \\ &\sum\limits_{j=1}^n (a_j + \delta \tilde{a}_j) \partial_{x}^{2j} U_{GL} + \sum_{j=1}^{\lfloor \frac{m-1}{2} \rfloor} \sum_{l=j}^{\lfloor \frac{m-1}{2}\rfloor} b_{jl} (\partial_{x}^{2j} U_{GL})(\partial_{x}^{2l} U_{GL}) + \text{h.o.t.}
    \end{split}
\end{equation}
The first non-trivial results appear, as usual, at the $\varepsilon^2E^0$-level:
\[
    - F_{u}^t X_{02} = \left(F_{uu}^t + P^t\right)|A|^2 -r F_\mu^t,
\]
and at the $\varepsilon^2 E^2$-level:
\[
    -\left(F_u^t + G(2k^c;-\delta)\right)X_{22} = \frac{1}{2}\left(F_{uu}^t + P^t \right)A^2.
\]
At the $\varepsilon^3E$-level we obtain
\[
    A_\tau = - \frac{1}{2}\omega_{kk}^t A_{\xi\xi} - r F_{\mu u}^t A + F_{uu}^t X_{02}A + \left(F_{uu}^t+P(2 k^c, k^c)\right)X_{22}A^* + \frac12F_{uuu}^t|A|^2A,
\]
where $A^*$ means $A$ conjugate. Consequently, using the two expressions we found at the $\varepsilon^2$-level, this yields the Ginzburg-Landau equation
\[
    A_\tau = - \frac{1}{2}\omega_{kk}^t A_{\xi\xi} - r \omega_{\mu}^t A + L |A|^2A
\]
(at leading order in $\varepsilon$), where the Landau coefficient $L$ is given by
\begin{equation}
    L = -\frac{F_{uu}^t\left(F_{uu}^t+P^t\right)}{F_{u}^t}-\frac{\left(F_{uu}^t+P(2k^c,k^c) \right)\left(F_{uu}^t + P^t\right)}{2\left(F_u^t + G(2k^c;-\delta) \right)}+\frac{1}{2} F_{uuu}^t.
\end{equation}
Using \eqref{genscalar-allnondegconds}, \eqref{eq:mut, kc, ut form:ndc:gs:se}, \eqref{eq:tilde mu u and omega:ndc:gs}, and \eqref{Fut}, we can obtain a simplified expression for $L$ (in terms of a leading order approximation in $\delta$) and express the Ginzburg-Landau equation in a form analogous to \eqref{eq:ginzburg-landau:a model}:
\begin{equation}\label{eq:ginzburg-landau with deltas:gs}
    A_\tau = \left(-\frac{1}{2}G_{kk}^* + \mathcal{O}(\delta)\right) A_{\xi\xi} + \frac{r}{\delta} \left(\frac{F^*_{uu} F^*_{\mu}}{G^*_\nu} + \mathcal{O}(\delta)\right) A + \frac{1}{\delta} \left(\frac{-F_{uu}^*(F_{uu}^* + P^*)}{G^*_\nu}+ \mathcal{O}(\delta)\right)|A|^2A
\end{equation}
(at leading order in $\varepsilon$). 

First, we note that the Turing bifurcation is supercritical, respectively subcritical, if
\[
P^* > -F_{uu}^*, \quad {\rm resp.} \quad P^* < - F_{uu}^*,
\]
since the non-degeneracy conditions stated in \eqref{genscalar-allnondegconds} dictate that $G_\nu^* F_{uu}^* > 0$. Consequently, as $P^*$ is associated with the nonlinearities involving spatial derivatives, we conclude that the Turing bifurcation near a co-dimension 2 Turing-fold point in the class of scalar systems \eqref{eq:general scalar eq:gs} cannot be supercritical in Swift-Hohenberg-type cases where the nonlinear terms do not contain spatial derivatives (i.e., when $b_{jl} \equiv 0$). Since stable periodic patterns emerge (as a first-order effect) from the Turing bifurcation only in the supercritical case, this strongly suggests that one can only expect the persistence of stable Turing patterns beyond the saddle-node point if \eqref{eq:general scalar eq:gs} includes nonlinear terms containing spatial derivatives (which motivated our inclusion of the $\eta(\partial_x^2U)$ term in \eqref{eq:scalar example}).

Finally, following the same argument as in \Cref{subsubsec:ginzburg-landau derivation:a model}, again, all $|A|^{2n}A$ terms ($n \geq 2$), which a priori appear as $\mathcal{O}(\varepsilon^{2n-2})$ correction terms in the Ginzburg-Landau approximation, grow to the same order as $|A|^2A$ (i.e., of order $\delta^{-1}$) when we increase $\varepsilon$ to $\varepsilon = \delta$. Consequently, the Ginzburg-Landau approximation procedure remains valid as long as $\varepsilon \ll \delta$. Therefore, we can (again) not cover the saddle-node bifurcation with this approach.

\subsubsection{The derivation of the AB-system}\label{subsubsec:coupled system derivation:gs}
Following \Cref{subsec:a model}, we set $\nu = - \delta$ and $\mu = \mu^t(-\delta) - r \delta^2$ and replace the standard Ginzburg-Landau Ansatz with its modified version:
\begin{equation}\label{eq:coupled system ansatz:csd:gs}
    \begin{split}
    U_{AB}(x, t) = u^* + E^0 \big[\delta B + \delta^{3/2}X_{02} + &\delta^{2} X_{03} +\delta^{5/2} X_{04}+\mathcal{O}\l\delta^3\r\big] \\
    + E \big[\delta A + \delta^{3/2} X_{12} + &\delta^{2} X_{13} + \delta^{5/2} X_{14} +\mathcal{O}\l \delta^3 \r \big] + \text{c.c.} \\
    +E^2 \big[&\delta^2 X_{23} + \delta^{5/2} X_{24} + \mathcal{O}\l\delta^3\r\big] + \text{c.c.} \\
     &\phantom{iiiii00000000000iiiiii+E_c^3}\:  + \text{h.o.t.},
\end{split}
\end{equation}
(cf. \eqref{eq: Coupled system Ansatz}) where $E = e^{i k^* x}$, and $A(\xi,\tau):\mathbb{R}\times\mathbb{R}^+ \to \mathbb{C}$ and  $B(\xi,\tau):\mathbb{R}\times\mathbb{R}^+ \to \mathbb{R}$ are two a priori unknown amplitude functions. Substitution into \eqref{eq:general scalar eq:gs} and subsequent expansion yields
\begin{equation}
    \begin{split}
        &F(U_{AB};\mu^t - \delta^2 r) + \sum_{j=1}^m (a_j + \delta \tilde{a}_j) \partial_x^{2j} U_{AB} + \sum_{j=1}^{\lfloor \frac{m-1}{2} \rfloor} \sum_{l=j}^{\lfloor \frac{m-1}{2}\rfloor} b_{jl} \l \partial_{x}^{2j} U_{AB}\r\l\partial_{x}^{2l} U_{AB}\r = \\
        &F^* + F^*_u U_{AB} + \frac{1}{2}F^*_{uu} U_{AB}^2 + (\hat{\mu} - r) \delta^2 F_{\mu}^*  + \sum\limits_{j=1}^n (a_j + \delta \tilde{a}_j) \partial_{x}^{2j} U_{AB} +
        \\ &\sum_{j=1}^{\lfloor \frac{m-1}{2} \rfloor} \sum_{l=j}^{\lfloor \frac{m-1}{2}\rfloor} b_{jl} \l \partial_{x}^{2j} U_{AB}\r\l\partial_{x}^{2l} U_{AB}\r + \text{h.o.t.}
    \end{split}
\end{equation}
Again, as in \Cref{subsec:a model}, the $\delta$- and $\delta^{3/2}$-levels cancel by construction, and we derive the AB-sytem by collecting the terms at the $\delta^2 E^i$-levels. Specifically, the $\delta^2E^0$-level gives
\begin{equation}
    B_\tau = -\frac{1}{2}\rho_{kk}^* B_{\xi\xi} + F_{\mu}^*(\hat{\mu}-r) + \frac{1}{2}F_{uu}^* B^2 + \left(F_{uu}^*+ P^*\right)|A|^2,
\end{equation}
and the $\delta^2E$-level yields
\begin{equation}
    A_\tau = - \frac{1}{2} \omega_{kk}^* A_{\xi\xi} + \sum_{j=1}^m\tilde{a}_j (-(k^*)^2)^{j} A + F_{uu}^* AB.
\end{equation}
Using \eqref{eq:eigenvalue_problem_scalar_system}, \eqref{genscalar-allnondegconds}, \eqref{eq:mut, kc, ut form:ndc:gs:se}, and \eqref{eq:tilde mu u and omega:ndc:gs}, we conclude that, near the co-dimension 2 Turing-fold bifurcation, the evolution of small-amplitude solutions of \eqref{eq:general scalar eq:gs} is described, at leading order in $\delta$, by
\begin{equation}\label{eq:coupled system:gs}
\begin{cases}
    A_\tau = - \frac{1}{2} G_{kk}^* A_{\xi\xi} - G^*_{\nu} A + F_{uu}^* AB \\
    B_\tau = -\frac{1}{2}\rho_{kk}^* B_{\xi\xi} + F_{\mu}^*\l-\frac{(G_\nu^*)^2}{2 F_{uu}^* F_{\mu}^*}-r\r + \frac{1}{2}F_{uu}^* B^2 + \left(F_{uu}^*+ P^*\right)|A|^2,
\end{cases}
\end{equation}
which has the same structure as the AB-system \eqref{eq:coupled system:a model} associated with the example model \eqref{eq:scalar example}. 

This AB-sytem covers a parameter region in $r$, which in essence is the scaled version of $\mu$, that includes both the Turing bifurcation (by construction) and the saddle-node bifurcation. To see this, consider the spatially homogeneous ODE associated with \eqref{eq:coupled system:gs}:
\begin{equation}
\label{eq:coupled system:gs:ODE}
\begin{cases}
    \dot{a} = - G^*_{\nu} a + F_{uu}^* ab \\
    \dot{b} = F_{\mu}^*\l-\frac{(G_\nu^*)^2}{2F_{uu}^* F_{\mu}^*}-r\r + \frac{1}{2}F_{uu}^* b^2 + \left(F_{uu}^*+ P^*\right)|a|^2,
\end{cases}
\end{equation}
which has critical points
\begin{equation}
\label{AsrBsr:gs}
\left(A_s(r),B_s^\pm(r)\right) = \left(0, \pm\sqrt{\frac{F_{\mu}^*}{F^*_{uu}}\left(\frac{(G_\nu^*)^2}{F_{uu}^* F_{\mu}^*} + 2r\right)}\right),
\end{equation}
corresponding directly to the $u^\pm(\mu)$ of \eqref{eq:general scalar eq:gs} (cf. \eqref{eq: CS saddle-node solution}). These critical points vanish in a saddle-node bifurcation as $r$ increases through
\begin{equation}
\label{rs-unsc}
r^s = -\frac{(G_\nu^*)^2}{2F_{uu}^* F_{\mu}^*} = \hat{\mu} > 0
\end{equation}
\eqref{eq:tilde mu u and omega:ndc:gs}, so that $\mu^* = \mu^* + \hat{\mu} \delta^2 - \hat{\mu} \delta^2 = \mu^t - r^s \delta^2 + \mathcal{O}(\delta^3)$ \eqref{eq:mut, kc, ut form:ndc:gs:se}.

\subsubsection{Recovering the Ginzburg-Landau equation by zooming in from the AB-sytem}
\label{subsubsec:Reducing the AB-system:gs:se}

\begin{figure}[t]
    \centering
    \includegraphics[width=0.7\linewidth]{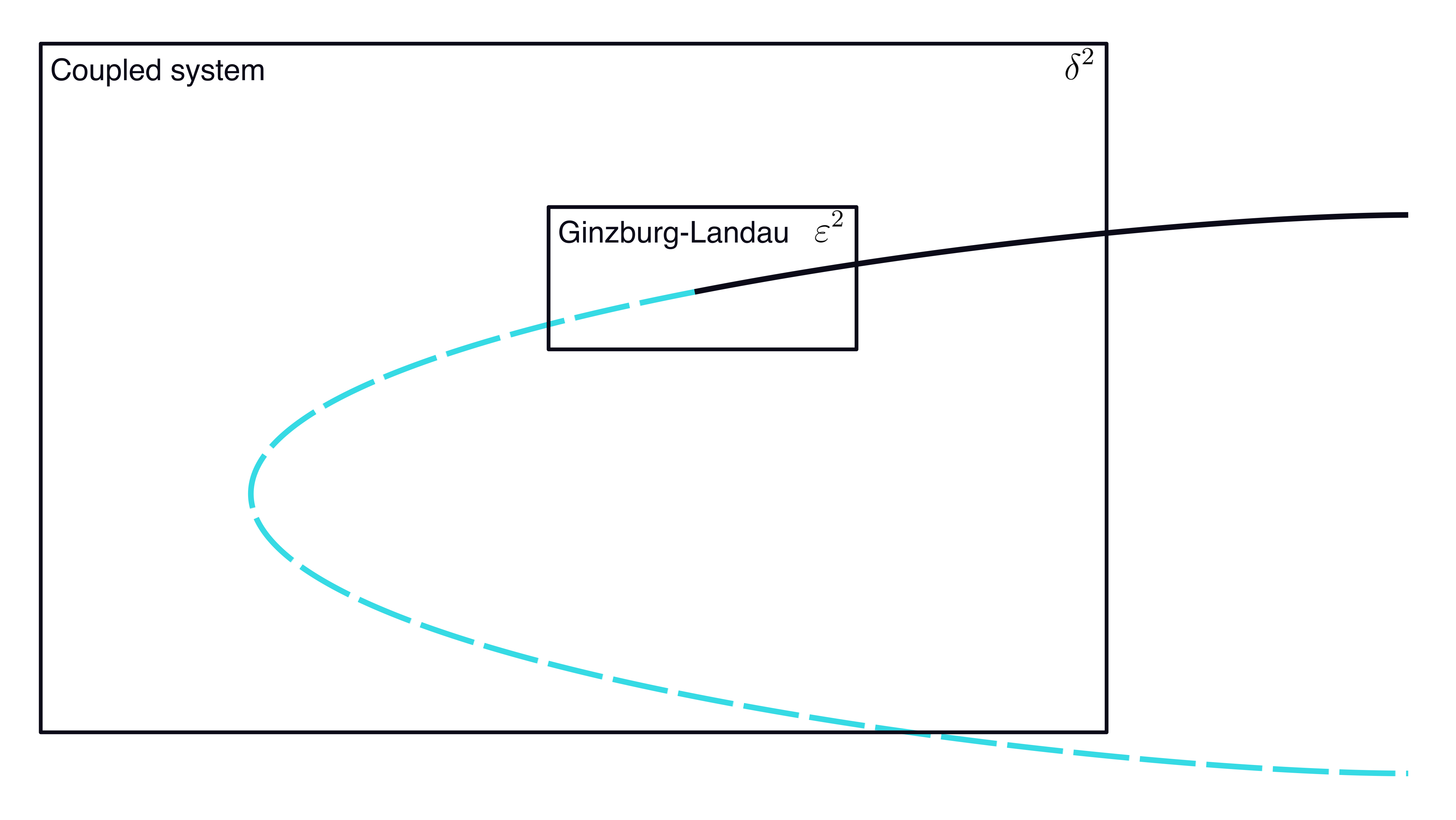}
    \caption{A sketch of the parameter region of validity of the coupled AB-sytem vs. the Ginzburg-Landau equation: in the Turing-fold setting, the coupled AB-sytem is, in essence, a `wider lens' that describes the nonlinear evolution of solutions of \eqref{eq:general scalar eq:gs} with small, $\mathcal{O}(\delta)$, amplitude for $\mu$ $\mathcal{O}(\delta^2)$ close to the Turing bifurcation, while the standard Ginzburg-Landau approach zooms in $\mathcal{O}(\varepsilon^2)$ close to Turing bifurcation and governs the dynamics solutions of solutions of $\mathcal{O}(\varepsilon)$ amplitude -- and $0 < \varepsilon \ll \delta \ll 1$.}
    \label{fig:gl vs cs validity:gs}
\end{figure}

As a system of modulation equations, the AB-system is, in essence, a `wider lens' to describe the behavior near the Turing bifurcation compared to the Ginzburg-Landau equation, since it covers a larger region in bifurcation parameter $\mu$ space than the Ginzburg-Landau equation, as illustrated in \Cref{fig:gl vs cs validity:gs}. In the region where the Ginzburg-Landau equation is valid, both the AB-system and the Ginzburg-Landau equation must describe the same behavior, implying that we should be able to reduce the AB-system to the Ginzburg-Landau equation by redoing and subsequently undoing the scalings by which we derived both modulation equations, in other (less formal) words, by zooming in from the AB-scalings into the Ginzburg-Landau scalings. 

The Ginzburg-Landau approximation is defined for solutions $\mathcal{O}(\varepsilon)$ near the homogeneous state $u^t$ of \eqref{eq:general scalar eq:gs}, with $0 < \varepsilon \ll \delta \ll 1$ (cf. \eqref{eq:scalar ginzburg-landau ansatz:ginzburg-landau derivation}). Thus, we set
\[
B = B_s^+(0) + \varepsilon^2 \tilde{B} = \sqrt{\frac{F_{\mu}^*}{F^*_{uu}} \frac{(G_\nu^*)^2}{F_{uu}^* F_{\mu}^*}} + \varepsilon^2 \tilde{B} = \frac{G_\nu^*}{F_{uu}} + \varepsilon^2 \tilde{B} = \tilde{u} + \varepsilon^2 \tilde{B},
\]
\eqref{AsrBsr:gs}, \eqref{eq:tilde mu u and omega:ndc:gs}, which reconfirms that $u^t = u^* + \tilde{u} \delta +\mathcal{O}(\delta^2)$. Introducing the Ginzburg-Landau scalings
\[
\tilde{A} = \frac{1}{\varepsilon} A, \quad \tilde{\xi} = \varepsilon \xi, \quad \tilde{\tau} = \varepsilon^2 \tau, \quad \text{ and } \quad \tilde{r} = \frac{1}{\varepsilon^2} r,
\]
transforms \eqref{eq:coupled system:gs} into 
\[
\begin{cases}
    \tilde{A}_{\tilde{\tau}} = - \frac{1}{2} G_{kk}^* \tilde{A}_{\tilde{\xi} \tilde{\xi}} + F_{uu}^* \tilde{A}\tilde{B} + \mathcal{O}\l\delta^{1/2},\varepsilon\r, \\
    0 = -F_{\mu} \tilde{r} + G^*_\nu \tilde{B} + \left(F_{uu}^*+ P^*\right)|\tilde{A}|^2 + \mathcal{O}\l\delta^{1/2},\varepsilon\r.
\end{cases}
\]
Solving the second equation for $\tilde{B}$ gives
\[
\tilde{B} = \frac{F_{\mu} \tilde{r} - \left(F_{uu}^*+ P^*\right)|\tilde{A}|^2}{G_\nu^*} + \mathcal{O}\l\delta^{1/2},\varepsilon\r, 
\]
and subsequent substitution into the first equation yields
\[
\tilde{A}_{\tilde{\tau}} = - \frac{1}{2} G_{kk}^* \tilde{A}_{\tilde{\xi} \tilde{\xi}} +
\tilde{r} \frac{F^*_{uu} F^*_{\mu}}{G^*_\nu} \tilde{A} - \frac{F_{uu}^*(F_{uu}^* + P^*)}{G^*_\nu} |\tilde{A}|^2 \tilde{A} + \mathcal{O}\l\delta^{1/2},\varepsilon\r,
\]
which is a rescaled version of \eqref{eq:ginzburg-landau with deltas:gs}. By finally undoing the scalings associated with the derivation of the AB-sytem,
\[
\hat{A} = \delta \tilde{A}, \quad \hat{r} = \delta^2 \tilde{r}, \quad \hat{\xi} = \frac{1}{\sqrt{\delta}} \tilde{\xi}, \quad  \text{ and } \quad \hat{\tau} = \frac{1}{\delta} \tilde{\tau},
\]
we indeed recover \eqref{eq:ginzburg-landau with deltas:gs} at leading order:
\[
\hat{A}_{\hat{\tau}} = - \frac{1}{2} G_{kk}^* \hat{A}_{\hat{\xi} \hat{\xi}} +
\frac{\hat{r}}{\delta} \frac{F^*_{uu} F^*_{\mu}}{G^*_\nu} \hat{A} - \frac{1}{\delta} \frac{F_{uu}^*(F_{uu}^* + P^*)}{G^*_\nu} |\hat{A}|^2 \hat{A}.
\]
Finally, we note the coefficient of the $|A|^2$ term in the (unscaled) AB-system \eqref{eq:coupled system:gs} is $F_{uu}^* + P^*$, while the Landau coefficient of the $|A|^2A$ term in the Ginzburg-Landau equation \eqref{eq:ginzburg-landau with deltas:gs} can be written as $-F_{uu}^*/(G^*_\nu \delta)$ times this factor. Consequently, since the (non-degeneracy) conditions stated in \eqref{scalar-SNcondition} and \eqref{eq:non-deg conditions turing:gs:se} imply that $F_{uu}^*G_\nu^* > 0$, these two coefficients must have opposite signs and vanish or change sign simultaneously. In \Cref{sec:reaction diffusion system}, we shall recover $-F_{uu}^*/(G^*_\nu \delta)$ as $\tilde{\omega}/(F_\mu^* \delta)$ (see \eqref{eq:mut, kc, ut form:ndc:gs:se}) and show that this relation also holds for $n$-component reaction-diffusion equations. In fact, as we will discuss in \Cref{sec:discussion}, this relationship is a universal feature of the Turing-fold bifurcation. In other words, the fact that $|A|^2$ coefficient of the AB-system and the Landau coefficient of the Ginzburg-Landau approximation always have opposite signs (if non-zero) is a universal feature of the Turing-fold bifurcation. 

\subsection{Rescaling the AB-system}\label{subsec:rescaling the ab system:gs:s}
To simplify the analysis of AB-systems associated with Turing-fold bifurcations, we normalize their coefficients to obtain a canonical form. While multiple canonical forms are possible, we adopt the form \eqref{eq:canonical form:csd:gs} as given in the Introduction. Note that this form is not as simplified as it could be, since the rescaling $A=|\beta|^{-1/2}\tilde{A}$ would transform the $\beta|A|^2$ term in \eqref{eq:canonical form:csd:gs} into $\text{sign}(\beta)|\tilde{A}|^2$ . However, this alternative form complicates the (presentation of the) analysis of \Cref{sec:coupled system}, so we prefer \eqref{eq:canonical form:csd:gs} over the $\pm |A|^2$-version. 

System \eqref{eq:coupled system:gs} can be transformed into \eqref{eq:canonical form:csd:gs} using the rescalings
\[
    A = \tilde{A}, \quad B = \frac{G^*_\nu}{F^*_{uu}} \tilde{B}, \quad R = -\frac{2 F_{uu}^* F_{\mu}^*}{(G_\nu^*)^2} r,\quad \tilde{\tau} = -G^*_\nu \tau, \quad \text{ and } \quad \tilde{\xi} = \sqrt{\frac{2 G^*_{\nu}}{G_{kk}^*}}\xi,
\]
so that
\begin{equation}
\label{eq:rescaleABalphaetc}    
    \alpha = \frac{1}{2} > 0, \quad d = \frac{2\rho^*_{kk}}{G^*_{kk}} > 0, \quad \text{ and } \quad \beta = -\frac{2 (F_{uu}^*+P^*) F_{uu}^*}{(G_\nu^*)^2}. 
\end{equation}
Since any AB-system associated with a Turing-fold bifurcation has eight unknown coefficients, $c_1, \cdots c_8$,
\begin{equation*}
    \begin{cases}
        A_\tau = c_1 A_{\xi\xi} + c_2 A + c_3 AB, \\
        B_\tau = c_4 B_{\xi\xi} + c_5 - c_6 r + c_7 B^2 + c_8|A|^2,
    \end{cases}
\; {\rm or} \quad
 \begin{cases}
        A_\tau = c_1 A_{\xi\xi} + c_2 A + c_3 AB, \\
        \frac{1}{c_5} B_\tau = \frac{c_4}{c_5} B_{\xi\xi} + 1 - \frac{c_6}{c_5} r + \frac{c_7}{c_5} B^2 + \frac{c_8}{c_5}|A|^2,
    \end{cases}    
\end{equation*}
it is not obvious that four rescaling can normalize five coefficients, i.e., that rescaling $A$, $B$, $\xi$, and $\tau$ can set $(c_1,c_2,c_3,c_6/c_5,c_7/c_5)$ to $(1,1,-1,1,-1)$ -- which is implicitly claimed in \eqref{eq:canonical form:csd:gs}. However, this is a generic property of any AB-system governing pattern formation near a Turing-fold point, as we will argue below. On the contrary, that $\alpha$ scales precisely to $\frac12$ \eqref{eq:rescaleABalphaetc} is not generic; it results from the nature of our (simplified) choice for the structure of \eqref{eq:general scalar eq:gs}, which we will also explain below.

To substantiate these claims, we consider the wider class of scalar PDEs that was already introduced in the Introduction (cf. \eqref{eq:illustration pde:rtab-Intro}),
\begin{equation}
\label{eq:illustration pde:rtab}
    \partial_t U = F(U;\mu) + \Psi(\partial_x^2 U, \dots, \partial_x^{2m} U, U;\nu),
\end{equation}
where $U(x,t): \mathbb{R}\times \mathbb{R}^+ \to \mathbb{R}$, $\mu, \nu \in \mathbb{R}$ are the bifurcation parameters, $2 < m \in  \mathbb{N}$, and $F: \mathbb{R}^2 \to \mathbb{R}$, $\Psi:\mathbb{R}^{m+2} \to \mathbb{R}$ are sufficiently smooth functions satisfying
\begin{equation*}
    \Psi(0,\dots,0, u;\nu) \equiv 0
\end{equation*}
for all $(u,\nu) \in \mathbb{R}^2$ (and thus clearly $\Psi_{u^i\nu^j}(0,\dots,0,u;\nu) \equiv 0$ for $(u,\nu, i, j) \in \mathbb{R}^2 \times \mathbb{N}^2_0$). Consequently, $\Psi$ does not contain nonlinear terms that do not contain spatial derivatives; these are all incorporated in $F$. Thus, the saddle-node is fully controlled by $\mu$ (and independent of $\nu$). 
On the other hand, \eqref{eq:illustration pde:rtab} is indeed more general than \eqref{eq:general scalar eq:gs}: unlike in \eqref{eq:general scalar eq:gs}, parameter $\nu$ may also drive nonlinear terms as $U U_{xx}$ in \eqref{eq:illustration pde:rtab}, and, more importantly, parameter $\mu$ will appear in the Turing analysis due to the presence of terms as $U U_{xx}$ (see below).

Now we follow the approach of \Cref{subsec:general setting:se} and suppose that a saddle-node bifurcation occurs at $\mu = \mu^* \in \mathbb{R}$, creating $u^\pm(\mu)$ as $\mu$ increases through $\mu^*$. Moreover, we assume that $F_u(u^+;\mu) < 0$ for $\mu > \mu^*$, in other words, that $u^+$ is stable in the associated ODE system $\dot{u} = F(u;\mu)$. Next, we suppose there exists a $\nu^* \in \mathbb{R}$, such that a Turing-fold bifurcation occurs at $(\mu^*, \nu^*) \in \mathbb{R}^2$, with critical wave number $k^*$ such that for $\nu < \nu^*$, $u^+$ destabilizes via a Turing bifurcation at $\mu^t > \mu^*$ with critical wave number $k^c$. Since $u^+ = u^+(\mu)$, it follows that both $\mu^t = \mu^t(\mu, \nu)$ and $k^c = k^c(\mu, \nu)$. Consequently, unlike in the case of model \eqref{eq:general scalar eq:gs}, the parameter $\mu$ enters into the $k$-dependent part of the dispersion relation (compare \eqref{eq:eigenvalue_problem_scalar_system} to \eqref{eq:eigenvalue_problem_scalar_system-2} below). 

For clarity of presentation, we omit the precise non-degeneracy conditions for the Turing-fold bifurcation in \eqref{eq:illustration pde:rtab} and only state that, under the appropriate non-degeneracy conditions, one can derive for $\nu = \nu^* - \delta$, with $0 < \delta \ll 1$,
\begin{equation*}
    u^t = u^* + \delta \tilde{u} + \mathcal{O}(\delta^2), \quad \mu^t = \mu^* + \delta^2 \hat{\mu} + \mathcal{O}(\delta^3), \quad \text{ and } \quad k^c = k^* + \delta \tilde{k} + \mathcal{O}(\delta^2)
\end{equation*}
(where again $u^t = u^+(\mu^t)$). The dispersion relation associated with $u^+$ is
\begin{equation}
\label{eq:eigenvalue_problem_scalar_system-2}
\omega(k;\mu,\nu) = F_u(u^+(\mu);\mu) + \sum_{j = 1}^m \left(k^2\right)^j\Psi_j(0,\dots, 0,u^+;\nu) \eqqcolon F_u(u^+(\mu);\mu) + 
G(u^+(\mu);k,\nu).
\end{equation}
By definition $\omega(k^c; \mu^t,\nu^* - \delta) = \omega_k(k^c; \mu^t,\nu^* - \delta) = F(u^t;\mu^t)\equiv 0$ for all $0 \leq \delta \ll 1$. Therefore, $G_k^* = 0$, and expanding gives
\begin{equation}
\label{eq:defH1}    
0 = \omega(k^c;\mu^t,\nu^* - \delta) = \delta \left((F_{uu}^* + G_u^*)\tilde{u}- G_\nu^*\right) + \mathcal{O}\left(\delta^2\right) \eqqcolon \delta H_1(\tilde{u}) + \mathcal{O}\left(\delta^2\right)
\end{equation}
and
\begin{equation}
\label{eq:defH2}    
0 = F(u^t;\mu^t) = \delta^2\left(\frac{1}{2}F_{uu}^* \tilde{u}^2 + F_\mu^* \hat{\mu}\right) + \mathcal{O}\left(\delta^3\right) \eqqcolon \delta^2 H_2(\tilde{u}, \hat{\mu}) + \mathcal{O}\left(\delta^3\right).
\end{equation}
Hence, $\tilde{u}$ and $\hat{\mu}$ satisfy $H_1(\tilde{u}) = 0$ and $H_2(\tilde{u}, \hat{\mu}) = 0$, respectively:
\begin{equation*}
    \tilde{u} = \frac{G_\nu^*}{F_{uu}^* + G_u^*} \quad \text{ and } \quad \hat{\mu} = - \frac{F_{uu}^* \tilde{u}^2}{2 F_\mu^*}.
\end{equation*}
The derivation of the AB-system will show that the $H_1$, and $H_2$ equations are inherently present in the AB-system. It is this specific fact, that allows us to normalize three coefficients by two transformations.

Following the procedure already outlined in \Cref{subsec:a model,subsec:general setting:se}, we (again) set $\mu = \mu^t - r \delta^2$ and introduce the (now standard) Ansatz
\begin{equation*}
    U_{AB}(x,t) = u^* + \delta B(\xi, \tau) + \delta \left(e^{i k^* x} A(\xi,\tau) + \text{c.c.}\right). 
\end{equation*}
Consequently, we expand
\begin{equation*}
    F(U_{AB};\mu^t- r \delta^2, \nu^* - \delta) + \Psi(\partial_x^2U_{AB}, \dots, \partial_x^{2m} U_{AB}, U_{AB}; \nu^* - \delta),
\end{equation*}
and collecting the terms at the $\delta^2 E^0$- and $\delta^2 E$-levels, where $E = e^{i k^* x}$. From the $\delta^2 E^0$-level we obtain:
\Cref{subsec:general setting:se},
\begin{equation}
\label{eq:genBeq}
    B_\tau = - \frac{1}{2}\rho_{kk}^*B_{\xi \xi} + H_2(B, \hat{\mu} - r) + c_8 |A|^2,
\end{equation}
and the $\delta^2 E$-level yields:
\begin{equation}
\label{eq:genAeq}
A_\tau = - \frac{1}{2} \omega_{kk}^* A_{\xi\xi} + H_1(B)A,
\end{equation}
where $H_1$ and $H_2$ are defined in \eqref{eq:defH1}, and \eqref{eq:defH2}. 

Next, we observe that 
\begin{equation*}
    \frac{H_1(\tilde{u}B)}{-G_\nu^*}A = A - AB \quad \text{ and } \quad \frac{H_2(\tilde{u} B, \hat{\mu}(1 - R))}{F_\mu^* \hat{\mu}} =  1 - R - B^2
\end{equation*}
which implies that the two transformations
\begin{equation*}
     B = \tilde{u} \tilde{B}, \quad \text{ and }  \quad r = \hat{\mu}R,
\end{equation*}
normalize three coefficients. Therefore, by applying the transformations
\begin{equation*}
    \tilde{\tau} = - G_\nu^* \tau, \quad \tilde{\xi} = \sqrt{\frac{2 G_\nu^*}{\omega_{kk}^*}}\xi, \quad B = \tilde{u} \tilde{B}, \quad \text{ and }  \quad r = \hat{\mu}R,
\end{equation*}
to the AB-system spanned by \eqref{eq:genAeq} and \eqref{eq:genBeq}, we recover the canonical form \eqref{eq:canonical form:csd:gs} for \eqref{eq:illustration pde:rtab}, with
\begin{equation}
\label{eq:scalingsalphaetc-gen}    \alpha = -\frac{F_\mu^* \hat{\mu}}{G_\nu^*\tilde{u}} = \frac{F_{uu}^*}{2(F_{uu}^* + G_{u}^*)}, \quad  d = -\frac{\rho_{kk}^* G_\nu^* \tilde{u}}{\omega_{kk}^*F_\mu^*\hat{\mu}} = \frac{\rho_{kk}^*}{\omega_{kk}^* \alpha}, \quad \text{ and } \quad \beta = \frac{c_8}{F_\mu^* \hat{\mu}}.
\end{equation}

Observe that, when $G(u^+(\mu);k,\nu)$ in dispersion relation \eqref{eq:eigenvalue_problem_scalar_system-2} is independent of $u$ (and thus $\mu$, via $u=u^+(\mu)$), we have $G^*_u = 0$, yielding $\alpha = \frac12$ in \eqref{eq:scalingsalphaetc-gen} -- which agrees with the scalings of \eqref{eq:rescaleABalphaetc} for model \eqref{eq:general scalar eq:gs} considered in \Cref{subsec:general setting:se}. Thus, indeed, $\alpha = \frac12$ for \eqref{eq:general scalar eq:gs} is not a generic property of the AB-system, but results from our special (simplified) choice of model \eqref{eq:general scalar eq:gs}. It should be noted though that, as in \eqref{eq:rescaleABalphaetc}, $\alpha d = \rho_{kk}^*/\omega_{kk}^*$ in this more general setting -- see also the brief observation on \eqref{rescalealphadbeta-Ncomp} in \Cref{subsec:AB-system derivation:rd}.

Finally, in AB-systems associated to a Turing-fold bifurcation in reaction-diffusion equations, as will be derived in \Cref{sec:reaction diffusion system}, the five coefficients can be normalized by the four rescalings (cf. \eqref{eq:rescalingsAB-syst}). Again, this is due to the nature of the derivation process, which yields the $AB$-system, and can be understood by uncovering an intrinsic structure of the asymptotic process, analogous to that presented here.

\section{Systems of reaction-diffusion equations}\label{sec:reaction diffusion system}
In this section, we consider the standard $n$-component reaction-diffusion system
\begin{equation}\label{eq:reaction-diffusion system:rds}
    \partial_t U = F(U;\mu,\nu) + D \partial_x^2 U,
\end{equation}
which we assume to be well-posed. Here, $n \in \mathbb{N}$, $U(x,t): \mathbb{R}\times \mathbb{R}^+ \to \mathbb{R}^n$, $\mu, \nu$ are real-valued parameters, $D$ is a diagonal diffusion matrix, and $F:\mathbb{R}^{n+2} \to \mathbb{R}^n$ is sufficiently smooth. Additionally, we adopt superscript notation for the entries of vectors and matrices -- for instance, $U = (U^1, \dots, U^n)^T$ and $F = (F^1, \dots, F^n)^T$ -- so that subscript notation may be used to denote partial derivatives with respect to the $j$-th coordinate, i.e., $F_j \coloneqq (F^1_j, \dots F^n_j)$. As in \Cref{sec:scalar}, we will omit parameter dependence when it is clear from context.

The spatially homogeneous ODE corresponding to \eqref{eq:reaction-diffusion system:rds} is given by
\begin{equation*}
    \dot{u} = F(u; \mu, \nu).
\end{equation*}
We assume that there exists a non-empty open interval $I \subset \mathbb{R}$ such that for each $\nu \in I$, a saddle-node bifurcation occurs at $\mu = \mu^s(\nu)$, resulting in the emergence of two solutions, $u^\pm(\mu,\nu)$, with $u^+ > u^-$, as $\mu$ increases past $\mu^s(\nu)$. Additionally, we assume that $u^+$ is stable in the ODE context; that is, for $\mu > \mu^s$, all the eigenvalues of the Jacobian matrix
\begin{equation*}
    F_u(u^+; \mu, \nu) :=
    \begin{pmatrix}
    F^1_{1}(u^+; \mu, \nu) & \dots & F^1_{n}(u^+; \mu, \nu) \\
    \vdots & \ddots & \vdots \\
    F^n_{1}(u^+; \mu, \nu) & \dots & F^n_{n}(u^+; \mu, \nu)
    \end{pmatrix}
\end{equation*}
lie in the left half of the complex plane.

The conditions for a saddle-node bifurcation ensure the presence of a unique zero eigenvalue \cite{homburg2024bifurcation}, implying that the Jacobian $F_u(u^+; \mu^s, \nu)$ is non-invertible. Let $v^s(\nu)$ denote the associated right zero eigenvector (i.e., $F_u(u^+;\mu^s,\nu)v^s = 0 $), and let $p^s(\nu)$ denote the associated left zero eigenvector. We assume the normalization $\langle v^s, v^s \rangle = \langle p^s, v^s \rangle = 1$, where $\langle \cdot, \cdot \rangle$ denotes the standard dot product in $\mathbb{R}^n$.

The signed non-degeneracy conditions for this specific (standard) saddle-node bifurcation are
\begin{equation}\label{eq:non deg saddle node:rds}
    \langle F_{uu}(u^+; \mu^s, \nu)(v^s,v^s),p^s \rangle < 0 \quad \text{ and } \quad \langle F_\mu(u^+;\mu^s,\nu) v^s, p^s \rangle > 0,
\end{equation}
see \cite{homburg2024bifurcation}, where
\begin{equation*}
    F_{uu}(u;\mu,\nu)(v,w) \coloneqq
    \begin{pmatrix}
        \left \langle \nabla^2_u F^1(u; \mu, \nu)v, w \right \rangle \\
        \vdots  \\
        \left \langle \nabla^2_u F^n(u; \mu, \nu)v, w \right \rangle \\
    \end{pmatrix},
\end{equation*}
with $\nabla^2_u F^j$ denoting the $n \times n$ Hessian of the $j$-th component of $F$.

\begin{figure}[t]
    \centering
    \includegraphics[width=0.6\textwidth]{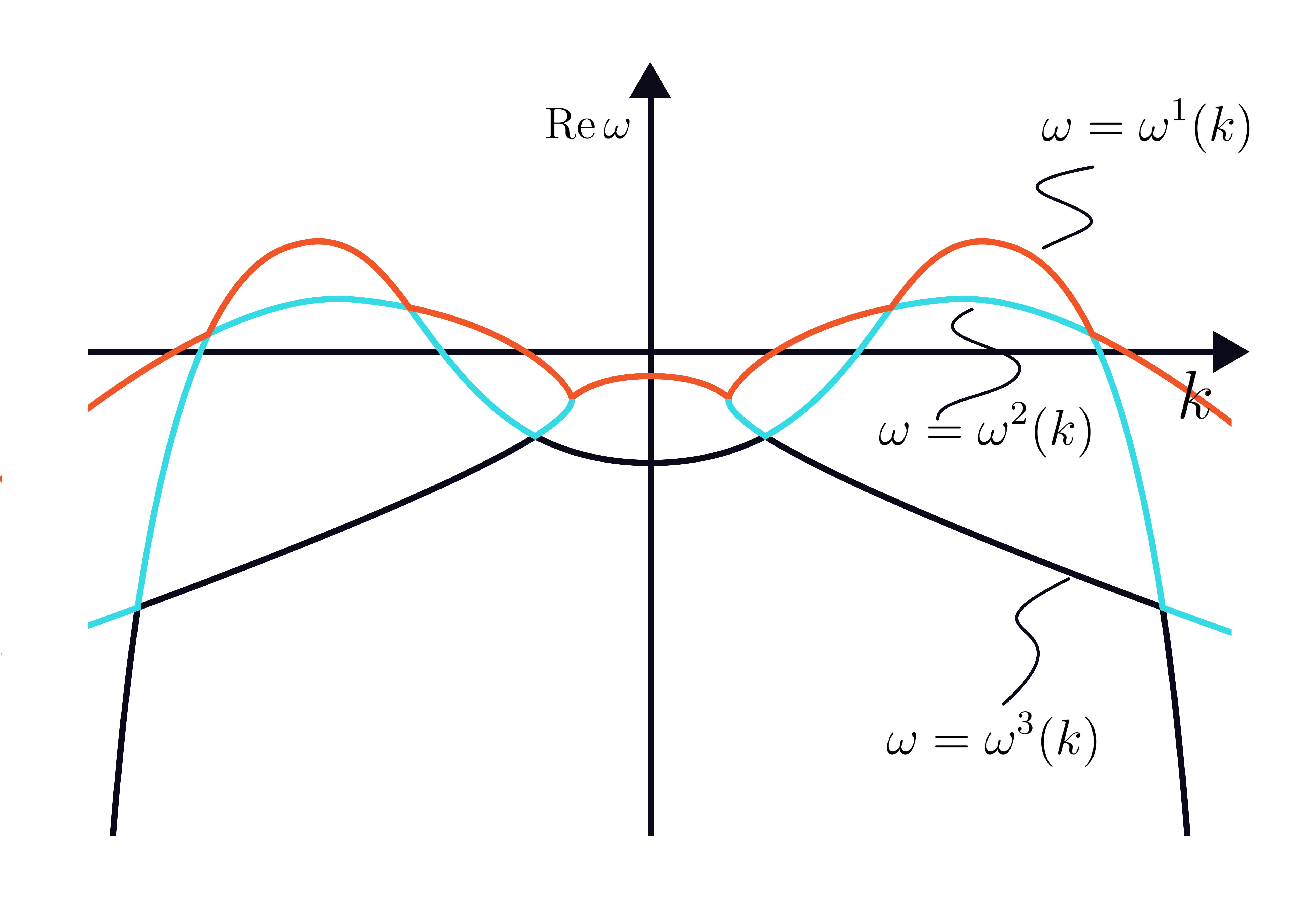}
    \caption{Typical eigenvalue curves of a 3-component reaction-diffusion system. Since the $\omega^j$'s are defined such that $\text{Re} \, \omega^{j-1} \leq \text{Re} \, \omega^j \leq \text{Re} \, \omega^{j+1}$ the eigenvalue curves are only piecewise differentiable.}
    \label{fig:Omega Curves}
\end{figure}

To determine the spectral stability of $u^+$ as a solution of \eqref{eq:reaction-diffusion system:rds}, we consider the perturbation
\begin{equation}\label{eq:perturbed solution:rds}
    U = u^+ + \left(\varepsilon v e^{ikx + \omega(k; \mu, \nu) t} + \text{c.c.}\right),
\end{equation}
where $v \in \mathbb{R}^n$. Substituting this perturbation into \eqref{eq:reaction-diffusion system:rds} and linearize around $u^+$ yields
\begin{equation}\label{eq:T def:rd}
        \omega(k;\mu, \nu) v = \left(F_u(u^+;\mu,\nu) - k^2 D \right) v \coloneqq T(k; \mu, \nu) v.
\end{equation}
Let $\omega^j(k; \mu, \nu)$ with $j = 1, \dots, n$, denote the ordered eigenvalues of the matrix $T(k, \mu, \nu)$, such that $\text{Re}\, \omega^j(k) \geq \text{Re}\, \omega^{j+1}(k)$ for all $k \in \mathbb{R}$ (cf. \Cref{fig:Omega Curves}). Henceforth, we shall denote $\omega^1$ by $\omega$. 

Just as in \Cref{subsec:general setting:se}, we assume that there exist $\nu^*, \bar{\nu} \in I$ with $\bar{\nu} < \nu^*$ such that for all $\nu \in (\bar{\nu},\nu^*)$, the steady state $u^+(\mu,\nu)$ undergoes a non-degenerate Turing bifurcation with critical wave number $k^c(\nu) > 0$ as $\mu$ decreases through a critical value $\mu^t(\nu)$ with $\mu^t(\nu) > \mu^s(\nu)$. The pair $(k^c(\nu),\mu^t(\nu))$ is determined by
\begin{equation}
\label{eq:genreactiondiffusionTuring:rds}    
\omega(k^c;\mu^t,\nu) = 0, \quad \omega_k(k^c;\mu^t,\nu) = 0, \quad \text{ and } \quad \text{Re}\, \omega(k; \mu^t, \nu) < 0 \quad \forall k \in \mathbb{R}\backslash\{-k^c,k^c\}.
\end{equation}
The signed non-degeneracy conditions for this specific Turing bifurcation are given by
\begin{equation*}
    \omega_{kk}(k^c; \mu^t, \nu) < 0 \quad \text{ and } \quad \omega_{\mu}(k^c;\mu^t,\nu) < 0.
\end{equation*}
Additionally, we assume that a Turing-fold bifurcation occurs at $(\mu^*,\nu^*) \in \mathbb{R}^2$, which implies that 
\[\mu^t(\nu^*) = \mu^s(\nu^*) = \mu^* \]
and that the (non-degeneracy) conditions outlined in \Cref{subsec:non-degeneracy conditions:rd} are satisfied. We thus (implicitly) assume that $\omega(k;\mu,\nu) \in \mathbb{R}$ near $k^c,\mu^t$; that is, that the Turing bifurcation yields stationary patterns. In other words, we assume that the destabilization at $k = \pm k^c$ is of (`pure') Turing type and not of Turing-Hopf or oscillatory/traveling Turing type (which would imply $\text{Im} \, \omega(k^c;\mu^t,\nu) \neq 0$) \cite{rademacher2007instabilities, scheel2003radially}. 

To aid the upcoming analysis, we introduce the following notation. First, the right and left zero eigenvectors of $T(k^c; \mu^t, \nu)$ are denoted by $v^t(\nu)$ and $p^t(\nu)$ respectively. Second, we define the characteristic polynomials
\begin{equation}\label{eq:P and Q definition:rds}
    P(\lambda; u, \mu, \nu, k) \coloneqq \det\left(F_u(u;\mu,\nu) - k^2 D - \lambda I \right) \quad \text{ and } \quad Q(\lambda; u,\mu, \nu) \coloneqq P(\lambda;u,\mu,\nu,0).
\end{equation}
Thirdly, we use superscript $t$, $*$, and $s$, to denote evaluation at
\begin{equation*}
    (u^t,\mu^t, \nu^* - \delta, k^c), \quad (u^*,\mu^*, \nu^*, k^*), \quad \text{ and } \quad (u^s,\mu^s, \nu^* - \delta, k^c),
\end{equation*}
respectively, where,
\begin{equation*}
    u^t \coloneqq u^+(\mu^t), \quad u^* \coloneqq u^+(\mu^*), \quad u^s \coloneqq u^+(\mu^s),
\end{equation*}
and $0 < \delta \ll 1$ is a small parameter. Fourthly, when $P$, $Q$, or one of their derivatives is marked by a superscript $t$ or $*$, we assume $\lambda = 0$; additionally, if it is marked by a superscript $s$, we assume $\lambda = \omega^s \coloneqq \omega(k^c;\mu^s,\nu^*-\delta)$. Finally, we denote
\begin{equation*}
    v^s_* \coloneqq v^s(\nu^*), \quad p^s_* \coloneqq p^s(\nu^*), \quad v^t_* \coloneqq v^t(\nu^*), \quad \text{ and } \quad p^t_* \coloneqq p^t(\nu^*).
\end{equation*}

\subsection{Non-degeneracy conditions}\label{subsec:non-degeneracy conditions:rd}

\begin{figure}[t]
    \centering
    \includegraphics[width=1\linewidth]{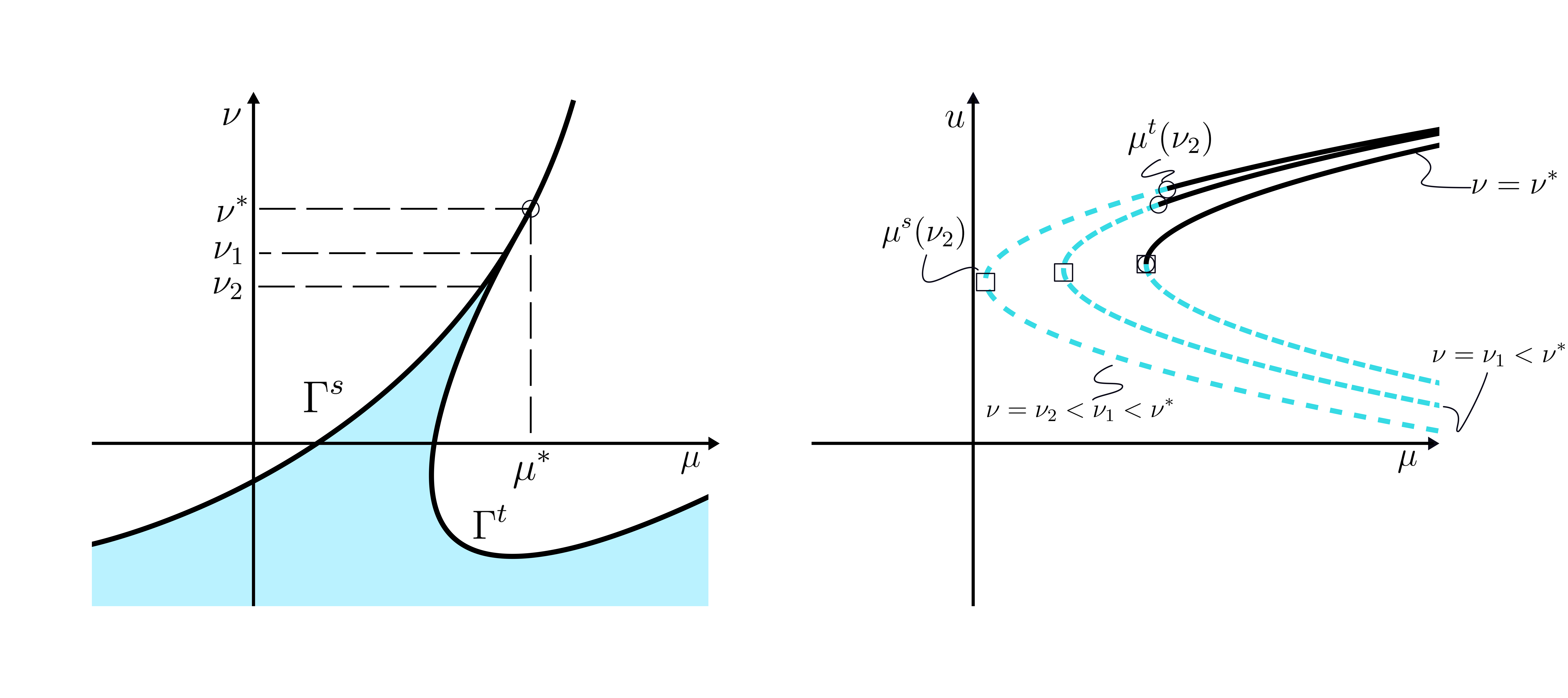}
    \caption{Typical bifurcation diagrams associated to a Turing-fold bifurcation in a reaction-diffusion system. In the left image, the curves $\Gamma^s$ (saddle-node points), and $\Gamma^t$ (Turing points) are plotted in the $(\mu, \nu)$-plane. These curves meet tangentially at the co-dimension 2 point $(\mu^*, \nu^*)$. Unlike in \Cref{sec:scalar}, $\mu^s$ may now depend on $\nu$, so $\Gamma^s$ is not necessarily a vertical line (cf. \Cref{fig:pde bifurcation diagram scalar example:a model}). The blue area beneath these curves represents the region where the Turing bifurcation has occurred, but the saddle-node bifurcation has not. In the right image, three typical bifurcation diagrams are shown, each corresponding to a different value of $\nu$.}
    \label{fig:eigenvalue picture:ndc:rd}
\end{figure}

In the reaction-diffusion setting, a co-dimension 2 Turing-fold bifurcation occurs under the following conditions:
\begin{equation}
\label{eq:reactiondiffusion-TuringFold}    
\begin{array}{rrrr}
{\rm (Fold)}_{\rm bif}:  &  F^* = 0, & \rho^* \coloneqq \omega(0;\mu^*,\nu^*) = 0; & \\
{\rm (Turing)}_{\rm bif}:  &  \omega^* = 0, & \omega_{k}^* = 0, & \text{Re} \, \omega(k; \mu^*, \nu^*) < 0 \; \forall \; k \in \mathbb{R}\backslash\{-k^*,0,k^*\}.
\end{array}
\end{equation}
Additionally, we consider such a bifurcation non-degenerate when the following (signed) conditions are also satisfied:
\begin{equation}
\label{eq:reactiondiffusionnondeg-TuringFold}    
\begin{array}{rrrr}
{\rm (Fold)}_{\rm n-d}: & \left\langle F_\mu^*, p^s_* \right\rangle > 0, &  \left \langle F_{uu}^*(v^s_*, v^s_*), p^s_* \right\rangle < 0, & \det \begin{pmatrix}
    F^*_u & F^*_\mu \\
    (Q^*_u)^T & Q_\mu^*
\end{pmatrix} \neq 0,
\\
 & \rho_{kk}^* \coloneqq \omega_{kk}(0;\mu^*,\nu^*) < 0, &\text{Re} \, \omega^2(0; \mu^*, \nu^*) < 0; \\
{\rm (Turing)}_{\rm n-d}:  &  \omega_{kk}^* < 0, & \left\langle P_u^*,v^s_*\right\rangle P_\lambda^* > 0, & \left \langle F_{uu}^*(v^s_*,v^t_*),p^t_* \right \rangle < 0, \\
& \left(P_\nu^* - P_\mu^* \tilde{\mu}^s - \left \langle P_u^*,\tilde{\mu}^s \right \rangle \right)P_\lambda^* > 0, & \text{Re} \, \omega^2(k^*; \mu^*, \nu^*) < 0; \\
\end{array}
\end{equation}
(see \eqref{eq:order of variables:ndc:rd} for the definition of $\tilde{\mu}^s$).  As in \Cref{sec:scalar}, the conditions \eqref{eq:reactiondiffusion-TuringFold} and \eqref{eq:reactiondiffusionnondeg-TuringFold} together imply the existence of a curve of saddle-node points $\Gamma^s \subset \mathbb{R}^2$ and a curve of Turing points $\Gamma^t \subset \mathbb{R}^2$, which meet tangentially at the co-dimension 2 point $(\mu^*, \nu^*) \in \Gamma^s \cup \Gamma^t$ -- see \Cref{fig:eigenvalue picture:ndc:rd}. 

Note that the condition on $\rho_{kk}^*$ reflects the spatial extension of the ODE-based non-degeneracy conditions in \eqref{eq:non deg saddle node:rds} -- as in \Cref{subsubsec:non-degeneracy conditions:gs}. Moreover, the conditions on $\omega^2$ ensure that the second eigenvalue branch of $T(k;\mu, \nu)$ does not interact with the bifurcations induced by the critical $\omega = \omega^1$ branch.

The specific signs are (again) chosen to mimic the setup of \Cref{sec:scalar}; that is, $u^+$ and $u^-$ are created as $\mu$ increases through $\mu^s$, and for $\nu \in (\bar{\nu}, \nu^*)$ with $\bar{\nu} < \nu^*$, the Turing bifurcation takes place as $\mu$ decreases through $\mu^t$. In fact, for $\nu = \nu^* - \delta$  with $0 < \delta \ll 1$ we have:
\begin{equation}\label{eq:order of variables:ndc:rd} 
\begin{array}{llll}
    & u^s = u^* + \tilde{u}^s \delta + \hat{u}^s \delta^2 + \mathcal{O}(\delta^3),  & \mu^s = \mu^* + \tilde{\mu}^s \delta + \hat{\mu}^s \delta^2 + \mathcal{O}(\delta^3), & v^s = v^s_* + \mathcal{O}(\delta), \\ & p^s = p^s_* + \mathcal{O}(\delta) , &  \omega^s = \tilde{\omega}^s \delta + \mathcal{O}(\delta^2), & u^t = u^s + \tilde{u}^t \delta + \mathcal{O}(\delta^2), \\  & \mu^t = \mu^s + \hat{\mu}^t \delta^2 + \mathcal{O}(\delta^3), & k^c = k^* + \tilde{k} \delta + \mathcal{O}(\delta^2), & \omega_\mu^t = \frac{\tilde{\omega}_\mu^t}{\delta} + \mathcal{O}(1),
\end{array}
\end{equation}
with
\begin{equation}\label{eq:values of variables:ndc:rd}
\begin{array}{lll}
    & \tilde{\omega}^s =  \frac{P_\nu^* - P_\mu^* \tilde{\mu}^s - \left \langle P_u^*, \tilde{u}^s \right \rangle}{P_\lambda^*}, 
    & \begin{pmatrix}
    \tilde{u}^s \\ \tilde{\mu}^s
\end{pmatrix} = \begin{pmatrix}
    F^*_u & F^*_\mu \\
    (Q^*_u)^T & Q_\mu^*
\end{pmatrix}^{-1}
\begin{pmatrix}
    F_\nu^* \\
    Q_\nu^*
\end{pmatrix},\\
    & \hat{\mu}^t = -\frac{1}{2} \frac{\left\langle F_{uu}^*(v^s_*, v^s_*), p^s_* \right\rangle\left(P_\lambda^* \tilde{\omega}^s\right)^2}{\left\langle F_{\mu}^*,p^s_*\right\rangle\left\langle P^*_u, v_*^s \right\rangle^2}, & \tilde{u}^t \coloneqq \bar{u}^t v^s_* =  \sqrt{-2\frac{\left\langle F_{\mu}^*,p^s_*\right\rangle}{\left\langle F_{uu}^*(v^s_*, v^s_*), p^s_* \right\rangle}\hat{\mu}^t}v^s_* = \frac{P_\lambda^* \tilde{\omega^s}}{\left \langle P^*_u,v^s_* \right \rangle}v^s_*, \\ & \tilde{\omega}^t_\mu = \frac{\bar{u}^t}{\hat{\mu}^t} \left \langle F_{uu}^*\left(v^s_*,v^t_*\right),p^t_* \right \rangle, & \tilde{k} = \frac{P_{k \nu}^* - P_{k \mu}^* \tilde{\mu}^s - \left \langle P_{ku}^*, \tilde{u}^t + \tilde{u}^s\right \rangle}{P_{kk}^*}.\\
\end{array}
\end{equation}

These results can, in essence, be deduced from direct expansion using basic properties of the saddle-node and Turing bifurcation. We observe:

\[
\begin{array}{rcl}
0 & = & F(u^s;\mu^s,\nu^* - \delta) \\
& = & F^* + F_u^*(u^s - u^*) + F_\mu^*(\mu^s - \mu^*) - F_\nu^*\delta + \text{h.o.t.} \\
& = & F_u^*(u^s - u^*) + F_\mu^*(\mu^s - \mu^*) - F_\nu^*\delta + \text{h.o.t.,}
\end{array}
\]
since $F^* = 0$ by \eqref{eq:reactiondiffusion-TuringFold}, and
\[
\begin{array}{rcl}
0 & = & Q(0;u^s,\mu^s,\nu^* - \delta) \\
& = & \left \langle Q_u^*, u^s-u^* \right \rangle+Q_\mu^*(\mu^s-\mu^*)- Q_\nu^*\delta + \text{h.o.t.}
\end{array}
\]
Consequently, for sufficiently small $|\delta|$, \eqref{eq:reactiondiffusionnondeg-TuringFold} ensures the existence of $u^s$ and $\mu^s$ such that:
\begin{equation*}
    u^s = u^* + \tilde{u}^s \delta + \mathcal{O}(\delta^2) \quad \text{ and } \quad \mu^s = \mu^* + \tilde{\mu}^s \delta + \mathcal{O}(\delta^2),
\end{equation*}
where $\tilde{u}^s$ and $\tilde{\mu}^s$ are determined by:
\begin{equation*}
    \begin{pmatrix}
    \tilde{u}^s \\ \tilde{\mu}^s
\end{pmatrix} = \begin{pmatrix}
    F^*_u & F^*_\mu \\
    (Q^*_u)^T & Q_\mu^*
\end{pmatrix}^{-1}
\begin{pmatrix}
    F_\nu^* \\
    Q_\nu^*
\end{pmatrix}.
\end{equation*}

For $0 < \mu - \mu^s \ll 1$ we have $0 < |u^\pm(\mu,\nu) - u^s| \ll 1$. Hence:
\begin{equation}\label{eq:ref once 1:ndc:rd}
\begin{array}{rcl}
    0 & = &F(u^\pm;\mu,\nu^*-\delta) \\
    & = & F^s + F_u^s (u^\pm - u^s) + \frac{1}{2} F_{uu}^s(u^\pm - u^s, u^\pm - u^s)  + F_\mu^s (\mu - \mu^s) + \text{h.o.t.} \\
    & = & F_u^s (u^\pm - u^s) + \frac{1}{2} F_{uu}^s(u^\pm - u^s, u^\pm - u^s)  + F_\mu^s (\mu - \mu^s) + \text{h.o.t.}
\end{array}
\end{equation}
since $F^s = 0$. Consequently, because $F_u^s v^s = 0$, we have: 
\begin{equation}\label{eq:ref once 2:ndc:rd}
    u^\pm - u^s = \pm c  v^s \sqrt{\mu - \mu^s} + \mathcal{O}(\mu - \mu^s),
\end{equation}
with $c > 0$. Substituting \eqref{eq:ref once 2:ndc:rd} into \eqref{eq:ref once 1:ndc:rd} and then taking the dot-product with $p^s$, we get:
\begin{equation}
    c^2\frac{1}{2}\left\langle F_{uu}^s(v^s, v^s), p^s \right\rangle = -\left\langle F_{\mu}^s,p^s\right\rangle.
\end{equation}
Therefore:
\begin{equation}\label{eq:u pm equation:ndc:rd}
    u^\pm = u^s \pm v^s\sqrt{-2\frac{\left\langle F_{\mu}^s,p^s\right\rangle}{\left\langle F_{uu}^s(v^s, v^s), p^s \right\rangle}(\mu-\mu^s)} + \mathcal{O}(\mu-\mu^s).
\end{equation}

Moreover, let $s^r(\mu, \nu)$ and $s^l(\mu, \nu)$ denote the right and left eigenvalues, respectively, associated with $\rho(\mu, \nu) \coloneqq \omega(0,\mu,\nu)$. Clearly, $s^r = v^s + \tilde{s}^r \sqrt{\mu - \mu^s} + \mathcal{O}(\mu - \mu^s)$, and $s^l = p^s + + \tilde{s}^l \sqrt{\mu - \mu^s} + \mathcal{O}(\mu - \mu^s)$. Therefore:
\begin{equation}\label{eq:rho things:ndc:rd}
\begin{array}{rcl}
    \rho(\mu, \nu^*-\delta) & \coloneqq & \left \langle F_u(u^+;\mu,\nu^* - \delta) s^r, s^l \right \rangle\\
    & = & \left \langle (F_u^s + dF_u^s(u^+ - u^s))s^r, s^l  \right \rangle + \mathcal{O}(\mu - \mu^s) \\
    & = &  \left \langle F_{uu}^s(v^s,v^s),p^s \right \rangle \sqrt{-2\frac{\left\langle F_{\mu}^s,p^s\right\rangle}{\left\langle F_{uu}^s(v^s, v^s), p^s \right\rangle}(\mu-\mu^s)} + \left\langle F_u^s \tilde{s}^r, p^s \right \rangle \sqrt{\mu - \mu^s}+ \mathcal{O}(\mu - \mu^s) \\
    & = & - \sqrt{-2 \left \langle F_{uu}^s(v^s, v^s), p^s \right \rangle \left \langle F_{\mu}^s, p^s \right \rangle(\mu - \mu^s)} + \mathcal{O}(\mu - \mu^s).
\end{array}
\end{equation}
Note that $\left\langle F_u^s \tilde{s}^r, p^s \right \rangle = 0$, since matrix multiplication is associative and $(p^s)^TF_u^s = 0$. Moreover, in the last line we have to add a minus sign in front of the square root, because:
\begin{equation*}
    \left\langle F_{uu}^s(v^s,v^s),p^s \right\rangle = \left\langle F_{uu}^*(v^s_*,v^s_*),p^s_* \right\rangle + \mathcal{O}(\delta) < 0 \quad \text{ and } \quad \left\langle F_{\mu}^s,p^s \right\rangle = \left\langle F_{\mu}^*,p^s_* \right\rangle + \mathcal{O}(\delta) > 0.
\end{equation*}
Here, we have used that $v^s = v^s_* + \mathcal{O}(\delta)$, and $p^s = p^s_* + \mathcal{O}(\delta)$. 

In conclusion, for sufficiently small $|\delta|$, the conditions proposed in \eqref{eq:reactiondiffusion-TuringFold} and \eqref{eq:reactiondiffusionnondeg-TuringFold} indeed guarantee the existence of unique saddle-node bifurcation point $(u^s,\mu^s)$ with $(u^s(\nu^*), \mu^s(\nu^*)) = (u^*,\mu^*)$. Moreover, $u^+$ and $u^-$ are created as $\mu$ increases through $\mu^s$, and $u^+$ is stable for $0 < \mu - \mu^s \ll 1$.

By definition, \eqref{eq:P and Q definition:rds},
\begin{equation*}
    P(0;u^+,\mu,\nu,k) = \prod_{j=1}^n \omega^j(k;\mu,\nu),
\end{equation*}
hence,
\begin{equation*}
    P_{kk}^* = \omega_{kk}^* \prod_{j=2}^n \omega^j(k^*;\mu^*,\nu^*) \neq 0,
\end{equation*}
since $\omega^* = \omega_k^* = 0$, $\omega_{kk} < 0$, and $\text{Re} \, \omega^j(k^*;\mu^*,\nu^*) < 0$ for $j \neq 1$ -- see \eqref{eq:reactiondiffusion-TuringFold}, and \eqref{eq:reactiondiffusionnondeg-TuringFold}. Consequently, we derive the following three equations via basic expansion,
\begin{align}
&\begin{array}{rcl} \label{eq:eq 1:ndc:rd}
    0 & = & P_k^t \coloneqq P_k(0;u^t,\mu^t,\nu^*-\delta,k^c) \\
    & = & P_{kk}^*(k^c-k^*) + P_{k\mu}^*(\mu^t - \mu^*) + \left\langle P_{ku}^*, u^t - u^* \right\rangle - P_{k\nu}^*\delta + \text{h.o.t.} \\
\end{array}  \\ \notag  \\
&\begin{array}{rcl} \label{eq:eq 2:ndc:rd}
    0 & = & P^s \coloneqq P(\omega^s;u^s,\mu^s,\nu^*-\delta, k^c) \\
    & = & \frac{1}{2} P_{kk}^*(k^c-k^*)^2 + P_{\mu}^*(\mu^s - \mu^*) + \left\langle P_{u}^*, u^s - u^* \right\rangle - P_{\nu}^*\delta + P_\lambda^*\omega^s + \text{h.o.t.}
\end{array}  \\ \notag \\
& \begin{array}{rcl}  \label{eq:eq 3:ndc:rd}
    0 & = & P^t \\
    & = & \left \langle P_u^s, u^t-u^s \right \rangle + P_\mu^s(\mu^t-\mu^s) - P_\lambda^s \omega^s + \text{h.o.t.} \\
    & = & \left \langle P_u^s, u^t-u^s \right \rangle - P_\lambda^s \omega^s + \text{h.o.t.},
\end{array}
\end{align}
where $P_\mu^s(\mu^t - \mu)$ in \eqref{eq:eq 3:ndc:rd} has been relegated to the h.o.t.-bucket, since $\mu^t - \mu^s = \mathcal{O}(\|u^t - u^s\|^2)$ -- see \eqref{eq:u pm equation:ndc:rd}. Rather providing a cumbersome, fully detailed argument -- one that would require repeatedly justifying that a quantity is of a certain order because it is composed of terms of various other orders -- we instead present a sketch of what follows from the three equations.

We have three equations, i.e., \eqref{eq:eq 1:ndc:rd}, \eqref{eq:eq 1:ndc:rd}, and \eqref{eq:eq 3:ndc:rd}. Furthermore, at leading order we are (essentially) dealing with only three unknown quantities: $\omega^s$, $k^c - k^*$, and $u^t - u^s$, because:
\begin{itemize}
    \item We have already determined $u^s - u^*$ and $\mu^s - \mu^*$.
    \item We can decompose $u^t - u^* = (u^t - u^s) + (u^s - u^*)$ and $\mu^t - \mu^* = (\mu^t - \mu^s) + (\mu^s - \mu^*)$.
    \item Although $u^t - u^s$ is a vector of size $n$, it has the form as presented in \eqref{eq:u pm equation:ndc:rd}; thus, effectively, determining $u^t-u^s$ to leading order reduces to finding the scalar quantity $\bar{u}^t$.
    \item Once $\bar{u}^t$ is known, $\hat{\mu}^t$ follows directly from \eqref{eq:u pm equation:ndc:rd}. 
\end{itemize}
Consequently, we obtain:
\begin{equation*}
    \tilde{\omega}^s = \frac{P_\nu^* - P_\mu^* \tilde{\mu}^s - \left \langle P_u^*, \tilde{u}^s \right \rangle}{P_\lambda^*}, \quad \tilde{k} = \frac{P_{k\nu}^* - P_{k\mu}^*\tilde{\mu}^s - \left\langle P_{ku}^*, \tilde{u}^t + \tilde{u}^s \right \rangle}{P_{kk}^*}.
\end{equation*}
\begin{equation}
    \hat{\mu}^t = -\frac{1}{2} \frac{\left\langle F_{uu}^*(v^s_*, v^s_*), p^s_* \right\rangle\left(P_\lambda^* \tilde{\omega}^s\right)^2}{\left\langle F_{\mu}^*,p^s_*\right\rangle\left\langle P^*_u, v_*^s \right\rangle^2}, \quad \text{ and } \quad
    \tilde{u}^t \coloneqq \bar{u}^tv^s_* = \sqrt{-2\frac{\left\langle F_{\mu}^*,p^s_*\right\rangle}{\left\langle F_{uu}^*(v^s_*, v^s_*), p^s_* \right\rangle}\hat{\mu}^t}v^s_* = \frac{P_\lambda^* \tilde{\omega^s}}{\left \langle P^*_u,v^s_* \right \rangle}v^s_*.
\end{equation}

Next, we show that $\omega^t_\mu$ is of the form $\omega_{\mu}^t = \tilde{\omega}_\mu^t \delta^{-1} + \mathcal{O}(1)$ with $\tilde{\omega}_\mu^t < 0$, using the following lemma \cite{Lancaster1964}.

\begin{lemma}\label{lem:lancaster:ndc:rd}
    Let $C:\mathbb{R} \to \mathbb{R}^{n\times n}$ be a matrix whose entries smoothly depends on a parameter $t\in\mathbb{R}$, let $\lambda(t)$ be an eigenvalue of $C(t)$, let $t_0 \in \mathbb{R}$ be such that $\lambda(t_0)$ is a simple eigenvalue of $C(t_0)$, and, finally, let $p(t)$ and $v(t)$ represent the corresponding left and right eigenvectors of $\lambda(t)$, respectively. Then, at $t_0$, we have
        \begin{equation}
        \lambda'(t_0) = \frac{\langle C'(t_0) v(t_0), p(t_0) \rangle}{\langle v(t_0), p(t_0) \rangle},
        \end{equation}
    where $C'(t)$ denotes the element-wise derivative of $C(t)$ with respect to $t$.
    \end{lemma}
By \eqref{eq:u pm equation:ndc:rd}, we have
\begin{equation}\label{eq:ut deriv:ndc:rd}
    u_\mu^t = -\frac{\left \langle F_\mu^*, p^s_* \right \rangle}{\bar{u}^t\delta\left \langle F_{uu}^*(v^s_*,v^s_*), p^s_* \right \rangle}v^s_* + \mathcal{O}(1)  = \frac{\bar{u}^t}{2 \hat{\mu}^t \delta}v^s_* + \mathcal{O}(1),
\end{equation}
so, by \Cref{lem:lancaster:ndc:rd}, we obtain
\begin{equation}\label{eq:once 5:ndc:rd}
    \omega^t_\mu = \left\langle (\frac{d}{d\mu} F_u^t) v^t, p^t\right\rangle = \left \langle (d F_{u}^* (u^t_\mu)) v^t_*, p^t_*\right \rangle + \mathcal{O}(1),
\end{equation}
where, given a function $A:\mathbb{R}^n \to \mathbb{R}^{n\times n}$, we let $(d A(u))v$ denote
\begin{equation}
    (dA(u)) v \coloneqq
    \begin{pmatrix}
    \nabla_u A^{11}(u) \cdot v & \cdots & \nabla_u A^{1n}(u) \cdot v\\
    \vdots & \ddots & \vdots \\
    \nabla_u A^{n1}(u) \cdot v & \cdots & \nabla_u A^{nn}(u)\cdot v \\
    \end{pmatrix}.
\end{equation}
Consequently, using Einstein summation notation, we see that for vectors $v, w \in \mathbb{R}^n$ we have
\begin{align}
    (dF_u(v))w =
    \begin{pmatrix}
        \nabla_u \frac{\partial F^1}{\partial u^1} \cdot v & \dots &\nabla_u \frac{\partial F^1}{\partial u^n} \cdot v \\
        \vdots & \ddots & \vdots \\
        \nabla_u \frac{\partial F^n}{\partial u^1} \cdot v & \dots &\nabla_u \frac{\partial F^n}{\partial u^n} \cdot v
    \end{pmatrix} w = \begin{pmatrix}
        w^j \frac{\partial^2 F^1}{\partial u^i \partial u^j} v^i \\
        \vdots \\
        w^j \frac{\partial^2 F^n}{\partial u^i \partial u^j} v^i \\
    \end{pmatrix}
    = F_{uu}(w,v) = F_{uu}(v,w).
\end{align} 
Thus, from \eqref{eq:once 5:ndc:rd}, it follows that
\begin{equation}
    \omega_{\mu}^t = \frac{\tilde{\omega}_\mu^t}{\delta} + \mathcal{O}(1) \quad \text{ with } \quad \tilde{\omega}_\mu^t = \frac{\bar{u}^t}{2\hat{\mu}^t} \left \langle F_{uu}^*\left(v^s_*,v^t_*\right),p^t_* \right \rangle.
\end{equation}
Moreover, \eqref{eq:reactiondiffusionnondeg-TuringFold} implies $\tilde{\omega}_{\mu}^t < 0$.

Finally, using the same arguments a presented in \Cref{subsubsec:non-degeneracy conditions:gs}, we may conclude that $k^* \neq 0$ and that the global conditions on $\omega(k;\mu,\nu)$ are satisfied for $\nu = \nu^* - \delta$ if $0 < \delta \ll 1$, since $\rho_{kk}^* < 0$ \eqref{eq:reactiondiffusionnondeg-TuringFold}.

\subsection{Extra notation, definitions and identities}\label{subsec:extra notations definitions and lemmas:rd}
To aid the upcoming analysis, we introduce additional notations, definitions, and approximations. First, both $F_u^t$ and $T^t$ \eqref{eq:T def:rd} are of the form
\begin{equation*}
    F_u^t = F_u^* + \tilde{F}_u^t\delta + \mathcal{O}(\delta^2) \quad \text{ and } \quad T^t = T^* + \tilde{T}^t\delta + \mathcal{O}(\delta^2),
\end{equation*}
with
\begin{equation}\label{eq:tilde T and tilde S:endl:rd}
    \tilde{F}_u^t = d F_u^*\left(\tilde{u}^t + \tilde{u}^s\right) + F^*_{u\mu} \tilde{\mu}^s - F^*_{u\nu}, \quad T^* = F_u^* - (k^*)^2D, \quad \text{ and } \quad \tilde{T}^t = \tilde{F}_u^t - 2 \tilde{k}k^* D.
\end{equation}

Second, using $\mathcal{O}(\|v^t - v^t_*\|) = \mathcal{O}(\|p^t - p^t_*\|) = \mathcal{O}(\delta)$ and \Cref{lem:lancaster:ndc:rd} we find:
\[
\begin{array}{rcl}
    0 & = & \omega^t_k \\
    & = & \left \langle T^t_k v^t, p^t \right \rangle \\
    & = & -2 k^c \left \langle D v^t, p^t \right \rangle \\
    & = & -2 k^* \left \langle D v^t_*, p^t_* \right \rangle + \mathcal{O}(\delta),
\end{array}
\]
and therefore,
\begin{equation}\label{eq:pDv = 0:endl:rd}
    0 = \left \langle D v^t, p^t \right \rangle =  \left \langle D v^t_*, p^t_* \right \rangle
\end{equation}
at leading order. 

Third, using \eqref{eq:rho things:ndc:rd}, we find that $\rho^t$ is of the form
\begin{equation}\label{eq:tilde omega 0:endl:rd}
    \rho^t = \tilde{\rho}\delta + \mathcal{O}(\delta^2) \quad \text{ where } \quad \tilde{\rho} = - \sqrt{-2 \left \langle F_\mu^*, p^s_* \right \rangle \left \langle F_{uu}^*(v^s_*,v^s_*),p^s_* \right \rangle \hat{\mu}^t}.
\end{equation}

Fourth, by definition \eqref{eq:T def:rd},
\begin{equation*}
    0 = \left \langle T^t v^t, p^t \right \rangle = \left \langle \tilde{T}^t v^t_*, p^t_* \right \rangle \delta + \mathcal{O}(\delta^2);
\end{equation*}
therefore, from \eqref{eq:pDv = 0:endl:rd} we conclude
\begin{equation*}
    \left \langle \tilde{F}_u^tv^t_*, p^t_* \right \rangle = 0.
\end{equation*}
Consequently, using \eqref{eq:tilde T and tilde S:endl:rd}, we obtain
\begin{equation}\label{eq:first identity in AB-system:endl:rd}
    \left\langle \left(\tilde{\mu}^s F_{u\mu}^* - F_{u\nu}^*\right)v^t_* + F_{uu}^*\left(\tilde{u}^s,v^t_*\right),p^t_*\right \rangle = - \bar{u}^t \left\langle F_{uu}^*(v^s_*, v^t_*),p^t_*\right \rangle
\end{equation}

Fifth, the $\mathcal{O}(\delta)$ expansion of $F^s = 0$ yields
\begin{equation}\label{eq:second identity in AB-system:endl:rd}
    F_u^* \tilde{u}^s + \tilde{\mu}^s F_\mu^* - F_\nu^* = 0,
\end{equation}
while the $\mathcal{O}(\delta^2)$ level gives
\begin{equation}\label{eq:third identity in AB-system:endl:rd}
    F_u^* \hat{u}^s + \hat{\mu}^s F_\mu^* + \frac{1}{2}F_{\nu\nu}^* + \frac{1}{2}F_{uu}^*(\tilde{u}^s,\tilde{u}^s) + \frac{1}{2}(\tilde{\mu}^s)^2 F_{\mu\mu}^* - F_{u\nu}^* \tilde{u}^s - \tilde{\mu}^sF_{\mu\nu}^* + \tilde{\mu}^s F^*_{u\mu}\tilde{u}^s = 0.
\end{equation}
Finally, from the $\mathcal{O}(\delta)$ expansion of $\left \langle F_u^s v^s, p^s \right \rangle = 0$ we obtain the identity
\begin{equation}\label{eq:fourth identity in AB-system:endl:rd}
    \left\langle F_{uu}^*(\tilde{u}^s,v^s_*) + \left(\tilde{\mu}^s F_{u\mu}^* - F_{u\nu}^*\right)v^s_*,p^s_*\right \rangle = 0.
\end{equation} 

\subsection{The Ginzburg-Landau approximation}
\label{subsec:ginzburg-landau derivation:rd}
As in \Cref{sec:scalar}, we begin our weakly nonlinear analysis of small-amplitude patters by deriving the classic Ginzburg-Landau approximation following the standard setup. First, we set $\nu = \nu^* - \delta$, and $\mu = \mu^t - r\varepsilon^2$ with $\varepsilon \ll \delta$. Second, based on \eqref{eq:perturbed solution:rds}, we introduce the standard Ginzburg-Landau Ansatz for a Turing bifurcation in a multi-component setting:
\begin{equation}\label{eq:Ansatz:gld:rd}
    \begin{split}
        E^0 \big[&\varepsilon^2 X_{02} + \varepsilon^3 X_{03} + \mathcal{O}\l\varepsilon^4\r \big] \\
        U_{\text{GL}}(x,t) = u^t+E \big[\varepsilon A v^t + &\varepsilon^2 X_{12} + \varepsilon^3 X_{13} + \mathcal{O}\l\varepsilon^4\r \big] + \text{c.c.} \\
        E^2 \big[&\varepsilon^2 X_{22} + \varepsilon^3 X_{23} + \mathcal{O}\l\varepsilon^4\r\big] + \text{c.c.} \\
         & \phantom{iiiii}\,\: E^3 \big[\varepsilon^3 X_{33} + \mathcal{O}\l \varepsilon^4\r \big] + \text{c.c.},
    \end{split}
\end{equation}
where, as in \Cref{sec:scalar}, $E = e^{ik^cx}$ and $A(\xi, \tau):\mathbb{R}\times\mathbb{R}^+\to \mathbb{C}$ is a scalar-valued amplitude function. However, unlike in \Cref{sec:scalar}, the unknown functions $X_0j(\xi,\tau) : \mathbb{R} \times \mathbb{R^+} \to \mathbb{R}^n$ and $ X_{ij}(\xi,\tau):\mathbb{R}\times\mathbb{R}^+\to \mathbb{C}^n$ for $i \geq 1$ are now vector-valued. Additionally, the dominant mode is also slightly modified and now takes the vector-valued form $\varepsilon e^{ik^c x} A v^t$.

Next, we determine these unknown functions by substituting \eqref{eq:Ansatz:gld:rd} into \eqref{eq:reaction-diffusion system:rds} and expanding around $(u^t,\mu^t)$. By design, the first non-trivial result appears at order $\varepsilon^2$. Specifically, at the $E^0 \varepsilon^2$-level, we obtain
\begin{equation*}
   F_u^t X_{02} + F_{uu}^t(v^t, v^t) |A|^2 - r F_\mu^t = 0,
\end{equation*}
at the $E \varepsilon^2$-level, we find
\begin{equation}\label{eq:epsilon 2 E 1:gld:rd}
    T^t X_{12} + 2ik^c D v^t A_{\xi} = 0,
\end{equation}
and at the $E^2 \varepsilon^2$-level, we have
\begin{equation*}
    \left(F_u^t - 4 (k^c)^2 D\right)X_{22} + \frac{1}{2}F_{uu}^t(v^t,v^t) A^2 = 0.
\end{equation*}
Although \eqref{eq:epsilon 2 E 1:gld:rd} is not invertible, it is nevertheless solvable since $\text{Im} T^t = \text{Ker}\left((T^t)^T\right)^\perp$ and $D v^t \in \text{Im}\left(T^t\right)$; see \eqref{eq:pDv = 0:endl:rd}. Therefore, $X_{12}$ is of the form $X_{12} = \eta_1 A_{\xi}$, where $\eta_1 \in \mathbb{C}^n$ solves
\begin{equation*}
    T^t \eta_1  = - 2 i k^c D v^t + A_1(\xi,\tau)v^t,
\end{equation*}
with $A_1(\xi,\tau)$ denoting the next order correction to the amplitude \cite{doelman2019pattern}. 

At the $E\varepsilon^3$-level, we obtain
\begin{equation*}
    A_\tau v^t = \eta_2 A_{\xi\xi} + \eta_3 A + F_{uu}^t(X_{02}, A v^t) + F_{uu}^t(X_{22},A^* v^t) + \frac{1}{2} |A|^2A F_{uuu}^t(v^t,v^t,v^t),
\end{equation*}
where $\eta_2, \eta_3 \in \mathbb{C}^n$, and $F_{uuu}^t(v^t,v^t,v^t)$ denotes the trilinear term arising from the third-order Taylor expansion of $F$. (The explicit values of the coefficients $\eta_i$ and the trilinear term are not essential for our analysis and are therefore omitted.). Consequently, by talking the dot-product with $p^t$, we retrieve the Ginzburg-Landau equation:
\begin{equation*}
    A_\tau = - \frac{1}{2}\omega_{kk}^t - r \omega_{\mu}^t A + L |A|^2A
\end{equation*}
where the Landau coefficient $L$ is given by
\begin{equation}\label{eq:full landau coefficient:rd}
    L = \left\langle - F_{uu}^t\left(\left(F_u^t\right)^{-1} F_{uu}^t\left(v^t,v^t\right),v^t\right) -\frac{1}{2} F_{uu}^t\left((F_{u}^t - 4 (k^c)^2 D)^{-1} F_{uu}^t\left(v^t,v^t\right),v^t\right)+\frac{1}{2}F_{uuu}^t\left(v^t,v^t,v^t\right), p^t \right\rangle
\end{equation}

Finally, since $F_u^t$ is nearly singular, we have
\begin{equation*}
    (p^s_*)^T (F_u^t)^{-1} = \frac{(p^s_*)^T}{\tilde{\rho}\delta} + \mathcal{O}(1), \quad \text{ and } \quad (F^t_u)^{-1} F_{uu}^t\left(v^t,v^t\right) = \frac{c}{\delta} v^s_* + \mathcal{O}(1),
\end{equation*}
for some $c \in \mathbb{R}$. Hence, using the associativity of matrix multiplication, we obtain
\begin{equation*}
\begin{array}{rcl}
    \frac{c}{\delta} + \mathcal{O}(1) & = & \left \langle (F^t_u)^{-1} F_{uu}^t\left(v^t,v^t\right),p^s_* \right \rangle \\
    & = & (p^s_*)^T (F_{u}^t)^{-1}F_{uu}^t(v^t,v^t) \\
    & = & \left((p^s_*)^T (F_{u}^t)^{-1}\right)F_{uu}^t(v^t,v^t) \\ 
    & = & \frac{\left \langle F_{uu}^*(v^t_*,v^t_*),p^s_*\right\rangle}{\tilde{\rho} \delta} + \mathcal{O}(1);
\end{array}
\end{equation*}
thus,
\begin{equation*}
    c = \frac{\left \langle F_{uu}^*(v^t_*, v^t_*), p^s_* \right \rangle}{\tilde{\rho}}.
\end{equation*}
Consequently, we find
\begin{equation}
\label{eq:leading order gl:gld:rd}
    A_\tau = \left(-\frac{1}{2} \omega_{kk}^* + \mathcal{O}(\delta)\right) A_{\xi\xi} + \frac{r}{\delta}\left(- \tilde{\omega}^t_\mu + \mathcal{O}(\delta)\right) A + \frac{1}{\tilde{\rho}\delta} \left(- \left\langle F_{uu}^*(v^s_*,v^t_*),p^t_*\right \rangle \left \langle F_{uu}^*(v^t_*,v^t_*),p^s_*\right\rangle +\mathcal{O}(\delta) \right) |A|^2A.
\end{equation} 

\subsection{The derivation of the AB-system}\label{subsec:AB-system derivation:rd}
Following the now-familiar setup, we set $\nu = \nu^* - \delta$, $\mu = \mu^t - r \delta^2$ and replace the standard Ginzburg-Landau Ansatz with its modified version
\begin{equation}\label{eq:Ansatz:abd:rd}
    \begin{split}
    U_{AB}(x, t) = u^* + E^0 \big[\delta (\tilde{u}^s + v^{s}_* B) + \delta^{3/2}X_{02} + &\delta^{2} (\hat{u}^s + X_{03}) +\delta^{5/2} X_{04}+\mathcal{O}\l\delta^3\r\big] \\
    + E \big[\delta v^{t}_* A + \delta^{3/2} X_{12} + &\delta^{2} X_{13} \phantom{aaaaai}+ \delta^{5/2} X_{14} +\mathcal{O}\l\delta^3\r \big] + \text{c.c.} \\
    +E^2 \big[&\delta^2 X_{23} \phantom{aaaaai} + \delta^{5/2} X_{24} + \mathcal{O}\l\delta^3\r\big] + \text{c.c.} \\
     &\phantom{iiiiii00000000000iii9i+E_c^3}  \phantom{aaaaai}+ \text{h.o.t.},
\end{split}
\end{equation}
where, again, $E = e^{i k^* x}$, $A(\xi, \tau):\mathbb{R}\times\mathbb{R}^+ \to \mathbb{C}$ and $B(\xi, \tau): \mathbb{R} \times \mathbb{R}^+ \to \mathbb{C}$ are the two unknown scalar amplitude functions, and $X_{ij}(\xi,\tau):\mathbb{R}\times \mathbb{R}^+ \to \mathbb{R}^n$ are vector-valued functions. 

Substituting \eqref{eq:Ansatz:abd:rd}into \eqref{eq:reaction-diffusion system:rds} and expanding yields
\begin{equation*}
    \begin{split}
        F(U_{AB}; \mu^t - r \delta^2, \nu) &= F^* + F^*_u(U_{AB}-u^*) - \delta F^*_{\nu}  + (\mu^t - r \delta^2 - \mu^*)F_{\mu}^* +  \frac{1}{2}(\mu^t - r\delta^2 - \mu^*)^2F_{\mu\mu}^* \\
        &+ \frac12 F_{uu}^*(U_{AB}-u^*,U_{AB}-u^*) + \frac{1}{2}\delta^2 F_{\nu\nu}^* - \delta F_{u\nu}^*(U_{AB}-u^*) \\
        &- \delta (\mu^t - r\delta^2 - \mu^*) F_{\mu \nu}^* + (\mu^t-r\delta^2 - \mu^*) F_{u\mu}^*(U_{AB}-u^*) + \text{h.o.t.}
    \end{split}
\end{equation*}
At the $\delta$-level we obtain
\begin{equation*}
    F_u^* \tilde{u}^s + \tilde{\mu}^s F_\mu^* - F_\nu^* + E T^* v^t_*A = 0,
\end{equation*}
which simplifies to $0 = 0$ by \eqref{eq:second identity in AB-system:endl:rd}. Consequently, the first non-trivial result appears at the $E^0 \delta^{3/2}$-level:
\begin{equation*}
    T^* X_{12} + 2 i k^* D v^t_* A_{\xi} = 0,
\end{equation*}
which is solvable by \eqref{eq:pDv = 0:endl:rd} (cf. \eqref{eq:epsilon 2 E 1:gld:rd}). Next, at the $E^0 \delta^2$-level we find
\begin{equation*}
        \begin{split}
        B_{\tau} v^s_* &= F_u^*(\hat{u}^s+X_{03}) + (\hat{\mu}^s+\hat{\mu}^t - r)F_{\mu}^* + \frac12 (\tilde{\mu}^s)^2 F_{\mu\mu}^* + \frac{1}{2}F_{uu}^*(\tilde{u}^s + v^s_* B, \tilde{u}^s + v^s_* B) \\ & + |A|^2 F_{uu}^*(v^t_*,v^t_*) + \frac12 F^*_{\nu\nu} - \tilde{\mu}^s F_{\mu\nu}^* - F_{u\nu}^* (\tilde{u}^s + v^s_* B) + \tilde{\mu}^s F_{u\mu}^*(\tilde{u}^s + v^s_* B) + D v^s_* B_{\xi\xi},\\
    \end{split}
\end{equation*}
which -- by applying \eqref{eq:third identity in AB-system:endl:rd}, taking the inner product with $p^s_*$, and then using \eqref{eq:fourth identity in AB-system:endl:rd} -- simplifies to
\begin{equation*}
    B_\tau = - \frac{1}{2}\rho_{kk}^* B_{\xi\xi} + (\hat{\mu}^t-r) \left\langle F_\mu^*, p^s_* \right \rangle + \frac{1}{2} \left \langle F_{uu}^*(v^s_*,v^s_*),p^s_* \right\rangle B^2 + \left \langle F_{uu}^*(v^t_*,v^t_*),p^s_*\right\rangle |A|^2.
\end{equation*}
Finally, the $E \delta^2$-level yields
\begin{equation*}
    A_{\tau} v^t_* = \eta A_{\xi\xi} + \left(\left(\tilde{\mu}^sF_{u\mu}^* - F_{u\nu}^*\right)v^t_* + F_{uu}^*(\tilde{u}^s,v^t_*)\right) A + F_{uu}^*(v^s_*,v^t_*)AB + T^* X_{13}.
\end{equation*}
for some $\eta \in \mathbb{R}^n$. Hence, by taking the inner product with $p^t_*$ and using \eqref{eq:first identity in AB-system:endl:rd} we arrive at:
\begin{equation*}
    A_{\tau} = -\frac{1}{2}\omega^*_{kk} A_{\xi\xi} - \bar{u}^t\left\langle F_{uu}^*(v^s_*,v^t_*),p^t_*\right\rangle A + \left\langle F_{uu}^*(v^s_*,v^t_*),p^t_* \right \rangle AB.
\end{equation*}
In conclusion, the AB-system for \eqref{eq:reaction-diffusion system:rds} is given by:
\begin{equation}
    \begin{cases}\label{eq:coupled system:absd:rd}
        A_{\tau} = -\frac{1}{2}\omega^*_{kk} A_{\xi\xi} - \bar{u}^t\left\langle F_{uu}^*(v^s_*,v^t_*),p^t_*\right\rangle A + \left\langle F_{uu}^*(v^s_*,v^t_*),p^t_* \right \rangle AB \\
        B_\tau = - \frac{1}{2}\rho_{kk}^* B_{\xi\xi} + (\hat{\mu}^t-r) \left\langle F_\mu^*, p^s_* \right \rangle + \frac{1}{2} \left \langle F_{uu}^*(v^s_*,v^s_*),p^s_* \right\rangle B^2 + \left \langle F_{uu}^*(v^t_*,v^t_*),p^s_*\right\rangle |A|^2.
    \end{cases}
\end{equation}
As already discussed in \Cref{subsec:rescaling the ab system:gs:s}, the derivation process allows the rescaling of \eqref{eq:coupled system:absd:rd} into the simpler form studied in \Cref{sec:coupled system} (see \eqref{eq:canonical form:csd:gs}) via:
\begin{equation}
\label{eq:rescalingsAB-syst}
    A = \tilde{A}, \quad B = \bar{u}^t \tilde{B}, \quad R = \frac{r}{\hat{\mu}^t}, \quad \tilde{\tau} = - \bar{u}^t \left \langle F_{uu}^*(v^s_*, v^t_*),p^t_* \right \rangle \tau, \quad \text{ and } \quad \tilde{\xi} = \sqrt{\frac{2 \bar{u}^t \left\langle F_{uu}^*(v^s_*,v^t_*),p^t_*\right \rangle}{\omega_{kk}^*}}\xi,
\end{equation}
so that:
\begin{equation}
\label{rescalealphadbeta-Ncomp}
    \alpha = \frac{ \left \langle F_{uu}^*(v^s_*,v^s_*),p^s_* \right \rangle}{2 \left \langle F_{uu}^*(v^s_*,v^t_*),p^t_* \right \rangle}, \quad d = \frac{2 \rho_{kk}^* \left \langle F_{uu}^*(v^s_*,v^t_*),p^t_* \right \rangle}{\omega_{kk}^* \left\langle F_{uu}^*(v^s_*,v^s_*), p^s_* \right \rangle}, \quad \text{ and } \quad \beta = \frac{\left \langle F_{uu}^*(v^t_*,v^t_*),p^s_*\right \rangle}{\hat{\mu}^t \left \langle F_{\mu}^*,p^s_* \right \rangle}
\end{equation}
It follows from \eqref{eq:reactiondiffusion-TuringFold} and \eqref{eq:reactiondiffusionnondeg-TuringFold} that $\beta$ and the $|A|^2$ coefficient in \eqref{eq:coupled system:absd:rd} have the same sign. Additionally, as in the higher-order scalar models of \Cref{subsec:general setting:se,subsec:rescaling the ab system:gs:s}, we have $\alpha d = \rho_{kk}^*/\omega_{kk}^*$ (cf. \eqref{eq:rescaleABalphaetc} and \eqref{eq:scalingsalphaetc-gen}).

\subsection{Reducing the AB-system to the Ginzburg-Landau system}
\label{subsec:reducing the AB-system back to the ginzburg landau system:rd}
In \Cref{subsec:ginzburg-landau derivation:rd} (see \eqref{eq:leading order gl:gld:rd}), we derived the Landau coefficient:
\begin{equation}\label{eq:landau coef:rabtgl:rd}
L = -\left\langle F_{uu}^*(v^s_*,v^t_*),p^t_*\right \rangle \left \langle F_{uu}^*(v^t_*,v^t_*),p^s_*\right\rangle \frac{1}{\tilde{\rho} \delta} + \mathcal{O}(1),
\end{equation}
and in \Cref{subsec:AB-system derivation:rd} (see \eqref{eq:coupled system:absd:rd}), we derived the $|A|^2$ coefficient:
\begin{equation}\label{eq:beta:rabtgl:rd}
    \beta = \left \langle F_{uu}^*(v^t_*,v^t_*),p^s_*\right\rangle.
\end{equation}
Both these coefficients arise in the context of the Turing-fold, as outlined in the beginning of \Cref{sec:reaction diffusion system}, for systems of the form \eqref{eq:reaction-diffusion system:rds}.  By combining \eqref{eq:reactiondiffusionnondeg-TuringFold}, \eqref{eq:values of variables:ndc:rd}, and \eqref{eq:tilde omega 0:endl:rd}, we see that the Landau coefficient \eqref{eq:landau coef:rabtgl:rd} and the $|A|^2$ coefficient \eqref{eq:beta:rabtgl:rd} are related via:
\begin{equation*}
   L = \frac{\tilde{\omega}_\mu^t}{\left\langle F_\mu^*,p^s_*\right \rangle \delta} \beta + \mathcal{O}(1).
\end{equation*}
Consequently, for sufficiently small $\delta > 0$ and non-zero $\beta$, the $L$- and $\beta$ coefficients have opposite signs.
   
As discussed in \Cref{subsubsec:Reducing the AB-system:gs:se}, the AB-system serves as a "wider lens" for describing the dynamics resulting from the Turing-fold interaction. Therefore, when zooming into the parameter region where both the Ginzburg-Landau and AB-system are valid, they should describe the same dynamics. This equivalence is (again) confirmed by applying the following (re)scalings:
    \begin{equation*}
        r = \frac{\varepsilon^2}{\delta^2} \tilde{r}, \quad B = \bar{u}^t + \frac{\varepsilon^2}{\delta} \tilde{B}, \quad A = \frac{\varepsilon}{\delta} \tilde{A}, \quad \xi = \frac{\varepsilon}{\sqrt{\delta}}\tilde{\xi}, \quad \text{ and } \quad \tau = \frac{\varepsilon^2}{\delta} \tilde{\tau},
    \end{equation*}
to \eqref{eq:coupled system:absd:rd}, yielding:
    \begin{equation*}
        \begin{cases}
            \varepsilon^3 \tilde{A}_{\tilde{\tau}} = - \frac{1}{2}\omega_{kk}^* \varepsilon^3 \tilde{A}_{\tilde{\xi}\tilde{\xi}} + \varepsilon^3\left\langle F_{uu}^*(v^s_*,v^t_*),p^t_*\right\rangle \tilde{A}\tilde{B} + \mathcal{O}(\varepsilon^4)\\
            0 = \varepsilon^2 \left\langle F_{uu}^*(v^s_*,v^s_*),p^s_*\right\rangle \bar{u}^t \tilde{B} - \frac{\varepsilon^2}{\delta} \tilde{r} \left\langle F_{\mu}^*, p^s_* \right \rangle + \frac{\varepsilon^2}{\delta} \left\langle F_{uu}^*(v^t_*,v^t_*),p^s_*\right\rangle|\tilde{A}|^2 + \mathcal{O}(\varepsilon^3).
        \end{cases}
    \end{equation*}
Consequently, by solving the second equation for $\tilde{B}$ and substituting the result in the first equation we recover the leading order part of Ginzburg-Landau equation \eqref{eq:leading order gl:gld:rd},
   \begin{equation*}
    \tilde{A}_{\tilde{\tau}} = -\frac{1}{2}\omega_{kk}^* \tilde{A}_{\tilde{\xi}\tilde{\xi}} - \frac{\tilde{r}}{\delta} \tilde{\omega}^t_\mu  \tilde{A} - \left\langle F_{uu}^*(v^s_*,v^t_*),p^t_*\right \rangle \left \langle F_{uu}^*(v^t_*,v^t_*),p^s_*\right\rangle \frac{1}{\tilde{\rho} \delta} |\tilde{A}|^2\tilde{A} + \mathcal{O}(\varepsilon).
\end{equation*}

\section{Periodic patterns in the AB-system}\label{sec:coupled system}
In \Cref{subsubsec:coupled system derivation:gs}, we sellected the following canonical form for the AB-system:
\begin{equation}\label{eq:canonical form:cs}
    \begin{cases}
        \phantom{\frac{1}{\alpha}} A_\tau = A_{\xi\xi} + A - AB \\
        \frac{1}{\alpha} B_\tau = d B_{\xi\xi} + 1 - R - B^2 + \beta |A|^2,
    \end{cases}
\end{equation}
where $\alpha, d > 0$, $\beta \in \mathbb{R}$, and $R \in \mathbb{R}$ is the bifurcation parameter (cf. \eqref{eq:canonical form:csd:gs}). Note that $\beta$ can also be removed from  \eqref{eq:canonical form:cs} by the additional rescaling $\tilde{A} = \sqrt{|\beta|}A$; however, this would necessitate distinguishing between the cases $\beta = \pm 1$ and $\beta = 0$. Consequently, we prefer to analyze the AB-system in the form \eqref{eq:canonical form:cs}. By construction, $R = 0$ corresponds to the Turing bifurcation, $R = 1$ to the saddle-node bifurcation (and thus to the tipping point), and $|R| \ll 1$ to the standard Ginzburg-Landau setting for the Turing bifurcation. 

Naturally, we are interested in the dynamics generated by \eqref{eq:canonical form:cs}. Here, we focus on the existence and stability of the simplest (and perhaps the most relevant) solutions of \eqref{eq:canonical form:cs}: the homogeneous states and the spatially periodic plane waves. Specifically, we consider solutions to \eqref{eq:canonical form:cs} of the form
\begin{equation}
\label{eq:periodic solutions:cs}
A(\xi,\tau) = \bar{A} e^{iK\xi} \quad \text{ and } \quad B(\xi,\tau) = \bar{B}, 
\end{equation}
with $\bar{A} > 0$, $\bar{B},  K \in \mathbb{R}$, and $(\bar{A}, \bar{B}) = (\bar{A}(K,R), \bar{B}(K,R))$. These solutions directly correspond to their counterparts (of the same form) in the Ginzburg-Landau equation. In fact, in the Ginzburg-Landau limit, the stable subfamily of these planar waves (parameterized by $R$) forms {\it the nose} of the Busse balloon, which opens/closes at the Turing bifurcation as $R$ passes through $0$. This `nose' is locally bounded in the $(K,R)$-plane by the Eckhaus parabola (for $|R| \ll 1$) -- see \cite{aranson2002world,busse1978non,cross1993pattern,eckhaus2012studies,mielke2002ginzburg}. Thus, the question whether stable spatially periodic patterns originating from a Turing bifurcation persist (and remain stable) beyond the saddle-node point at which the underlying spatially homogeneous state disappears -- i.e., {\it Can tipping be evaded by pattern formation?} \cite{rietkerk2021evasion} -- reduces to whether \eqref{eq:canonical form:cs} exhibits patterns of the form \eqref{eq:periodic solutions:cs} that appear as stable patterns as $R$ increases through $0$ and remain so as $R$ passes through $R=1$.

\subsection{The existence of `planar' spatially periodic patterns}
\label{subsec:existence:cs}
\begin{figure}[t]
    \centering
\begin{subfigure}{0.49\textwidth}
    \includegraphics[width=\textwidth]{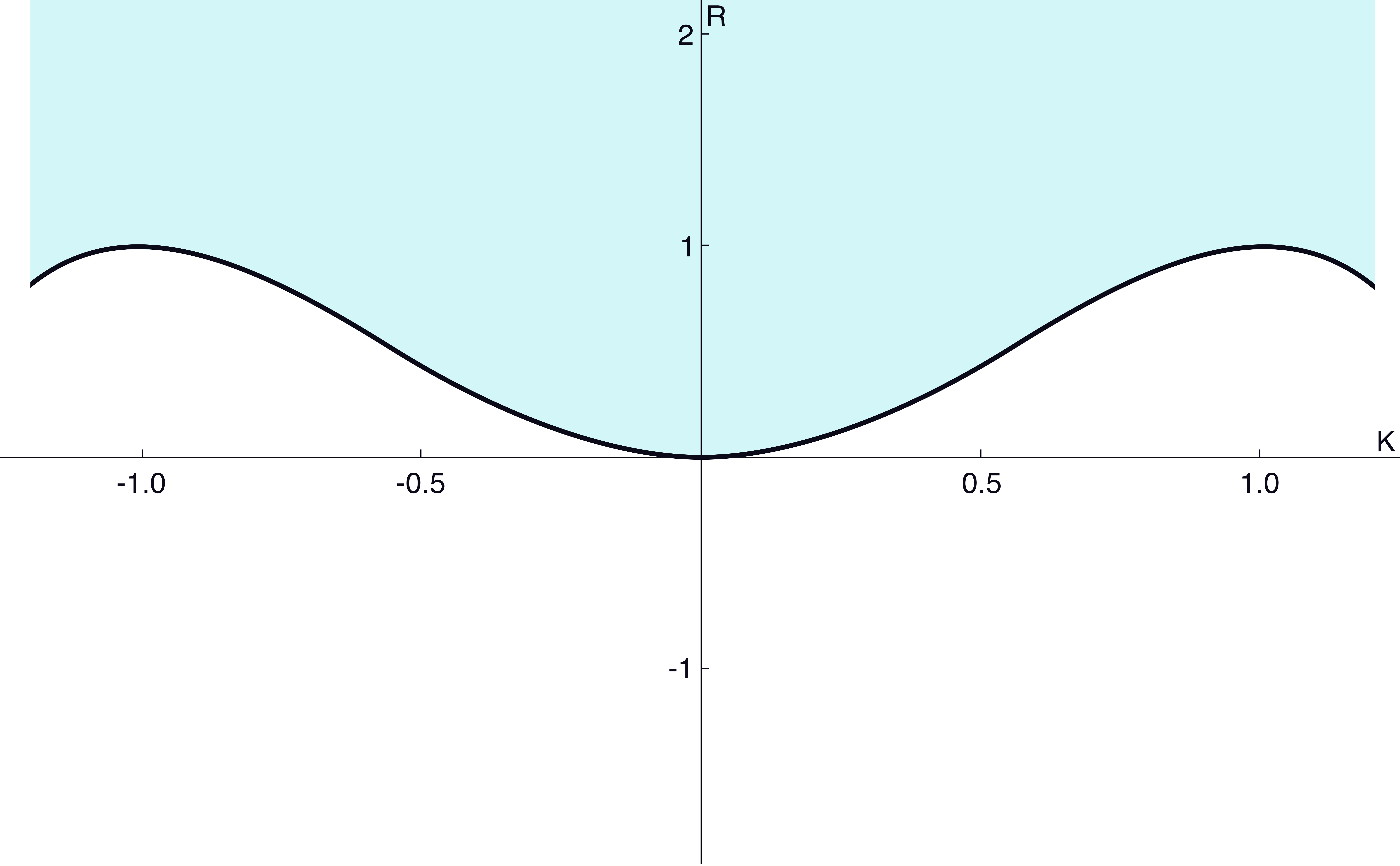}
    \caption{$\beta > 0$}
    \label{fig:beta positive:cs}
\end{subfigure}
\hfill
\begin{subfigure}{0.49\textwidth}i
    \includegraphics[width=\textwidth]{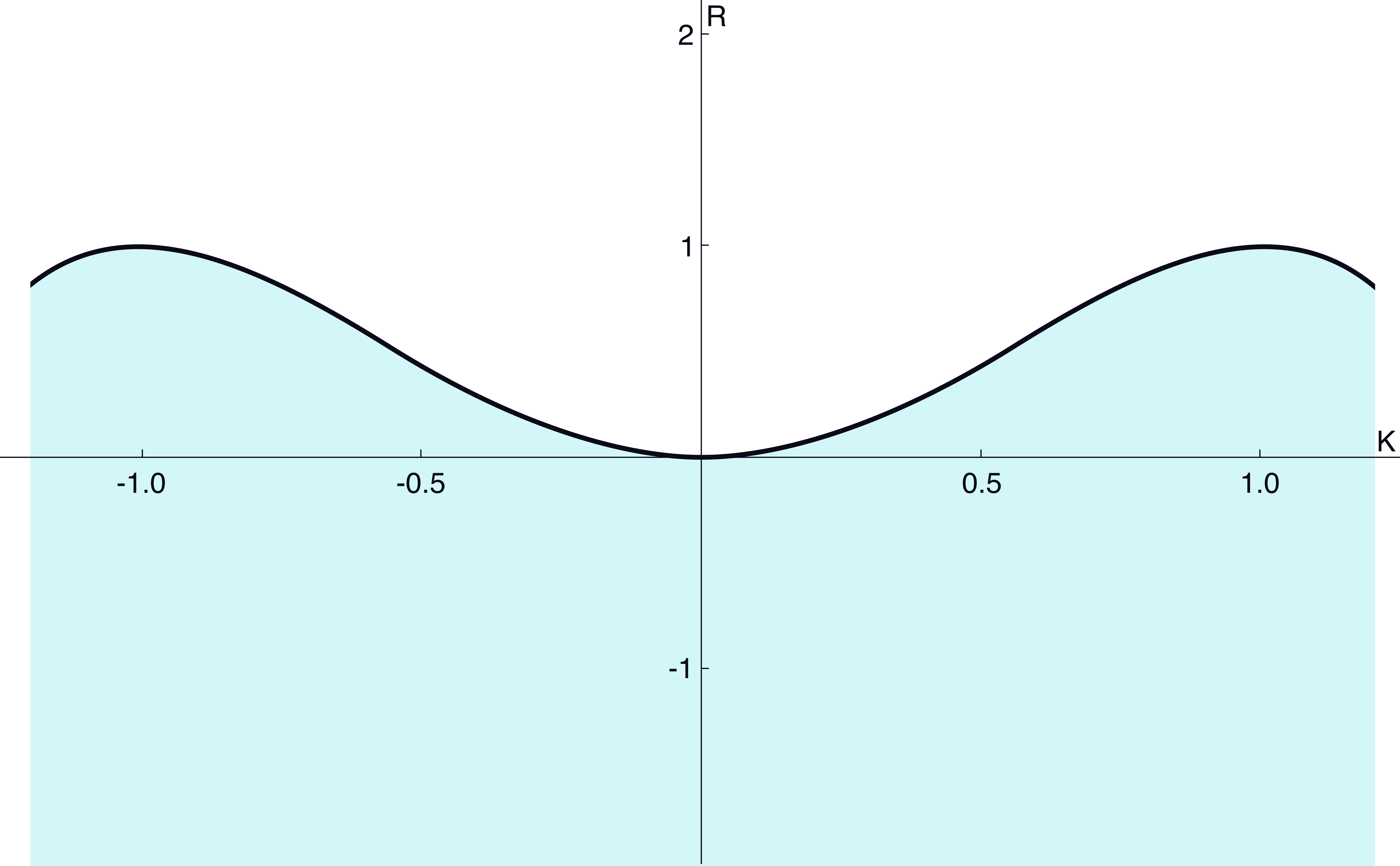}
    \caption{$\beta < 0$}
    \label{fig:beta negative:cs}
\end{subfigure}
    \caption{The regions (in blue) in the $(K,R)$-plane in which planar waves of the structure \eqref{eq:periodic solutions:cs} exist in the AB-system \eqref{eq:canonical form:cs}: $\beta > 0$ on the left and $\beta < 0$ on the right. Note that the boundary of the domains is in both cases given by the curve $ R_e(K)$ \eqref{defReK}.}
    \label{fig: Existence regions AB systems}
\end{figure}

Substituting \eqref{eq:periodic solutions:cs} into \eqref{eq:canonical form:cs} yields:
\begin{equation}\label{eq:barA and barB eq:ex:cs}
\begin{cases}
    0 = \left(-K^2 + 1 - \bar{B}\right)\bar{A} \\
    0 = 1 - R - \bar{B}^2 + \beta \bar{A}^2.
\end{cases}
\end{equation}
Clearly, \eqref{eq:barA and barB eq:ex:cs} has two kinds of solutions $(\bar{A}, \bar{B})$. First, the homogeneous states:
\begin{equation}
\label{eq: CS saddle-node solution}
\l A_s(R), B^{\pm}_s(R)\r = \l 0, \pm \sqrt{1 - R}\r,
\end{equation}
which are associated with the saddle-node bifurcation and (for $R < 1$) correspond to the $\mathcal{O}(\delta)$ coefficient in the Taylor expansion of $u^\pm$ (cf. \eqref{AsrBsr:gs}, \eqref{rs-unsc} for the equivalent unscaled expressions). Second, the $K$-parameterized family of planar waves:
\begin{equation}
\label{eq: CS periodic solutions}
\l A_p(K, R), B_p(K, R)\r = \l \sqrt{\frac{1}{\beta}\left((1-K^2)^2 + R - 1\right)}, 1 - K^2 \r,
\end{equation}
which originates from the Turing bifurcation. Thus, the spatially periodic patterns $(A_p(K, R), B_p(K, R))$ exist for $R \geq 1 - (1-K^2)^2$ if $\beta > 0$ and for $R \leq 1 - (1-K^2)^2$ if $\beta < 0$. We denote the boundary of the existence domain by the curve $R_e=R_e(K)$:
\begin{equation}
\label{defReK}
    R_e(K) = 1 - (1-K^2)^2.
\end{equation}
Observe that $R_e(K) \leq 1$, which for $\beta > 0$ implies the existence of spatially periodic patterns (for all wave numbers $K \in \mathbb{R}$) beyond the tipping point $R=1$; conversely, for $\beta < 0$ no planar waves of the type \eqref{eq:periodic solutions:cs} exists for $R > 1$, as illustrated in \Cref{fig: Existence regions AB systems}.

\subsection{Stability of the spatially periodic patterns}
\label{subsec:stability:cs}
To analyze the stability of these patterns, we follow the standard Ginzburg-Landau procedure by substituting
\begin{equation*}
A(\xi,\tau) = \left(\bar{A} + a(\xi,\tau) \right)e^{i K \xi} \quad \text{ and } \quad
B(\xi, \tau) = \bar{B} + b(\xi,\tau)
\end{equation*}
into \eqref{eq:canonical form:cs}, and linearizing:
\begin{equation*}
    \begin{cases}
        a_\tau = a_{\xi\xi} + 2iKa_\xi + a(1-K^2)   - (\bar{B}a + \bar{A}b) \\
        b_\tau = d\alpha b_{\xi\xi} +\alpha \left( -2 \bar{B} b + \beta \bar{A}(a + a^*)\right).
    \end{cases}
\end{equation*}
Next, we decompose $a$ into its real and imaginary parts, $a = U + i V$, and set $b = W$:
\begin{equation*}
    \begin{cases}
        U_\tau = U_{\xi\xi} - 2K V_{\xi} +(1-K^2-\bar{B})U  - \bar{A} W \\
        V_\tau = V_{\xi\xi} + 2 K U_{\xi} +(1-K^2-\bar{B})V \\
        W_\tau = d \alpha W_{\xi\xi} + 2\alpha (\beta \bar{A} U - \bar{B} W ).\\
    \end{cases}
\end{equation*}
Finally, by assuming that the solutions are of the form
\begin{equation*} 
    \begin{pmatrix}
        U \\ V \\ W
    \end{pmatrix} =
        e^{ik \xi + \lambda \tau}
    \begin{pmatrix}
        u \\ v \\ w
    \end{pmatrix}
\end{equation*}
we obtain the following spectral problem:
\begin{equation}\label{eq: CS eigenvalueproblem for stability}
    \lambda
    \begin{pmatrix}
        u \\ v \\ w
    \end{pmatrix} =
    \begin{pmatrix}
        1 - K^2 - \bar{B} - k^2 & - 2i k K & - \bar{A} \\
        2 i k K  & 1 - K^2 - \bar{B} - k^2 & 0 \\
        2 \alpha \beta \bar{A} & 0 & - 2 \alpha \bar{B} - d \alpha k^2
    \end{pmatrix}
    \begin{pmatrix}
        u \\ v \\ w
    \end{pmatrix}.
\end{equation}

\subsubsection{Stability of the homogeneous states $(A_s,B_s^+)$ and $(A_s,B_s^-)$}
\label{subsubsec: Analysis of fold branches}
\begin{figure}[t]
    \centering
    \includegraphics[width=\textwidth]{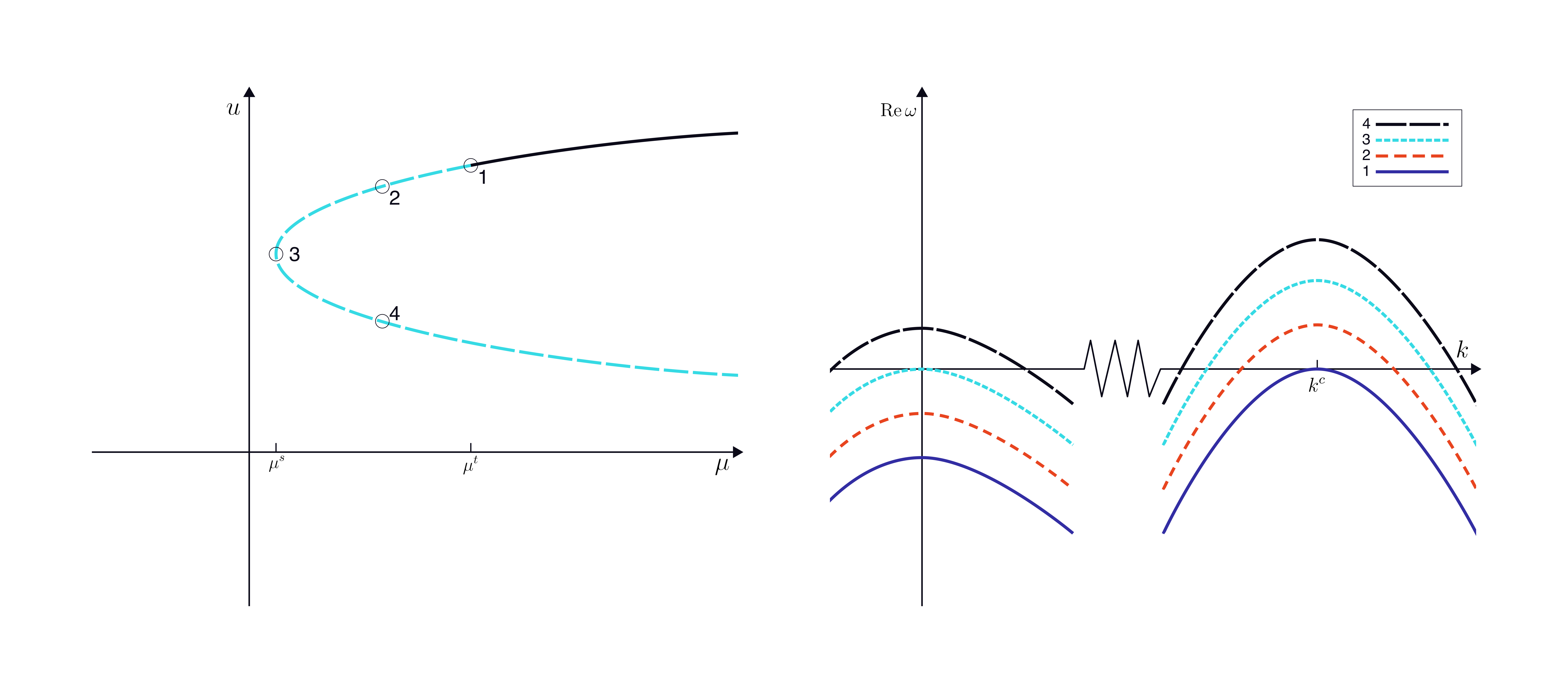}
    \caption{In the right figure, the $\text{Re} \, \omega$ dispersion curve is plotted as a function of $k$ for various values of $(u^\pm,\mu)$. The left figure shows the corresponding values of $(u^\pm, \mu)$.}
    \label{fig: Omega curves}
\end{figure}

For the homogeneous states $(A_s,B_s^+)$ and $(A_s,B_s^-)$ \eqref{eq: CS saddle-node solution}, the eigenvalue problem \eqref{eq: CS eigenvalueproblem for stability} simplifies to:
\begin{equation*}
    \lambda
    \begin{pmatrix}
        u \\ v \\ w
    \end{pmatrix} =
    \begin{pmatrix}
    1 \mp \sqrt{1-R} - k^2 & 0 & 0 \\
    0 & 1 \mp \sqrt{1-R} - k^2 & 0 \\
    0 & 0 & \mp 2 \alpha \sqrt{1-R} - d \alpha k^2
    \end{pmatrix}
    \begin{pmatrix}
        u \\ v \\ w
    \end{pmatrix},
\end{equation*}
yielding the eigenvalues:
\begin{equation*}
    \lambda_1(k)^\pm = \lambda_2^\pm(k) = 1 \mp \sqrt{1-R} - k^2 \quad \text{ and } \quad \lambda_3^\pm = \mp 2\alpha \sqrt{1-R} - d \alpha k^2.
\end{equation*}
Hence: 
\begin{itemize}
    \item For $R < 0$, $(A_s,B_s^+)$ is stable, while $(A_s,B_s^-)$ is unstable.
    \item For $R\in(0,1)$, both $(A_s,B_s^+)$ and $(A_s,B_s^-)$ are unstable.
\end{itemize}
These stability results are fully consistent with the original setup. At the saddle-node bifurcation point ($R = 1$), both trivial states are unstable; additionally, as $R$ decreases through $0$, $u^+$ regains stability via a Turing bifurcation. In \Cref{fig: Omega curves}, we indicate the (in)stability of the corresponding $u^\pm$-states as a function of $\mu$ by plotting the dispersion curves $\text{Re} \, \omega(k)$ near $k=0$ and $k=k^c$.

\subsubsection{Stability of the periodic states $(A_p, B_p)$ }
The eigenvalue problem for $(A_p, B_p) = (A_p(K,R),B_p(K,R))$ \eqref{eq: CS periodic solutions} is:
\begin{equation}
    \lambda
    \begin{pmatrix}
        u \\ v \\ w
    \end{pmatrix} =
    \begin{pmatrix}
        - k^2 & - 2i k K & - A_p \\
        2 i k K  & - k^2 & 0 \\
        2 \alpha \beta A_p & 0 & - 2 \alpha B_p - d\alpha k^2
    \end{pmatrix}
    \begin{pmatrix}
        u \\ v \\ w
    \end{pmatrix} =: \mathcal{A}(k) \begin{pmatrix}
        u \\ v \\ w
    \end{pmatrix}.
\end{equation}
This system has a cubic characteristic polynomial:
\begin{equation}
\label{def:Qklambda}
    Q(k,\lambda) = \det(\mathcal{A}(k) - \lambda I),
\end{equation}
which for any $k \in \mathbb{R}$ has three zeroes $\lambda_j(k) = \lambda_j(k; K,R)$, $j=1,2,3$, and satisfies
\begin{equation}
\label{detAk}
Q(k,0) = \det(\mathcal{A}(k)) = - \alpha k^2 \left[\l 2B_p + dk^2\r\l k^2-4K^2\r + 2 \beta A_p^2 \right].
\end{equation}

We will approach the stability problem in three steps. First, we consider the `trivial' ODE case $k \equiv 0$. Second, we assume $|k|\ll 1$ to study the sideband stability. Finally, we consider $k \in \mathbb{R}$, which corresponds to the Turing stability of the trivial state $(U,V,W) \equiv (0,0,0)$ of the linearized system \eqref{eq: CS eigenvalueproblem for stability}.

\paragraph{$\bullet$ ODE stability.} For $k = 0$, the characteristic polynomial simplifies to
\[
Q(0, \lambda) = -\lambda\l\lambda^2 + 2\alpha \l1-K^2\r\lambda + 2\alpha \beta A_p^2\r;
\]
hence, the eigenvalues are: 
\[
\lambda_1(0) = 0, \quad \lambda_{2,3} = \lambda_{\pm}(0) = -\alpha\left(\left(1-K^2\right)\pm \sqrt{\left(1-K^2\right)^2-2\frac{\beta}{\alpha}A_p^2}\right).
\]
Since $\alpha > 0$, we conclude that the special ODE case ($k=0$) imposes strong conditions on the potential stability of the $(A_p(K,R),B_p(K,R))$-patterns: $\beta > 0$, and $-1 < K < 1$. Particularly, the AB-system \eqref{eq:canonical form:cs} does not have stable spatially periodic patterns of type \eqref{eq:periodic solutions:cs} if $\beta < 0$ -- i.e., $\beta < 0$ corresponds to the subcritical Turing bifurcation in the Ginzburg-Landau approximation. 

From now on we assume that the $(A_p,B_p)$-solutions are ODE stable, i.e. $\beta > 0$ and $|K| < -1$.

\paragraph{$\bullet$ Sideband stability.}
We assume that $k^2 = \varepsilon$ with $0 < \varepsilon \ll 1$ sufficiently small. (This $\varepsilon$ is of course unrelated to the $\varepsilon$ used in the derivation of the Ginzburg-Landau approximations.) This ensure that $\lambda_{2,3}(\sqrt{\varepsilon}) < 0$, since $\lambda_{2,3}(0) < 0$. Therefore, since $\lambda_1(0)=0$, the (sideband) stability of $(A_p,B_p)$ is determined by $\lambda_1(k)$. Consequently, we set $\lambda_1 = \varepsilon \tilde{\lambda}$, which simplifies the characteristic polynomial \eqref{def:Qklambda}:
\[
    Q\l \sqrt{\varepsilon}, \varepsilon \tilde{\lambda}\r = \varepsilon \l 8\alpha B_p K^2 - 2 \alpha \beta A_p^2\l 1 + \tilde{\lambda}\r\r + \mathcal{O}\l\varepsilon^2\r;
\]
hence,
\[
    \tilde{\lambda} = \frac{8 \alpha B_p K^2}{2 \alpha \beta A_p^2} - 1 = \frac{- 5 K^4 + 6 K^2 - R}{\beta A_p^2}.
\]
Consequently, if $\beta > 0$, $K^2 < 1$, and $R > -5 K^4 + 6 K^2$, then $(A_p,B_p)$ is sideband stable. We denote this sideband stability curve by $R_s = R_s(K)$:
\begin{equation}
\label{defRsK}
    R_s(K) = -5 K^4 + 6 K^2.
\end{equation}
It is worth noting that this (sideband) stability result also corroborates the Ginzburg-Landau nature of the $R \to 0$ limit in the AB-system. To see this, we define $K_e = K_e(R)  > 0$ and $K_s = K_s(R) > 0$ such that $R_e(K_e) = R$ \eqref{defReK} and $R_s(K_s) = R$ \eqref{defRsK}. For $0 < R \ll 1$, we find at leading order in $R$ that $K_e(R) = \sqrt{R/2}$ and $K_s(R) = \sqrt{R/6}$; thus:
\[
\lim_{R \downarrow 0} \frac{K_s(R)}{K_e(R)} = \frac{1}{\sqrt{3}},
\]
which recovers the (classic) $1:1/\sqrt{3}$ Ginzburg-Landau ratio between the width of the existence $K$-interval and the stable $K$-subinterval, as determined by the Eckhaus parabola \cite{eckhaus2012studies}.

\paragraph{$\bullet$ Turing stability.}
Although the Turing instability is a generic instability mechanism of spatially periodic solutions of systems of reaction-diffusion equations -- and thus of the Ginzburg-Landau equation -- on the real line, it does not appear as a boundary of the family of stable spatially periodic `planar waves' in the real Ginzburg-Landau equation \cite{rademacher2007instabilities}, i.e., a Ginzburg-Landau equation with real coefficients (like \eqref{eq:ginzburg-landau:a model}, \eqref{eq:ginzburg-landau with deltas:gs}, \eqref{eq:leading order gl:gld:rd}). However, it does occur as an instability mechanism in the complex Ginzburg-Landau equation \cite{matkowsky1993stability} (i.e., Ginzburg-Landau equation with complex coefficients), and, as we will show here, also in the AB-system.

\begin{figure}[tp]
    \centering
    \begin{subfigure}{0.8\textwidth}
    \includegraphics[width=\textwidth]{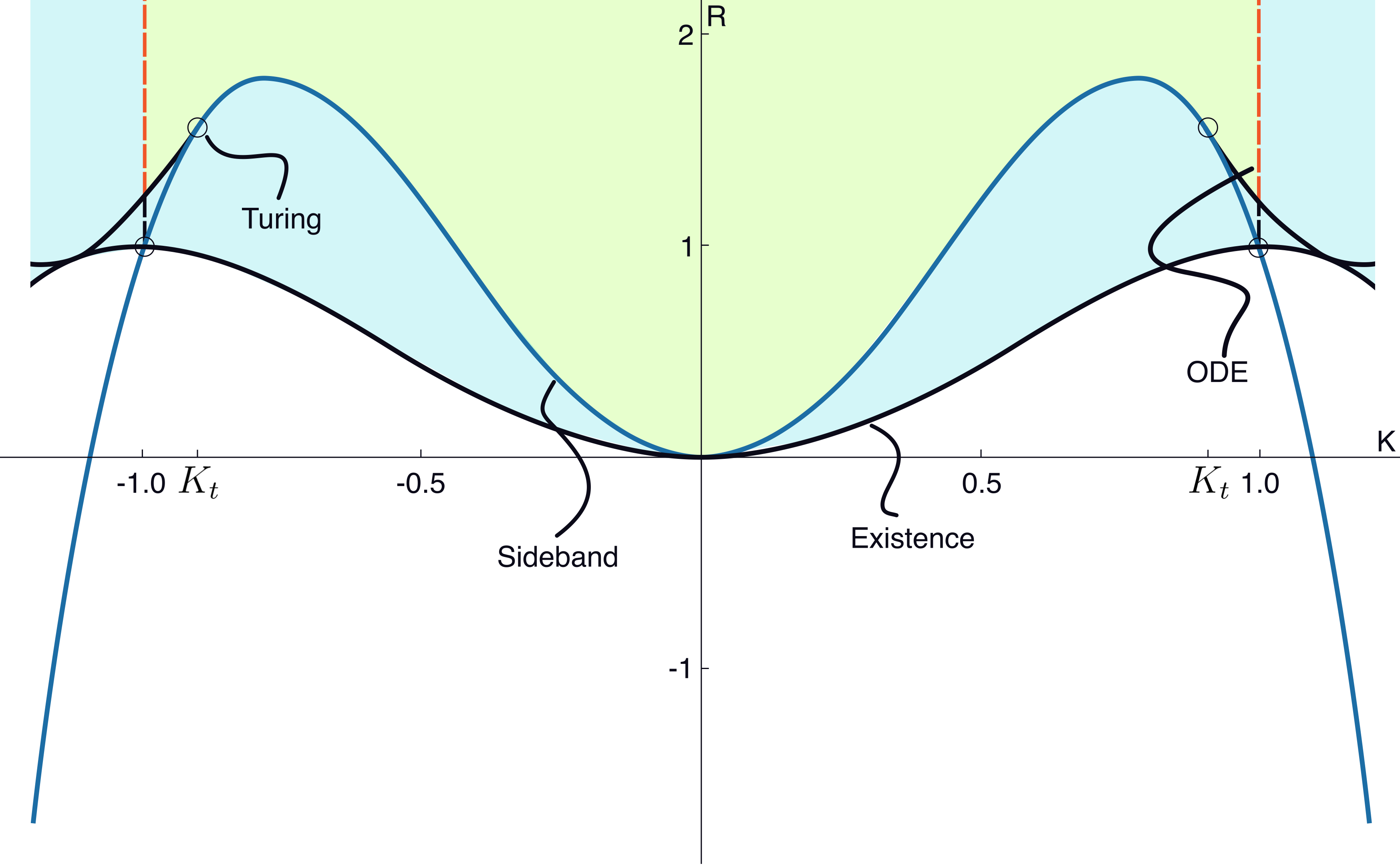}
    \caption{$d = 0.1$}
    \label{fig:small d:cs}
\end{subfigure}
\hfill
\begin{subfigure}{0.8\textwidth}
    \includegraphics[width=\textwidth]{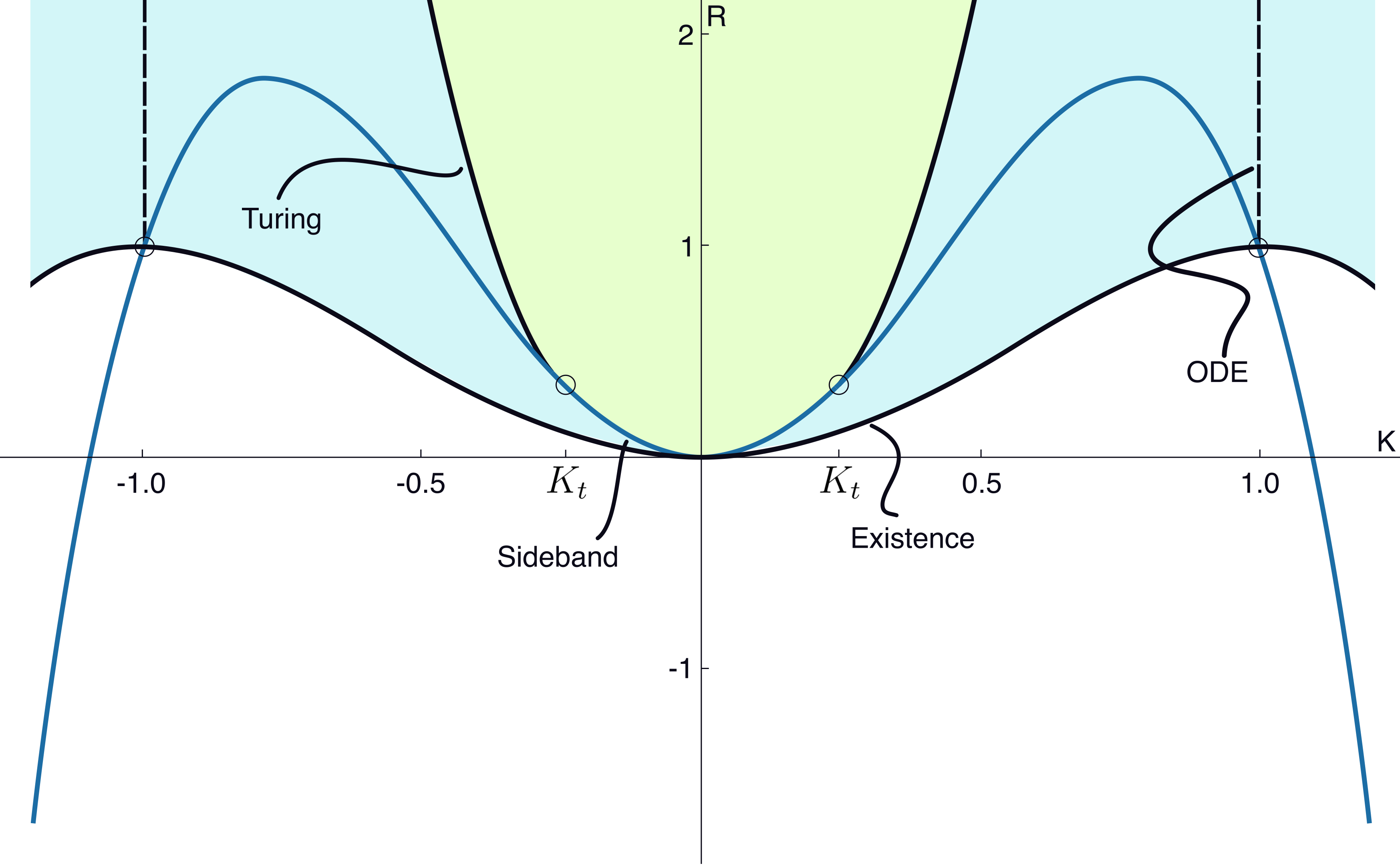}
    \caption{$d=12$}
    \label{fig:large d:cs}
\end{subfigure}
    \caption{Bifurcation diagrams for solutions $(A_p, B_p)$ in $(K,R)$-space for $\beta > 0$ for $d = 0.1$ (a) and $d = 12$ (b). In the white region below the existence curve, the solutions do not exist. In the light-blue regions, the solutions are unstable and in the light-green region, the solutions are stable.}
    \label{fig: Bifurcation Diagram Coupled System}
\end{figure}

The state $(A_p, B_p)$ destabilizes by a Turing bifurcation with critical wave number $k_c$ ($k_c \neq 0$) if $\lambda(k_c) = \lambda'(k_c) = 0$, where $\lambda(k)=\lambda_j(k)$ is one of the three branches determine by $Q(k,\lambda(k)) = 0$ \eqref{def:Qklambda}. Using
\[
0 \equiv \frac{d}{dk} [Q(k,\lambda(k))] = \frac{\partial Q}{\partial k} (k,\lambda(k)) + \frac{\partial Q}{\partial \lambda} (k,\lambda(k)) \lambda'(k), 
\]
we see that a (potential) Turing bifurcation takes place when
\[
Q(k_c,0) = \det(\mathcal{A}(k))|_{k=k_c} = 0 \quad \text{ and } \quad \frac{\partial Q}{\partial k} (k_c,0) =
\det(\mathcal{A}(k))'|_{k=k_c} =0
\]
\eqref{def:Qklambda}. Consequently, because $k_c \neq 0$ and by \eqref{detAk}, we have:
\[
2 Q(k_c,0) - k \frac{\partial Q}{\partial k}(k_c,0) = 4 k_c^2 \alpha \left(B_p + d(k_c^2 - 2K^2) \right)= 0.
\]
Therefore, we conclude that
\begin{equation}
\label{kcTuringAB}
k_c^2 = 2 K^2 - \frac{B_p}{d} = \frac{(2d + 1)K^2 - 1}{d}
\end{equation}
\eqref{eq: CS periodic solutions}; hence, since $d > 0$, a Turing bifurcation can only take place if 
\[
K^2 > \frac{1}{2d + 1} \in (0,1).
\]
Additionally, using
\[
0 = \det(\mathcal{A}(k_c)) = -\frac{1}{d}\left(1+(2d-1)K^2\right)^2 + 2(1-K^2)^2 + 2(R-1),
\]
we see that a Turing bifurcation occurs at $R = R_t(K)$:
\begin{equation}
\label{defRtK}
R_t(K) = 1 + \frac{1}{2d}\left(1 + \left(2d-1\right)K^2\right)^2 - \left(1-K^2\right)^2.
\end{equation}
Next, from \eqref{defRsK}, we deduce:
\[
R_t(K) - R_s(K) = \frac{1}{2d}\left(1 - 2(2d + 1)K^2 + (2d + 1)^2K^4\right) = \frac{1}{2d}\left(1 - (2d + 1)K^ 2\right)^2,
\]
which implies that the sideband curve $R=R_s(K)$ and Turing curve $R=R_t(K)$ are tangent at their (co-dimension 2) point(s) of intersection, 
\[
(K_{st},R_{st}) = \left(\frac{\pm 1}{\sqrt{2d+1}}, \frac{12d+1}{(2d+1)^2} \right).
\]
Note that the boundary of the region of stable periodic patterns consists of two non-empty symmetric (for $K \to -K$) Turing components. These components  are smooth connections between $(K_{st},R_{st})$ on the sideband instability curve $R_s(K)$ and $(\pm 1,R_{t}(\pm 1)) = (\pm 1,1+2d)$ on the ODE instability boundaries $K = \pm 1$. Therefore, these Turing regions shrink to an arbitrarily small region as $d \downarrow 0$ (since $(K_{st},R_{st}) \to (\pm 1,R_{t}(\pm 1)) = (\pm 1, 1)$) and become arbitrary large as $d \to \infty$ (since $(K_{st},R_{st}) \to (0,0)$) -- see \Cref{fig: Bifurcation Diagram Coupled System}. It is also worth noting that the region in which the $(A_p,B_p)$-pattern is destabilized by the classical Ginzburg-Landau sideband curve shrinks to $(K,R) = (0,0)$ in the latter limit (although one needs to be careful with a Ginzburg-Landau interpretation of the limit $d \gg 1$ since it is implicitly assumed in the derivation procedure that $d$ is $\mathcal{O}(1)$ with respect to the $\varepsilon$ introduced in the Ginzburg-Landau Ansatz (cf. \eqref{eq:scalar ginzburg-landau ansatz:ginzburg-landau derivation})).

To wrap up the stability analysis, we must address a few remaining issues. First, we need to determine the direction of the destabilization as the Turing curve is crossed. In other words, is $(A_p,B_p)$ stable for $R < R_t(K)$ or for $R > R_t(K)$? This question can be answered by approximating $Q(k_c,0)$ for $(A_p, B_p)=(A_p(K,R_t + \varepsilon), B_p(K,R_t + \varepsilon))$ with $0 < \varepsilon \ll 1$: 
\begin{equation*}
  Q(k_c, 0) = 2 \alpha (1 - (1+2d)K^2) \varepsilon/d +\mathcal{O}(\varepsilon^2);
\end{equation*}
thus, $Q(k_c,0) > 0$ if $R < R_t$ and $Q(k_c,0) < 0$ if $R > R_t$. Consequently, since, by definition, $Q(k_c, 0) = \lambda_1(k_c) \lambda_2(k_c) \lambda_3(k_c)$, we conclude that $(A_p,B_p)$ is stable for $R > R_t$. Second, we need to confirm that all curves $\lambda_j(k)$ satisfy Re$(\lambda_j(k)) < 0$ for all $k \neq \pm k_c$; luckily, this follows directly from inspection. 

The statement of Claim \ref{thm:close enough:t:cs} now follows (formally) by combining the derivations (and rescalings) of the AB-system in Sections \ref{sec:scalar} and \ref{sec:reaction diffusion system} with the above spectral stability analysis. For Claim \ref{thm:supercrit:t:cs} one also needs to incorporate the results on the relation between the `zoomed-in' Ginzburg-Landau equations and their AB-systems from Sections \ref{subsubsec:Reducing the AB-system:gs:se} and \ref{subsec:reducing the AB-system back to the ginzburg landau system:rd}.

\begin{remark}
\label{rem:alphainstab}
\rm    
The $\alpha$ coefficient of AB-system \eqref{eq:canonical form:cs} does not have a direct impact on the spectral stability of the $(A_p,B_p)$-patterns, since it does not appear in any of the boundary components of the stable pattern region: $K=\pm 1$ (ODE), $R=R_s(K)$ (sideband, \eqref{defRsK}), or $R_t(K)$ (Turing, \eqref{defRtK}). While somewhat unexpected, $\alpha$ may impact the (nonlinear) nature of the associated ODE, sideband, or Turing bifurcations, as already demonstrated in \Cref{fig:Intro-qp,fig:Intro-chaos} (only the values of $\alpha$ differ in the simulations shown there); see also \Cref{fig:ChaosAB}.
\end{remark}

%


\section{Simulations}\label{sec:simulations}

In this section, we present a series of simulations designed to achieve three goals. First, in \Cref{subsec:corroborating validity:s}, we corroborate the validity of the AB-system as a modulation equation by comparing the periodic solutions arising from the Turing bifurcation at $\mu = \mu^t$ in \eqref{eq:scalar example} with the approximations provided by its associated AB-system, given in \eqref{eq:coupled system:a model}. Specifically, we analyze the convergence rate of the $L^2$-norm difference between the small-amplitude periodic patterns associated with the Turing bifurcation in \eqref{eq:scalar example} and their corresponding approximations by the AB-system as $\delta \to 0$. Second, in \Cref{subsec:interesting behavior in the AB-system:s}, we investigate the evolution of solutions initialized near a specific periodic pattern at parameter values where the periodic solution is unstable. Finally, in \Cref{subsec:predicting behavior using the AB-system:s}, we utilize the insights from \Cref{subsec:interesting behavior in the AB-system:s} to predict behavior in PDEs exhibiting a Turing-fold bifurcation, thereby demonstrating its direct relevance for understanding and predicting the dynamics driven by the interplay between Turing and saddle-node bifurcations.

\begin{remark}
\label{rem:simulations}
\rm
{
In this section, we restrict ourselves to a preliminary numerical analysis of the AB-system and its use as a modulation equation. All simulations were carried out in Mathematica using the NDSolve function, with the MethodOfLines method with either periodic or homogeneous Neumann boundary conditions.
}
\end{remark}
\subsection{Corroborating the validity of the plane waves}
\label{subsec:corroborating validity:s}
The periodic solutions $A(\xi,\tau) = \bar{A} e^{iK\xi}, \; B(\xi,\tau) = \bar{B}$ of AB-system \eqref{eq:coupled system:a model} associated to model \eqref{eq:scalar example}, 
\begin{equation*}
    \begin{cases}
        A_{\tau} = 4 A_{\xi \xi} + A - 2 A B, \\
        B_{\tau} = B_{\xi \xi} + \frac{1}{4} - r - B^2 + 2(\eta - 1)|A|^2.
    \end{cases}
\end{equation*}
are given by
\begin{equation*}
    A_s(r) = 0, \quad B_s^\pm(r) = \pm \sqrt{\frac14 - r},
\end{equation*}
and
\begin{equation}\label{eq:per approx AB:cv:s}
A_p(K, r, \eta) = \sqrt{\frac{1}{2 (\eta - 1)}\left(B_p^2 + r - \frac14 \right)}, \quad B_p(K, r) = \frac12 \left(1 - 4 K^2\right).
\end{equation}

The solutions \eqref{eq:per approx AB:cv:s} correspond to periodic solutions $U_p(x; \delta, K,r, \eta)$ in \eqref{eq:scalar example}, emerging from the Turing bifurcation at $\mu = \mu^t = -1 + \frac{\delta^2}{4} + \mathcal{O}(\delta^3)$. In this section, our primary objective is to examine the convergence rate $\|U_p - U_{AB}\|_{L^2}$ as $\delta \to 0$, where
\begin{equation*}
  U_{AB}(x; \delta, K,r, \eta) = 1 + \delta B_p(K,r) + 2 \delta A_p(K,r, \eta) \cos\left(\left(1+\sqrt{\delta} K)x\right)\right).
\end{equation*}

We have the simulated model
\begin{equation}\label{eq:scalar example:vc:s}
    \partial_t U = \left(1 + \frac{\delta^2}{4} - r \delta^2\right) U + 2U^2 - U^3 + (1-\delta) \partial_x^2 U + 2 \partial_x^4 U + \partial_x^6 U + \eta (\partial_x^2 U)^2,
\end{equation}
(cf. \eqref{eq:scalar example}) on the spatial domain $\left[0, L\right]$, where $L \coloneqq n \cdot \frac{2 \pi}{1 + \sqrt{\delta} K}$. The simulations were initialized with
\begin{equation*}
    U_0(x; K, r, \eta) = U_{AB}(x;K, r, \eta),
\end{equation*}
and evolved until convergence to a periodic solution $U_p$ was attained. Subsequently, we computed the normalized $L^2$-norm of the difference between $U_p$ and $U_{AB}$, and postulated that, as $\delta \to 0$,
\begin{equation*}
    \|U_p - U_{AB}\|_{L^2} \sim C(K,r, \eta) \delta^{e(K,r,\eta)},
\end{equation*}
for some positive constants $C(K,r, \eta), e(K,r,\eta) \in \mathbb{R}^+$. To estimate the exponent $e(K,r,\eta)$ $\delta$, we varied $\delta$ and used the approximation 
\begin{equation}
    e(K,r,\eta) \approx \frac{\log(\|U_p(x;\delta_1, K, r, \eta) - U_{AB}(x;\delta_1, K, r, \eta)\|_{L^2} - \log(\|U_p(x;\delta_2, K, r, \eta) - U_{AB}(x;\delta_2, K, r, \eta)\|_{L^2}}{\log(\delta_1) - \log(\delta_2)}.
\end{equation}
The results can be found in \Cref{tab:combined convergence rates:s}. 

\begin{table}[ht]
\centering
\begin{subtable}{0.33\textwidth}
\centering
\begin{tabular}{|c|c|c|}
\hline
$\delta$ & $\|U - U_{AB}\|_{L^2}$ & $e$ \\
\hline\hline
0.02 & 0.002 & 2.071 \\
0.04 & 0.007 & 1.987 \\
0.06 & 0.016 & 1.954 \\
0.08 & 0.029 & 1.934 \\
0.10 & 0.044 & 1.920 \\
0.12 & 0.063 & 1.910 \\
0.14 & 0.084 & 1.904 \\
0.16 & 0.108 & 1.902 \\
0.18 & 0.135 & 1.905 \\
0.20 & 0.165 & 1.912 \\
\hline
\end{tabular}
\caption{$K = 0$}
\label{tab:K=0 simulation convergence rate:cv:s}
\end{subtable}%
\begin{subtable}{0.33\textwidth}
\centering
\begin{tabular}{|c|c|c|}
\hline
$\delta$ & $\|U - U_{AB}\|_{L^2}$ & $e$ \\
\hline\hline
0.02 & 0.002 & 1.680 \\
0.04 & 0.008 & 1.761 \\
0.06 & 0.017 & 1.806 \\
0.08 & 0.029 & 1.837 \\
0.10 & 0.044 & 1.852 \\
0.12 & 0.061 & 1.856 \\
0.14 & 0.082 & 1.861 \\
0.16 & 0.105 & 1.863 \\
0.18 & 0.130 & 1.868 \\
0.20 & 0.159 & 1.874 \\
\hline
\end{tabular}
\caption{$K = 0.1$}
\label{tab:K=0.1simulation convergence rate:cv:s}
\end{subtable}%
\begin{subtable}{0.33\textwidth}
\centering
\begin{tabular}{|c|c|c|}
\hline
$\delta$ & $\|U - U_{AB}\|_{L^2}$ & $e$ \\
\hline\hline
0.02 & 0.006 & 1.478 \\
0.04 & 0.016 & 1.451 \\
0.06 & 0.029 & 1.476 \\
0.08 & 0.044 & 1.491 \\
0.10 & 0.061 & 1.503 \\
0.12 & 0.081 & 1.516 \\
0.14 & 0.103 & 1.533 \\
0.16 & 0.126 & 1.550 \\
0.18 & 0.152 & 1.564 \\
0.20 & 0.180 & 1.593 \\
\hline
\end{tabular}
\caption{$K = 0.3$}
\label{tab:K=0.3 simulation convergence rate:cv:s}
\end{subtable}
\caption{Convergence rates from the numerical experiments described in \Cref{subsec:corroborating validity:s}, for various values of $K$, with $t_{\text{max}} = 10000$, $r = 4$, $\eta = 2$, and $L = 6\pi$.}
\label{tab:combined convergence rates:s}
\end{table}

It can be verified that the next-to-leading-order approximation of the AB-system is given by:
\begin{equation}\label{eq:refined AB system:cv:s}
\begin{cases}
A_{\tau} = 4 A_{\xi \xi} + A - 2 A B - 12 i \sqrt{\delta} A_{\xi\xi\xi} + \mathcal{O}(\delta), \\
B_{\tau} = B_{\xi \xi} + \frac{1}{4} - r - B^2 + 2(\eta - 1)|A|^2 + 4 \eta i \sqrt{\delta} \left( A A^*_\xi - A^* A_\xi \right) + \mathcal{O}(\delta)
\end{cases}
\end{equation}
(see also Remark \ref{rem:degenerate}). This refined formulation demonstrates that the numerical outcomes align closely with theoretical expectations. Specifically, when $ K = 0$ we have $A_{\xi\xi\xi} = 0$ and $A_\xi = 0$, which directly corresponds to having $X_{02}(\xi,\tau) \equiv X_{12}(\xi,\tau) \equiv 0$ in the expansion/Ansatz \eqref{eq: Coupled system Ansatz}. Thus, in this case, the convergence rate is anticipated to approach $e = 2$, as confirmed by Table \ref{tab:combined convergence rates:s}(a). By contrast, for $ K \neq 0 $, the presence of the $\sqrt{\delta} $-order terms -- or equivalently, the fact that $X_{02}, X_{12}$ in \eqref{eq: Coupled system Ansatz} do not vanish -- implies a convergence rate of $ e = \frac{3}{2} $. These predictions are also fully consistent with the results obtained from our numerical experiments (see Table \ref{tab:combined convergence rates:s}(b), (c)). 

Thus, the simulations not only corroborate the validity of the approximation by the AB-system -- which can only be expected to hold for finite time -- but also the persistence of stable patterns of the AB-system into the underlying model for all time (cf. Remark \ref{rem:validity}). 

\subsection{The dynamics of the AB-system}
\label{subsec:interesting behavior in the AB-system:s}

\begin{figure}[tp]
    \centering
    
    \begin{subfigure}[b]{0.3\textwidth}
        \centering
        \includegraphics[width=\textwidth]{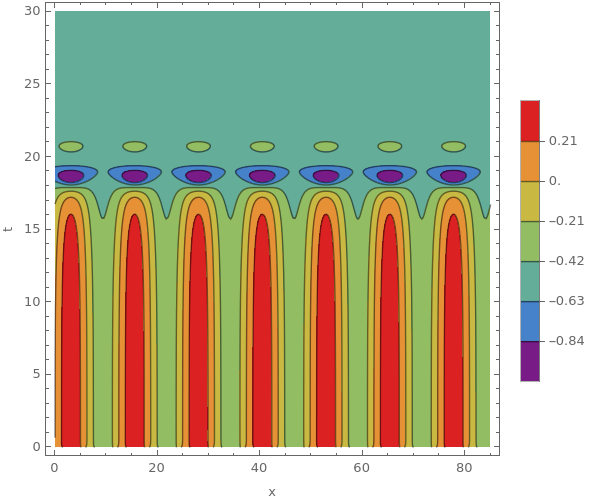}
        \caption{$\text{Re}\,(A(\xi,\tau))$; $(K,R) = (1.01, 2)$}
        \label{fig:K101U}
    \end{subfigure}
    \hfill
    \begin{subfigure}[b]{0.3\textwidth}
        \centering
        \includegraphics[width=\textwidth]{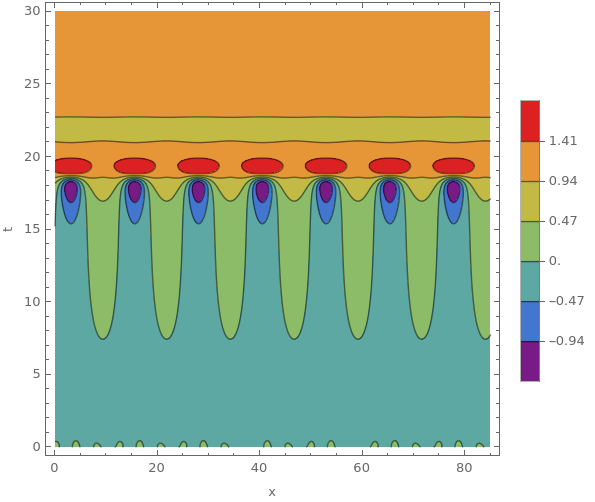}
        \caption{$B(\xi,\tau)$; $(K,R) = (1.01, 2)$}
        \label{fig:K101W}
    \end{subfigure}
    \hfill
    \begin{subfigure}[b]{0.315\textwidth}
        \centering
        \includegraphics[width=\textwidth]{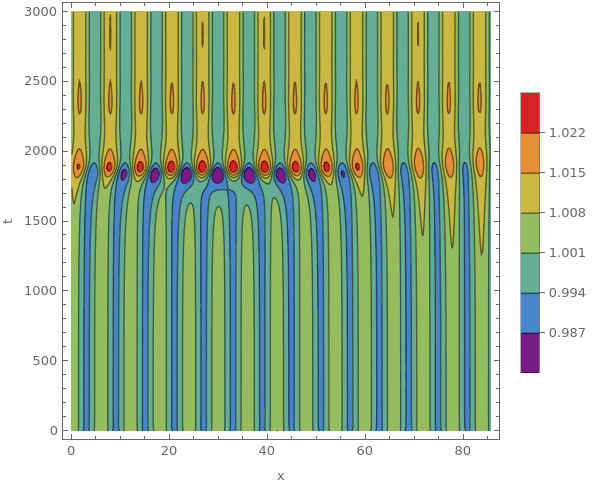}
        \caption{$U_{AB}(x,t)$; $(K,R) = (1.01, 2)$}
        \label{fig:K101AB}
    \end{subfigure}
    
    \vspace{1em}
    
    \begin{subfigure}[b]{0.3\textwidth}
        \centering
        \includegraphics[width=\textwidth]{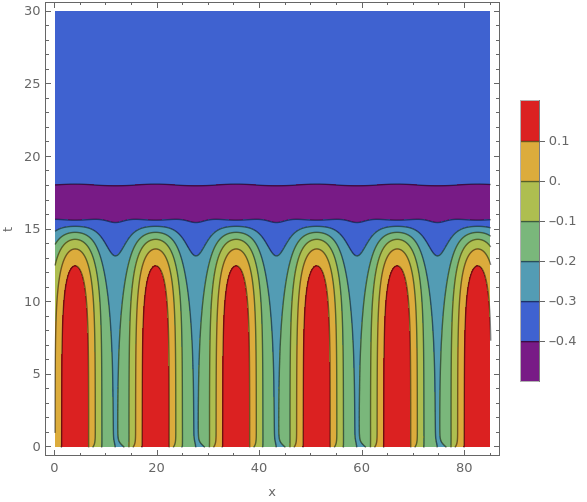}
        \caption{$\text{Re}\,(A(\xi,\tau))$; $(K,R) = (0.8, 1.2)$}
        \label{fig:K08U}
    \end{subfigure}
    \hfill
    \begin{subfigure}[b]{0.3\textwidth}
        \centering
        \includegraphics[width=\textwidth]{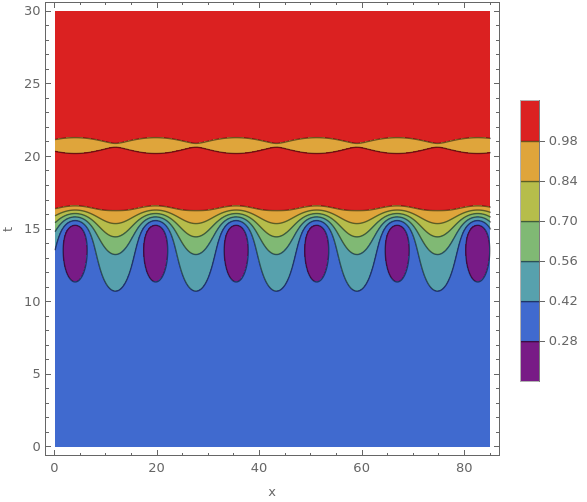}
        \caption{$B(\xi,\tau)$; $(K,R) = (0.8, 1.2)$}
        \label{fig:K08W}
    \end{subfigure}
    \hfill
    \begin{subfigure}[b]{0.315\textwidth}
        \centering
        \includegraphics[width=\textwidth]{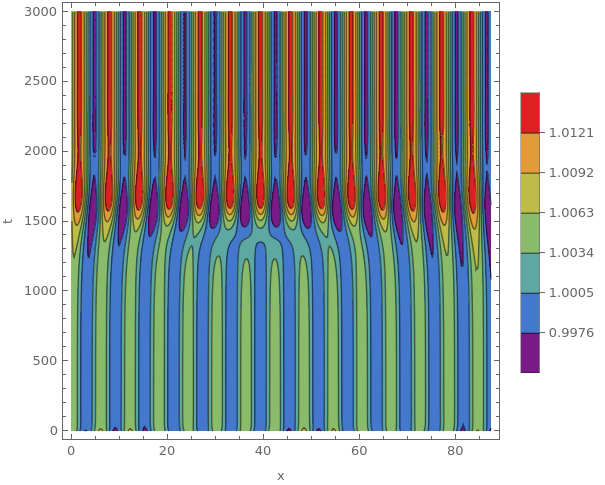}
        \caption{$U_{AB}(x,t)$; $(K,R) = (0.8, 1.2)$}
        \label{fig:K08AB}
    \end{subfigure}
    
    \vspace{1em}
    
    \begin{subfigure}[b]{0.3\textwidth}
        \centering
        \includegraphics[width=\textwidth]{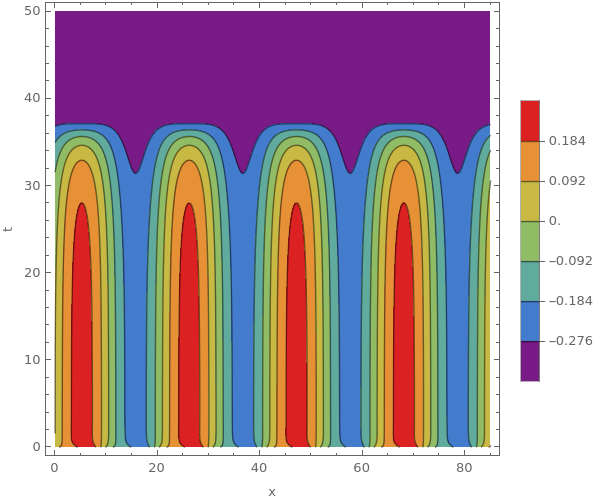}
        \caption{$\text{Re}\,(A(\xi,\tau))$; $(K,R) = (0.6, 1)$}
        \label{fig:K06U}
    \end{subfigure}
    \hfill
    \begin{subfigure}[b]{0.3\textwidth}
        \centering
        \includegraphics[width=\textwidth]{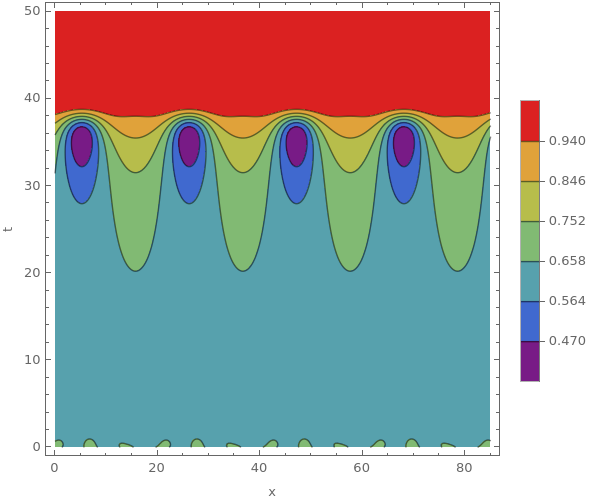}
        \caption{$B(\xi,\tau)$; $(K,R) = (0.6, 1)$}
        \label{fig:K06W}
    \end{subfigure}
    \hfill
    \begin{subfigure}[b]{0.315\textwidth}
        \centering
        \includegraphics[width=\textwidth]{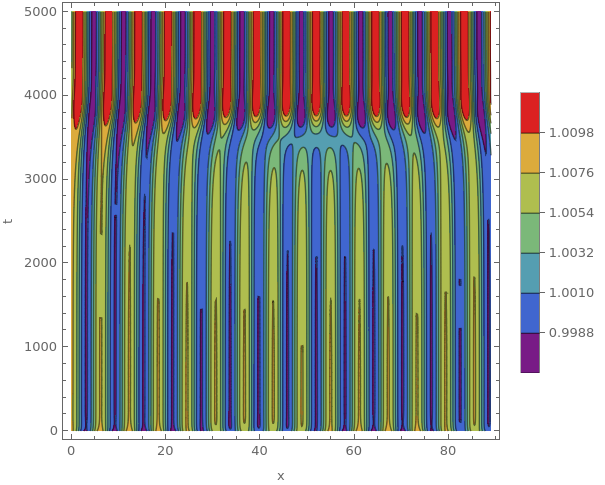}
        \caption{$U_{AB}(x,t)$; $(K,R) = (0.6, 1)$}
        \label{fig:K06AB}
    \end{subfigure}
       \caption{Simulations of the evolution of $\text{Re}\,(A(x,t))$ (first column), $B(x,t)$ (second column) as solutions to \eqref{eq:canonical form:cs}, and $U_{AB}(x,t)$ (third column) as defined in \eqref{eq:AB An}, are shown for parameter values $\alpha = d = \frac12$, $\beta = 8$ and $\delta = 0.01$.  Each row corresponds to a fixed pair $(K, R)$, indicating the initial condition: an unstable plane wave with $(K, R)$. The simulations are carried out on the domain $\xi \in (0, L)$ with periodic boundary conditions, where $L = \tfrac{600\pi}{K}$. For clarity, only subintervals of equal length are displayed. The plots of $U_{AB}(x,t)$ will later be compared to simulations of an underlying system in \Cref{subsec:predicting behavior using the AB-system:s}.}
    \label{fig:ABsimulations ODE and sideband instability}
\end{figure}

Turing bifurcations give rise to families of spatially periodic solutions whose stability is determined by the (bifurcation) parameters in the underlying system. While the Ginzburg–Landau equation and AB-system capture only the onset of these patterns, they typically persist well beyond the regimes where these reduced models are valid. The full stability landscape of such periodic solutions in parameter–wavenumber space (in our case, $\mu$–$k$ space) is known as the {\it Busse balloon} \cite{busse1978non,doelman2012hopf}. In this section, we consider the natural and relevant question: {\it What happens to solutions near a specific periodic pattern when that pattern destabilizes as the bifurcation parameter moves out of the Busse balloon and into an unstable regime?}

We present simulations of the dynamics exhibited by the AB-system \eqref{eq:canonical form:cs} for $\beta > 0$ that start from initial conditions near a periodic state $(A_p(K,R), B_p(K,R))$ \eqref{eq: CS periodic solutions} that is marginally unstable -- i.e., for parameter combinations just outside the Busse balloon but close to its boundary. As summarized in \Cref{sec:coupled system}, periodic solutions $(A_p, B_p)$ may exist and be unstable for $\beta > 0$ due to one of three distinct mechanisms: sideband instability, ODE instability ($K > 1$), or Turing instability -- see \Cref{fig: Bifurcation Diagram Coupled System}.

In \Cref{fig:ABsimulations ODE and sideband instability}, we show simulations of \eqref{eq:canonical form:cs} initialized near periodic solutions that are (marginally) unstable due to either ODE or sideband instability. These simulations reveal that such instabilities typically lead to pattern reselection: once the initial periodic solution becomes unstable, the system evolves toward a stable periodic pattern. In other words, the system evolves back into the Busse balloon. This behavior is well-known; see \cite{kashyap2025laminar,rietkerk2021evasion,siteur2014beyond} and the references therein for numerically-obtained away-from-onset examples and \cite{asch2025slow,kramer1985eckhaus} for analytical studies of this reselection mechanism in the context of the Ginzburg-Landau equation. 

\begin{figure}[t]
    \centering
    \begin{subfigure}[t]{0.28\textwidth}
        \centering
        \includegraphics[width=\textwidth]{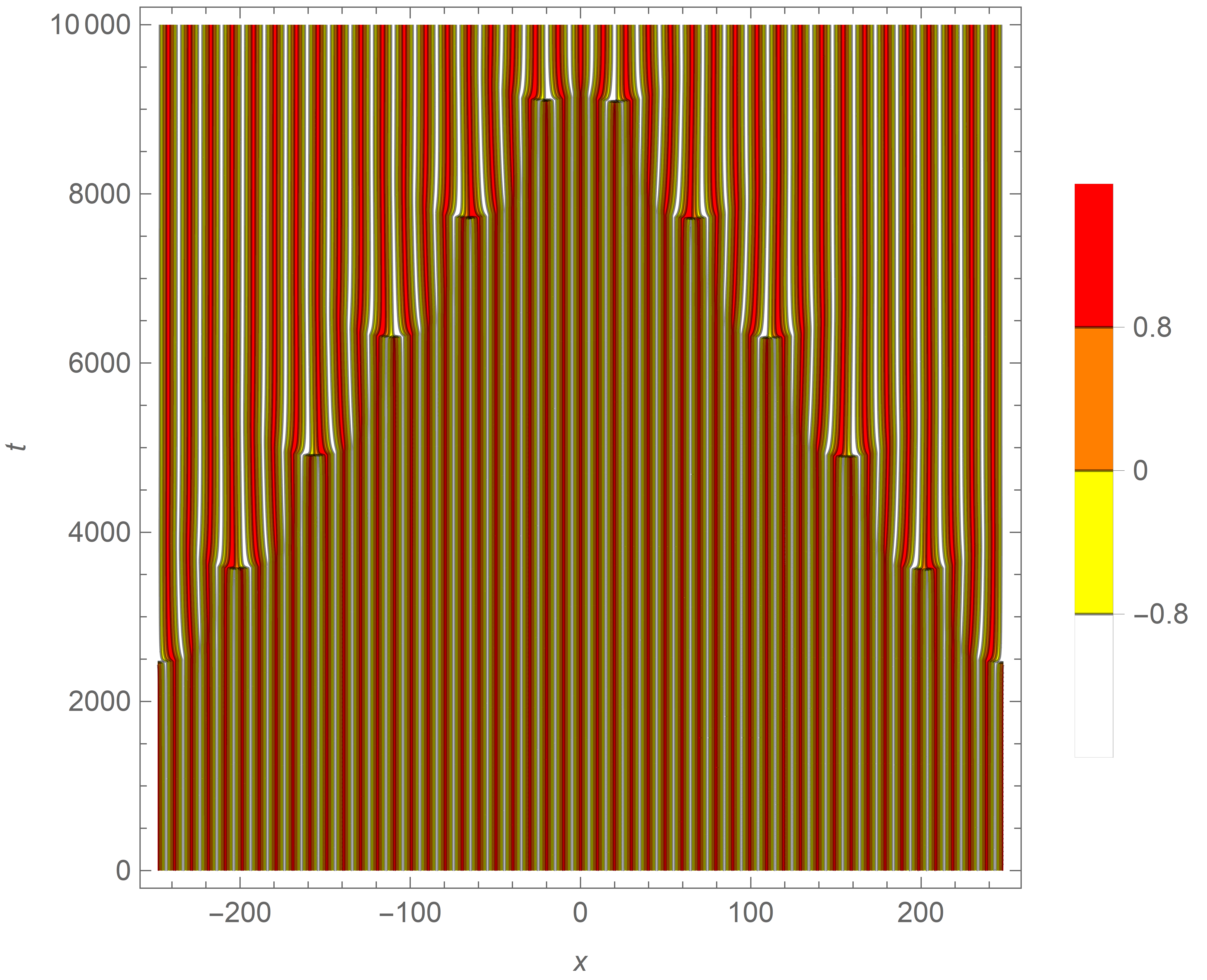}
        \caption{Contour plot of Re$\,(A(\xi,\tau))$ for $\tau \in (0, 10000)$.}
        \label{fig:Extra-Contour}
    \end{subfigure}
    \hfill
    \begin{subfigure}[t]{0.35\textwidth}
        \centering
        \includegraphics[width=\textwidth]{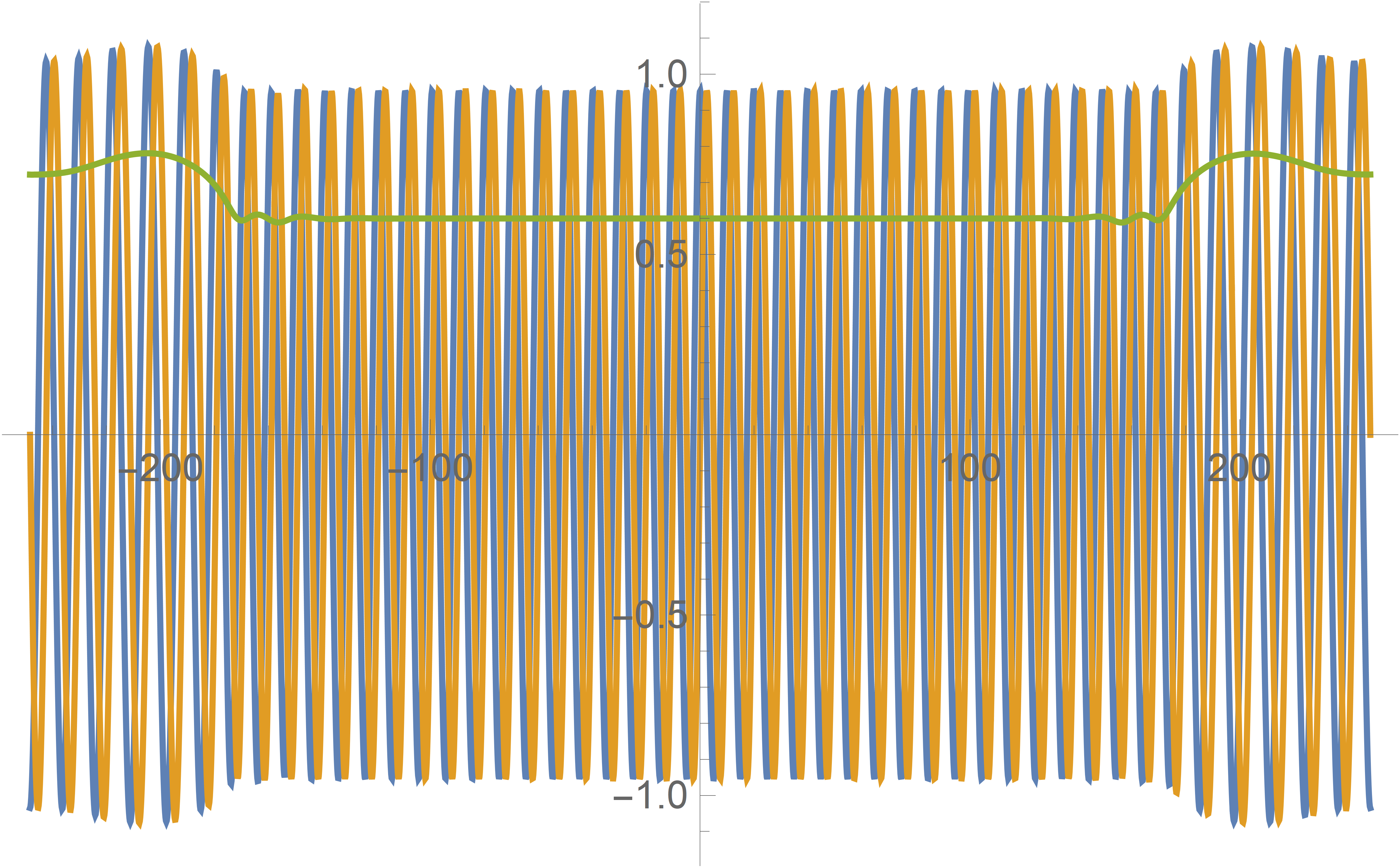}
        \caption{Re$\,(A(\xi,\tau))$ (blue), Im$\,(A(\xi,\tau))$ (orange) and Re$B(\xi,\tau)$ (green) at $\tau = 4000$.}
        \label{fig:Extra-4000}
    \end{subfigure}
    \hfill
    \begin{subfigure}[t]{0.35\textwidth}
        \centering
        \includegraphics[width=\textwidth]{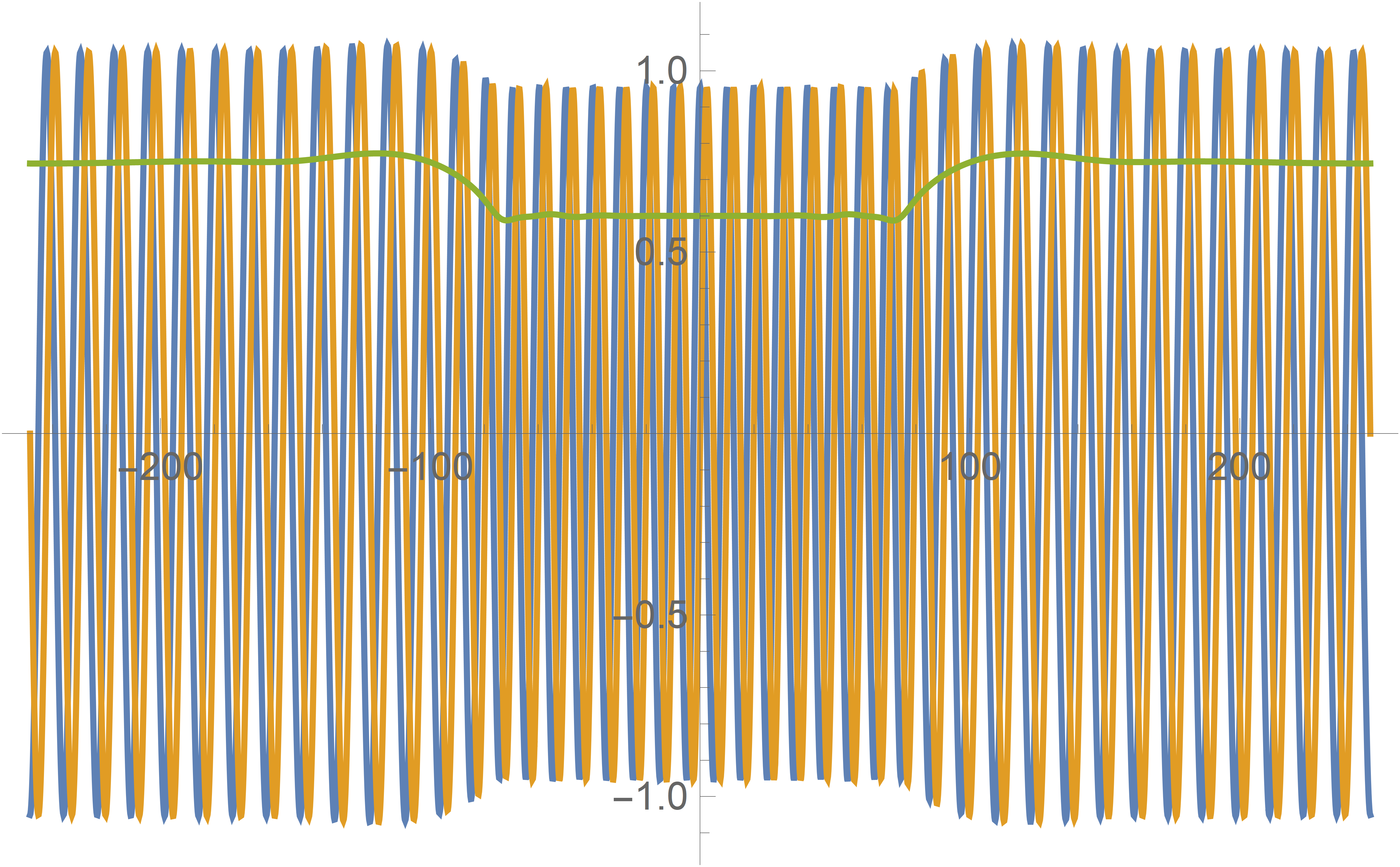}
        \caption{Re$\,(A(\xi,\tau))$ (blue), Im$\,(A(\xi,\tau))$ (orange) and Re$B(\xi,\tau)$ (green) at $\tau = 7000$.}
        \label{fig:Extra-7000}
    \end{subfigure}
    \caption{Simulation of \eqref{eq:canonical form:cs} with $\alpha = 1.4$, $d=1/3$, $\beta =1$, initialized with a plane wave solution \eqref{eq: CS periodic solutions} where $K = \sqrt{2/5}$ and $R = R_s(\sqrt{2/5}) - 0.05 = 8/5-0.005$ \eqref{defRsK}; this value of $R$ ensures that the pattern is weakly unstable with respect to the sideband instability. The simulation is done on the interval $\xi \in (-50 \pi \sqrt{2/5}, 50 \pi \sqrt{2/5})$) with homogeneous Neumann boundary conditions.}
    \label{fig:Extra}
\end{figure}

Note that each of the simulations shown in \Cref{fig:ABsimulations ODE and sideband instability} has the {\it Stoke wave} as endstate: since $A(\xi,\tau)$ has become constant, the final attracting plane wave \eqref{eq:periodic solutions:cs} has $K=0$, the wave number of the periodic pattern in the underlying equation is exactly $k^*$ \eqref{eq:Ansatz:abd:rd-Intro}, the {\it linearly most unstable wave}. This is certainly not necessary (as is also known from the above-mentioned literature). In \Cref{fig:Extra} we show a simulation of the evolution of a plane wave that is marginally unstable with respect to the sideband instability back into the Busse balloon, i.e., towards a stable plane wave (with $K\neq0$): the transformation seems to occur locally near invading destabilization fronts that are initiated at the boundary of the domain.

Nevertheless, the simulations shown so far indicate that the dynamics exhibited by the AB-system may be very similar to that of the real Ginzburg-Landau equation: it seems to be restricted to evolution within the family of plane waves. However, this is not the case. The Turing instability does not occur in the real Ginzburg-Landau equation; therefore, it is natural to look for non-real-Ginzburg-Landau-type behavior near (and especially just beyond) the Turing instability boundary (cf. \Cref{fig: Bifurcation Diagram Coupled System}). Our simulations, which are certainly not exhaustive ( \Cref{rem:sims}), indicate that the nature of the Turing bifurcation may change between super- and sub-critical as parameters in \eqref{eq:canonical form:cs} vary and that this plays a crucial role in the kind of dynamics exhibited by AB-system \eqref{eq:canonical form:cs}. (It should be noted that the Turing-bifurcation terminology used here is overly simplified: periodic patterns have a spectrum attached to $k=0$ (by translational invariance), this is not the case for classical Turing bifurcations appearing from homogeneous states.) 

\begin{figure}[tp]
    \centering
    \begin{subfigure}[b]{0.45\textwidth}
        \centering
        \includegraphics[width=\textwidth]{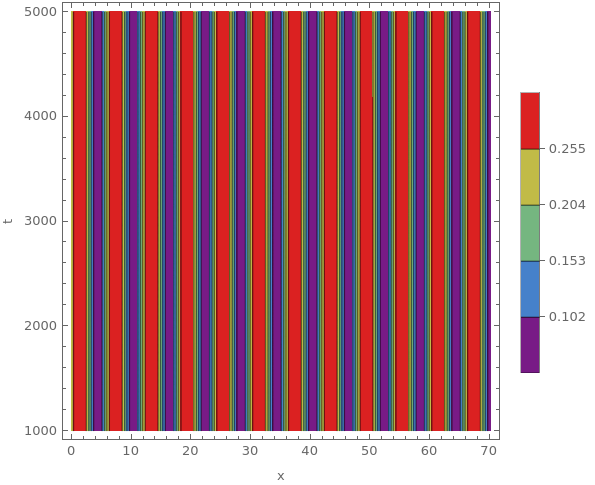}
        \caption{$\alpha = 0.8$}
        \label{eq:chaos subfiga}
        \end{subfigure}
    \hfill
        \begin{subfigure}[b]{0.45\textwidth}
        \centering
        \includegraphics[width=\textwidth]{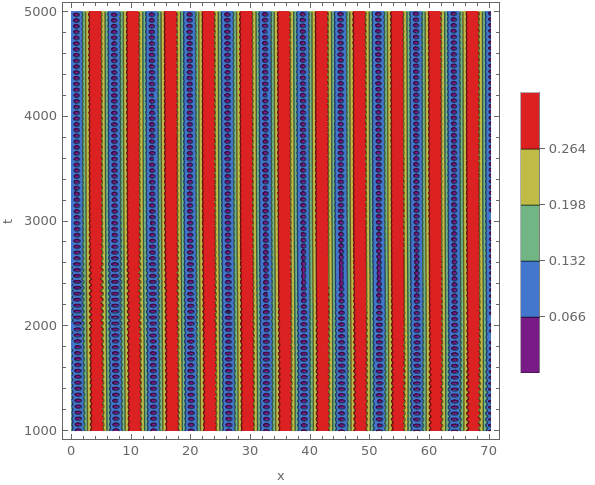}
        \caption{$\alpha = 0.7745$}
        \label{eq:chaos subfigb}
        \end{subfigure}
    
    \vspace{0.5cm} 
        \begin{subfigure}[b]{0.45\textwidth}
        \centering
        \includegraphics[width=\textwidth]{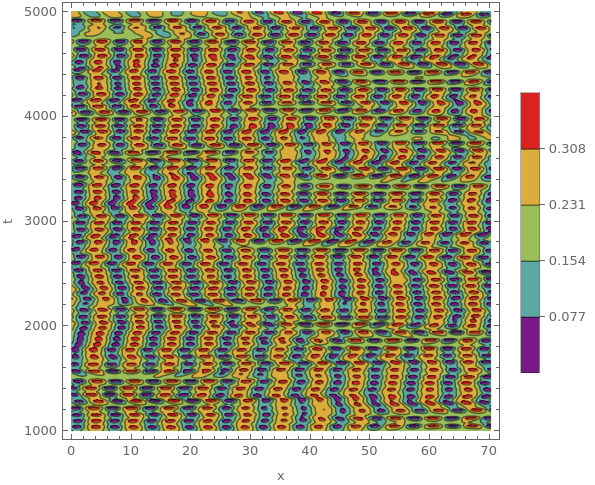}
        \caption{$\alpha = 0.75715$}
        \label{eq:chaos subfigc}
        \end{subfigure}
    \hfill
        \begin{subfigure}[b]{0.45\textwidth}
        \centering
        \includegraphics[width=\textwidth]{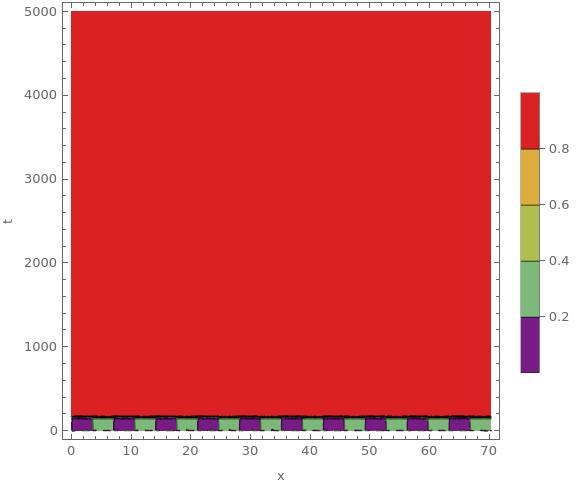}
        \caption{$\alpha = 0.7$}
        \label{eq:chaos subfigd}
        \end{subfigure}    
    \caption{Simulations of \eqref{eq:canonical form:cs} showcasing the evolution of $B(\xi,\tau)$ for decreasing values of $\alpha$, with $\beta = 1$ and  $d = 1/3$, starting from a plane wave which is marginally unstable, with respect to the Turing instability: $K = \sqrt{4/5}$ and $R = R_t\left(\sqrt{4/5}\right) - 0.01$ (see \eqref{defRtK}). The simulations were conducted on the interval $(0,L)$, with $L = 200 \pi / \sqrt{4/5}$, under periodic boundary conditions, for $\tau \in (0, 5000)$. Only selected subdomains of the full simulations are shown.}
    \label{fig:ChaosAB}
\end{figure}

We have observed that pattern reselection typically occurs near the part of the (local) Busse balloon boundary where the Turing bifurcation is subcritical. However, as previously observed in \Cref{fig:Intro-qp} and \ref{fig:Intro-chaos}, more intricate dynamics arise when the Turing bifurcation is supercritical. In \Cref{fig:ChaosAB}, we consider the same setting as in \Cref{fig:Intro-qp,fig:Intro-chaos}, i.e., we take $K = \sqrt{4/5}$, $R = R_t(\sqrt{4/5}) - 0.01$ \eqref{defRtK}, $\beta = 1$, $d = 1/3$, and decrease $\alpha$ in four steps from $0.8$ to $0.7$.
\Cref{eq:chaos subfigb} shows that the stationary quasi-periodic pattern of \Cref{fig:Intro-qp,eq:chaos subfiga}, which both have $\alpha =0.8$, undergoes a temporal Hopf bifurcation as $\alpha$ decreases (Remark \ref{rem:sims}). As $\alpha$ decreases further, the dynamics of AB-system \eqref{eq:canonical form:cs} becomes irregular, leading to the complex dynamics shown in \Cref{fig:Intro-chaos,eq:chaos subfigc} at $\alpha =0.75715$. This spatio-temporal chaos-like behavior only occurs in a relatively thin region in parameter spaces. As $\alpha$ decreases further, the nature of the Turing bifurcation changes from super- to subcritical and we observe the expected reselection process by which the system `jumps back' into the Busse balloon (\Cref{eq:chaos subfigd}).    

\begin{remark}
\label{rem:sims}
\rm 
We (re)emphasize that our simulations are by no means exhaustive. They provide illustrative examples of the dynamics that may arise in the AB-system; however, we cannot (and do not) claim that these capture the full range of possible behaviors, nor do we claim that the `spatio-temporal chaos-like' behavior shown in \Cref{fig:Intro-chaos,eq:chaos subfigc} and upcoming \Cref{fig:chaos in Arjen system} indeed exhibit the (mathematical) characteristics of chaotic dynamics. Even the precise location (in parameter space) of the Hopf bifurcation that initiated the (time-)periodic behavior shown in \Cref{eq:chaos subfigb} is modified in a subtle manner by the accuracy and the integration time of the simulation. A more comprehensive numerical and analytical investigation is needed to better understand the dynamics of AB-systems. Naturally, this will be the topic of future investigations. 
\end{remark}

\subsection{The persistence of the dynamics of AB-system into the underlying system}\label{subsec:predicting behavior using the AB-system:s}
The primary goal of the perturbation analysis presented in this paper is to gain general insights into local, small-amplitude spatial dynamics near a Turing-fold bifurcation. Our analysis allowed us to derive a system of coupled modulation equations that, assuming it is a valid approximation, captures this behavior near the co-dimension 2 bifurcation point. As such, this system should allow us to predict the dynamics of the underlying systems by analyzing the AB-system with appropriately translated coefficients. In \Cref{subsubsec:pattern reselection:pred}, we demonstrate that pattern reselection occurs where the AB-modulation equation predicts it. Furthermore, in \Cref{subsubsec:chaos:pbab:s}, we show that the onset of chaos-like dynamics exhibited by the AB-system also aligns with the parameter regime identified in \Cref{subsec:interesting behavior in the AB-system:s}.

\subsubsection{Pattern (re)selection}\label{subsubsec:pattern reselection:pred}

In \Cref{fig:ABsimulations ODE and sideband instability}, we demonstrate pattern (re)selection resulting from ODE or sideband instabilities in the canonical AB-system \eqref{eq:canonical form:cs}, for specific choices of $\alpha$, $\beta$, and $d$. In \Cref{subsubsec:pattern reselection:pred}, we present simulations of the scalar example \eqref{eq:scalar example}, using $\mu = 1 + \frac{\delta^2}{4}(1 - R)$, $\nu = 1 - \delta$, suitable values of $\delta$, $R$, $\eta$, $K$, and a domain size chosen to ensure that the associated AB-system matches the one used in \Cref{fig:ABsimulations ODE and sideband instability}.

By introducing the rescaled variable $\tilde{\xi}=\frac{1}{2}\sqrt{\delta} x$ and using the scaled Ansatz
\begin{equation}\label{eq:AB An}
    U_{AB}(x, t) = 1  + \frac{\delta}{2} \tilde{B}( \tilde{\xi}, \tau) + \delta e^{i x} A(\tilde{\xi}, \tau) + \text{c.c.},
\end{equation}
we bring the AB-system \eqref{eq:coupled system:a model}, associated with \eqref{eq:scalar example}, into its canonical form \eqref{eq:canonical form:cs}:
\begin{equation*}
    \begin{cases}
        \phantom{2}A_\tau = A_{\tilde{\xi\xi}} + A - A \tilde{B}, \\
        2 \tilde{B}_\tau = \frac{1}{2} A_{\tilde{\xi}\tilde{\xi}} + 1 - R - \tilde{B}^2 + 8 (\eta - 1) |A|^2.
    \end{cases}
\end{equation*}
Using this form, we conclude that the parameter choices $\alpha = d = \frac12$ and $\beta = 8$ in the simulations of \eqref{eq:canonical form:cs} shown in \Cref{fig:ABsimulations ODE and sideband instability} indeed correspond to the AB-system associated to \eqref{eq:scalar example} if we set $\eta = 2$. 

\begin{figure}[tp]
\centering
    \begin{subfigure}[t]{0.3\textwidth}
        \centering
        \includegraphics[width=\textwidth]{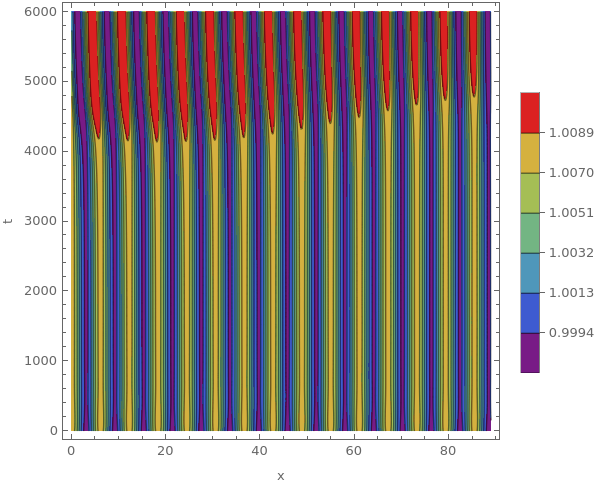}
        \caption{$(K,R)=(0.6,1)$}
        \label{fig:K06real}
    \end{subfigure}
\hfill
    \begin{subfigure}[t]{0.3\textwidth}
        \centering
        \includegraphics[width=\textwidth]{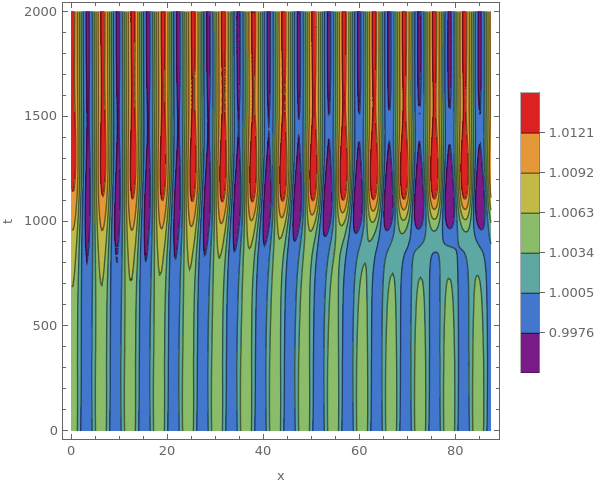}
        \caption{$(K,R)=(0.8,1.2)$}
        \label{fig:K08real}
    \end{subfigure}
\hfill
    \begin{subfigure}[t]{0.3\textwidth}
        \centering
        \includegraphics[width=\textwidth]{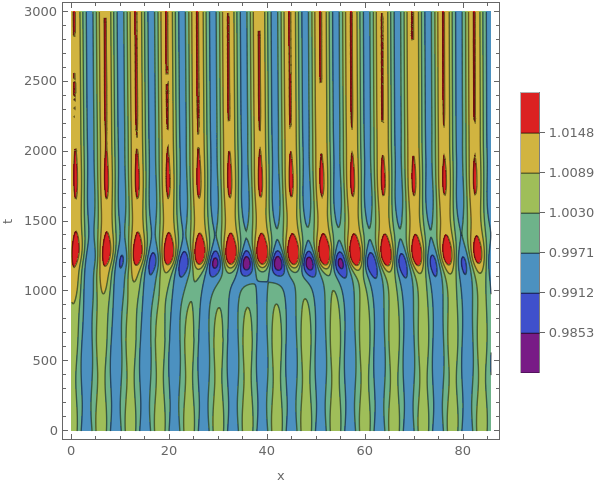}
      \caption{$(K,R)=(1.01,2)$}
        \label{fig:K101real}
    \end{subfigure}
\centering
    \begin{subfigure}[t]{0.3\textwidth}
        \centering
        \includegraphics[width=\textwidth]{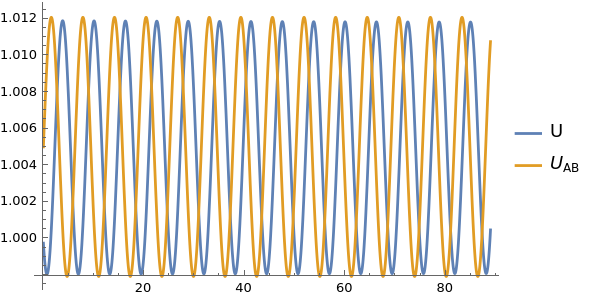}
        \caption{$(K,R)=(0.6,1)$}
        \label{fig:K06t5000}
    \end{subfigure}
\hfill
    \begin{subfigure}[t]{0.3\textwidth}
        \centering
        \includegraphics[width=\textwidth]{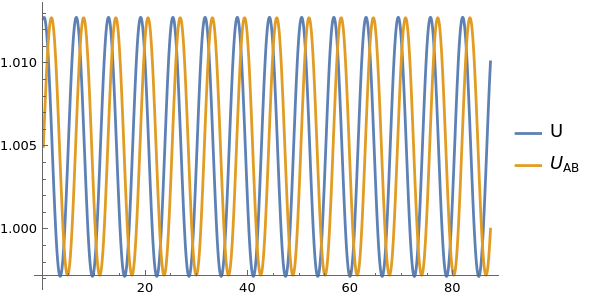}
        \caption{$(K,R)=(0.8,1.2)$}
        \label{fig:K08t5000}
    \end{subfigure}
\hfill
    \begin{subfigure}[t]{0.3\textwidth}
        \centering
        \includegraphics[width=\textwidth]{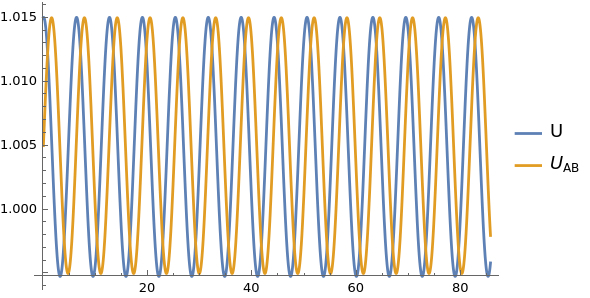}
      \caption{$(K,R)=(1.01,2)$}
        \label{fig:K101t5000}
    \end{subfigure}
\caption{Simulations of the evolution of solutions $U(x,t)$ to \eqref{eq:scalar example} are shown for parameter values corresponding to those used in \Cref{fig:ABsimulations ODE and sideband instability}; specifically, $\eta = 2$, $R = 1$, $\delta = 0.01$, and domain size $L = \frac{600 \pi}{1 + \frac{1}{2}\sqrt{\delta} K}$. The first row displays contour plots of $U(x,t)$ for $(K,R)=(0.6,1)$, $(K,R)=(0.8,1.2)$, and $(K,R)=(1.01,2)$, initialized with (approximately) the same initial condition as $U_{AB}(x,t)$ in \Cref{fig:ABsimulations ODE and sideband instability}. These plots should thus be compared to the contour plots of $U_{AB}(x,t)$ in the third column of \Cref{fig:ABsimulations ODE and sideband instability}. In the second row, the endstates $U(x,5000)$ (in blue) and $U_{AB}(x,5000)$ (in orange) of the simulations are compared. Note that both the amplitude and the period of $U_{AB}(x,5000)$ approximate those of $U(x,5000)$ very well; however, the phases differ slightly, which is due to the periodic boundary conditions.} 
\label{fig:reselectionK = 0.60.81.01}
\end{figure}

\Cref{fig:reselectionK = 0.60.81.01} shows the outcome of simulations of model \eqref{eq:scalar example} corresponding to those (of the AB-system) in \Cref{fig:ABsimulations ODE and sideband instability}. Each of the three cases exhibits a pattern (re)selection process that is similar to their counterparts in \Cref{fig:ABsimulations ODE and sideband instability}. In fact, the selected patterns in the scalar PDE model \eqref{eq:scalar example} accurately match those selected by the AB-system -- see the second row of \Cref{fig:reselectionK = 0.60.81.01}.

\subsubsection{Complex dynamics}
\label{subsubsec:chaos:pbab:s}

\begin{figure}[t]
    \centering
    \begin{subfigure}[b]{0.48\textwidth}
        \centering
        \includegraphics[width=\textwidth]{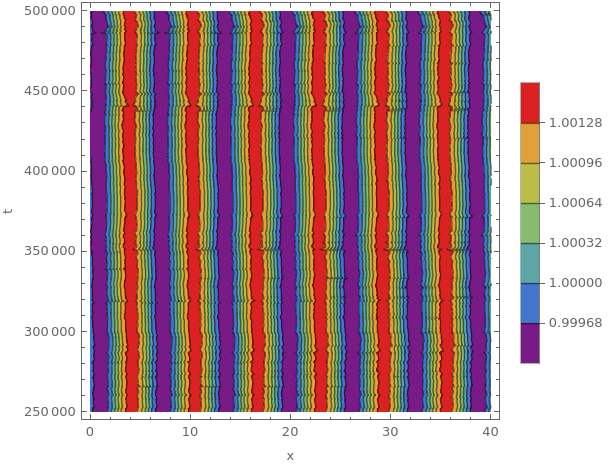}
        \caption{contour plot of $U(x,t)$}
        \label{fig:Ucontour plot:chaosArjenSystem}
    \end{subfigure}
    \hfill
    \begin{subfigure}[b]{0.48\textwidth}
        \centering
        \includegraphics[width=\textwidth]{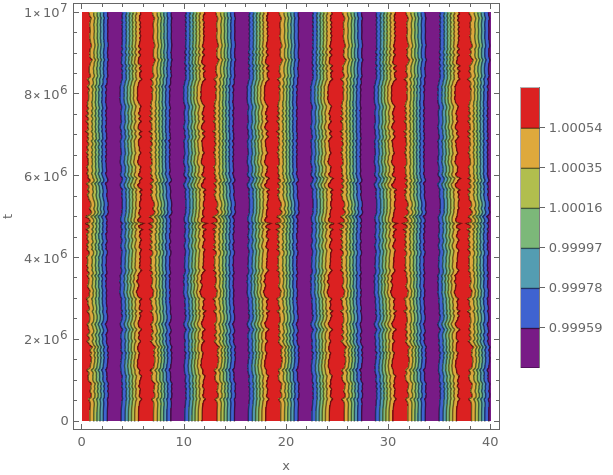}
        \caption{contour plot of $U_{AB}(x,t)$}
        \label{fig:UABcontour plot:chaosArjenSystem}
    \end{subfigure} 
    \vspace{0.1cm} 
    \begin{subfigure}[b]{0.48\textwidth}
        \centering
        \includegraphics[width=\textwidth]{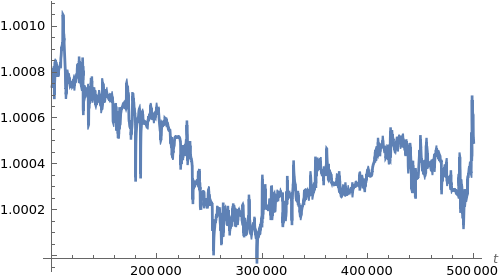}
        \caption{Plot of $U(40,t)$}
        \label{fig:UPlot:chaosArjenSystem}
    \end{subfigure}
    \hfill
    \begin{subfigure}[b]{0.48\textwidth}
        \centering
        \includegraphics[width=\textwidth]{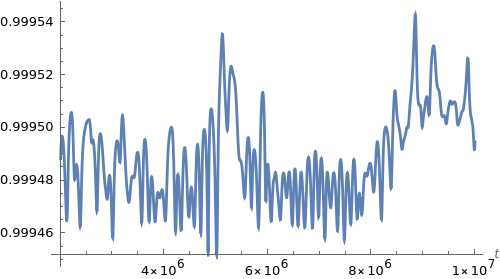}
        \caption{Plot of $U_{AB}(40,t)$}
        \label{fig:UABPlot:chaosArjenSystem}
    \end{subfigure}    
\caption{Numerical simulations comparing the evolution of a solution $U(x,t)$ of \eqref{eq:Arjen 6th order} (left column), with its AB-system approximation $U_{AB}(x,t)$ obtained from simulations of the AB-system \eqref{eq:Arjen system AB} (right column), where $U_{AB}(x,t)$ is similar to \eqref{eq:AB An}. The parameters in \eqref{eq:Arjen system AB} are set to $\alpha = 0.76$ and $\beta = 2$, which correspond (via \eqref{parsextendedmodel}) to $\gamma = -0.684...$, $\eta = 0.893...$ in \eqref{eq:Arjen 6th order}, and $d = 0.328...$ in \eqref{eq:Arjen system AB}. The simulation of AB-system \eqref{eq:Arjen system AB} is initialized with the marginally unstable plane wave with $K = \sqrt{4/5}$ and $R = R_t(K) - 0.01$ \eqref{defRtK}, and the simulation of \eqref{eq:Arjen 6th order} is initialized with the associated approximation \eqref{eq:AB An} and $\delta = 0.0005$. The simulation of \eqref{eq:Arjen system AB} is performed on the $\xi$-interval $(0,200 \pi/K)$ and that of \eqref{eq:Arjen 6th order} -- of which all parameters are now determined by \eqref{munuextended} -- on the associated $x$-interval (both with periodic boundary conditions). Only subintervals are shown (with identical lengths in $x$).}
\label{fig:chaos in Arjen system}
\end{figure}

We observed chaos-like behavior in our (limited!) simulations of the canonical AB-system only for values of $\alpha$ greater than $1/2$. We found in \Cref{subsec:rescaling the ab system:gs:s}, by rescaling, that $\alpha \equiv \frac12$ for all general higher-order scalar models of type \eqref{eq:general scalar eq:gs} (cf. \eqref{eq:rescaleABalphaetc}) and subsequently showed that $\alpha$ can be varied away from $\frac12$ in the broader class of scalar models 
\eqref{eq:illustration pde:rtab-Intro}/\eqref{eq:illustration pde:rtab} (cf. \eqref{eq:scalingsalphaetc-gen}). To test whether the observed chaos-like dynamics of AB-systems as shown in \Cref{fig:Intro-chaos,fig:ChaosAB} persists into the underlying PDE, we therefore considered the following model system of type \eqref{eq:illustration pde:rtab-Intro}/\eqref{eq:illustration pde:rtab}, that reduces to our basic model \eqref{eq:scalar example} when $\gamma = 0$:
\begin{equation}\label{eq:Arjen 6th order}
    \partial_t U = \mu U + 2U^2 - U^3 + \nu \partial_x^2 U + 2 \partial_x^4 U + \partial_x^6 U + \gamma U \partial_x^2U + \eta (\partial_x^2 U)^2.
\end{equation}
Naturally, the ODE problem is again given by:
\begin{equation*}
    \dot{u} = \mu u + 2 u^2 - u^3 \eqqcolon F(u;\mu),
\end{equation*}
yielding the same 3 spatially homogeneous steady states as in \eqref{model-homstates}. The eigenvalue problem, however, is now given by:
\begin{equation*}
    \omega(k; \mu,\nu) = \mu + 4u^+ - 3\left(u^+\right)^2 - \left(\nu + \gamma u^+\right)k^2 + 2k^4 - k^6 \eqqcolon F_u\left(u^+;\mu\right) + G(u^+; k,\nu).
\end{equation*}
We are interested in a Turing-fold bifurcation, i.e., we need to find $(\nu,k) = \left(\nu^*,k^*\right)$ such that $G^* = G_k^* = 0$; a straightforward shows:
\begin{equation}
\nonumber
    k^* = 1 \quad \text{ and } \quad \nu^* = 1 - \gamma.
\end{equation}
Consequently, at $(\mu^*, \nu^*) = (-1, 1 - \gamma)$ a Turing-fold bifurcation occurs. As usual, we set $\nu = \nu^* - \delta$ and find:
\begin{equation*}
    u^t = u^* + \tilde{u}\delta + \mathcal{O}(\delta^2), \quad 
    \mu^t = \mu^* + \hat{\mu}\delta^2 + \mathcal{O}(\delta^3), \quad \text { and } \quad 
    k^c = k^* + \tilde{k}\delta + \mathcal{O}(\delta^2),
\end{equation*}
with:
\begin{equation}
\nonumber
    \tilde{u} = \frac{1}{2+\gamma}, \quad \hat\mu = \frac{1}{(2+\gamma)^2}, \quad \text{ and } \quad \tilde{k} = \frac{1}{2(2 + \gamma)}.
\end{equation}
For $\mu = \mu^t - r \delta^2$, the associated AB-system is given by:
\begin{equation}
\nonumber
    \begin{cases}
        A_\tau = 4 A_{\xi\xi} + A - (2 + \gamma) AB, \\
        B_\tau = B_{\xi\xi} + \frac{1}{(2 + \gamma)^2} - r - B^2 -2(1 + \gamma - \eta) |A|^2.
    \end{cases}
\end{equation}
which, under the rescalings
\begin{equation}
\nonumber
    B = \frac{\tilde{B}}{2 + \gamma}, \quad r = \frac{R}{(2+\gamma)^2}, \quad \text{ and } \quad \xi = 2 \tilde{\xi},
\end{equation}
is scaled into its canonical form:
\begin{equation}\label{eq:Arjen system AB}
    \begin{cases}
         \phantom{(2 + \gamma)}A_\tau = A_{\tilde{\xi}\tilde{\xi}} + A - A \tilde{B}, \\
        (2 + \gamma) \tilde{B}_\tau = \frac{(2 + \gamma)}{4}\tilde{B}_{\tilde{\xi}\tilde{\xi}} +1 - R- \tilde{B}^2 - 2 (2 +\gamma)^2(1 + \gamma - \eta) |A|^2,
    \end{cases}
\end{equation}
Thus, for
\begin{equation}
\label{munuextended}
\mu = -1 + \frac{1-R}{(2+\gamma)^2} \delta^2 + \mathcal{O}(\delta^3), \; \; \nu = 1 - \gamma - \delta,  
\end{equation}
and $0 < \delta \ll 1$, the parameters of \eqref{eq:Arjen 6th order} are related to those of the AB-system \eqref{eq:canonical form:cs} in canonical form by:
\begin{equation}
\label{parsextendedmodel}
\alpha = \frac{1}{2+\gamma}, \quad d = \frac{1}{4}(2+\gamma), \quad \text{ and } \quad
\beta = - 2 (2 +\gamma)^2(1 + \gamma - \eta).
\end{equation}
In \Cref{fig:chaos in Arjen system}, we present the outcomes of two simulations. The first column shows that the extended model system \eqref{eq:Arjen 6th order} indeed can exhibit the spatio-temporal kind of chaos we already found for the AB-system. The parameters of \eqref{eq:Arjen 6th order} used for this simulation are based on those of the $AB$-system \eqref{eq:Arjen system AB} (see \eqref{munuextended} and \eqref{parsextendedmodel}). The second column of \Cref{fig:chaos in Arjen system} shows the outcome of the corresponding simulations of \eqref{eq:Arjen system AB}, presented in the form of $U_{AB}(x,t)$ as defined in \eqref{eq:AB An}. 

These simulations indicate that the complex dynamics exhibited by the AB-system of coupled modulation equations near a co-dimension 2 Turing-fold bifurcation are indeed reflected by the underlying model. In other words, the intricate dynamics shown in \Cref{fig:Intro-chaos,fig:ChaosAB,fig:chaos in Arjen system} can be expected to be observable near a Turing-fold bifurcation in large classes of coupled systems of evolutionary PDEs (on cylindrical domains). However, the second row of \Cref{fig:chaos in Arjen system} also hints at a subtle but potentially significant discrepancy between the full system and its modulation approximation: the time scale of the {\it turbulent flow} in the underlying model seems differ from the $\tau$ time-scale of its associated AB-system. It is not clear how incorporating higher-order terms in the Ansatz (cf. \eqref{eq:Ansatz:abd:rd}) used to derive $U_{AB}(x,t)$ \eqref{eq:AB An} could account for this discrepancy, since all terms in the Ansatz evolve on the same $\tau$-scale as the AB-system. Moreover, the fine-scale structure (or amplitude) of the {\it turbulent irregularities} may be sufficiently small to remain consistent with standard leading-order modulation approximations over finite times (see \cite{mielke2002ginzburg, schneider2017nonlinear} and the references therein). From an intuitive point of view, it is conceivable that the {\it slow turbulence} captured by the AB-system {\it cascades} to shorter time-scales in the underlying system. This would be an intriguing novel distinction between the underlying system and the approximating system of modulation equations -- one that can only arise in settings where both systems exhibit irregular dynamics (see \Cref{sec:discussion}).

\section{Conclusion and Discussion}
\label{sec:discussion}
In this work, we have investigated the formation and dynamics of patterns near a Turing-fold bifurcation -- i.e., the co-dimension 2 bifurcation at which a Turing and a saddle-node/fold/tipping bifurcation coincide, as can naturally appear in (systems of) evolutionary PDEs (defined on an unbounded domain). Thus, a description of the Turing-fold point involves at least two bifurcation parameters, $\mu$ and $\nu$, together with the bifurcation point $(\mu^*, \nu^*)$ characterized by the tangential intersection of two bifurcation curves in the $(\mu, \nu)$-plane: a curve $\Gamma^s$ of spatially homogeneous saddle-node points that give rise to a branch of homogeneous {\it background states} $u^+$ that are stable against non-spatial (ODE) perturbations, and a curve $\Gamma^t$ of non-degenerate Turing bifurcations that destabilize $u^+$. In other words that adhere to the chosen sign conventions, for $0 < \nu^* - \nu \ll 1$ the ODE-stable branch $u^+$, created via the saddle-node bifurcation, loses stability through a Turing bifurcation with critical wave number $k^c(\nu)$ -- under the natural non-degeneracy assumption that $k^c(\nu^*) \neq 0$.

We developed our approach by first showing -- better: confirming -- that the standard Ginzburg–Landau formalism is insufficient to capture the full dynamics associated with the Turing-fold bifurcation, since it cannot incorporate the interaction of the Turing bifurcation with the fold point. However, this analysis guided our reshaping of the standard formalism, which lead to the derivation of a new system of coupled modulation equations coined the AB-system, for a complex amplitude $A(\xi, \tau)$ representing the Turing bifurcation and a real amplitude $B(\xi, \tau)$ representing the saddle-node/fold bifurcation. We set up this procedure (and thus derived the AB-system) in the context of systems of reaction-diffusion equations and of higher-order scalar phase-field models. In both cases, we carefully listed the required non-degeneracy conditions associated with this co-dimension 2 bifurcation. We showed that any AB-system obtained from a formal derivation can be reduced to a canonical form with only three independent coefficients. Moreover, also clarified its relationship to the `zoomed-in' Ginzburg-Landau equation associated (only) with the Turing bifurcation. Naturally, we analyzed the existence and stability of small-amplitude spatially periodic `plane wave' patterns corresponding directly to their (classical) counterparts in the Ginzburg-Landau equation.

This allowed us to answer the motivating question posed in the Introduction:  {\it Under which conditions do stationary spatially-periodic `plane waves' that emerge at the Turing bifurcation persist as stable patterns beyond the fold point -- i.e., when can a tipping point be evaded through pattern formation?} Our analysis shows that this is completely determined by the sign of the coefficient $\beta$ in the AB-system -- i.e., stable patterns persist beyond the fold -- if and only if $\beta > 0$ (\Cref{thm:close enough:t:cs}).  Moreover, this sign can be explicitly determined, and we provided an explicit expression for $\beta$ in the fully general setting of a Turing-fold bifurcation in $n$-component system of reaction-diffusion equations. Furthermore, the signs of the $\beta$ coefficient of the AB-system and that of the Landau coefficient of the `zoomed-in' Landau equation must be reversed: stable patterns persist beyond the fold point if and only if the associated Turing bifurcation is supercritical (\Cref{thm:supercrit:t:cs}). This also yields an unexpected bonus: the sign of $\beta$, or equivalently that of the Landau coefficient near the Turing-fold point, is orders of magnitude easier to compute than the Landau coefficient in the general setting (cf. \eqref{eq:beta:rabtgl:rd} with \eqref{eq:full landau coefficient:rd}). All in all, our methods enable one to conclude whether a tipping point can be evaded in a given ecosystem model, if the Turing and saddle-node bifurcations take place for parameter combinations that can be assumed to be `sufficiently close'.  

Naturally, we followed up our analysis with a numerical investigation of the dynamics of the (supercritical) AB-system, with a focus on the question whether the AB-system indeed gives a reliable approximation of the dynamics of the underlying system (near a Turing-fold point). We demonstrated excellent agreement between the theoretically predicted and numerically observed convergence rate of the $L^2$-norm difference between the small-amplitude periodic patterns associated with the Turing bifurcation in the underlying scalar model system and the corresponding approximations provided by the AB-system. Moreover, we found that even the subtle transient pattern reselection dynamics -- where unstable plane waves outside the Busse balloon transition back into the region of stable solutions -- were reproduced with remarkable accuracy. Finally, we showed (numerically) that the AB-system may exhibit complex, spatio-temporal chaos-like dynamics and that qualitatively similar dynamics can be observed in the underlying model -- under identical parameter values.

The nature of our approach, and its similarity to that by which the Ginzburg-Landau equation is obtained, strongly suggests that the approximation of the PDE-dynamics near a co-dimension 2 Turing-fold point by an AB-system can be expected to be as widely applicable as that of the Ginzburg-Landau equation near a co-dimension 1 Turing-type instability. In fact, our methods are not special for reaction-diffusion systems or phase-field models; they can be directly employed in the context of other types of models, and our present focus on these models does not represent a true limitation, but rather is chosen to enhance the transparency and accessibility of our research.

Based on our present insights, we identify three promising directions for future research. 

First, there are various relevant bifurcations similar to the Turing-fold bifurcation for which coupled systems of modulation equations like the AB-system can (and should) be systematically derived and analyzed. Naturally, the Turing bifurcation could also interact with other types of classical ODE-type bifurcations, such as the transcritical and pitchfork bifurcation or higher-order catastrophes like the cusp, swallowtail, etc. \cite{homburg2024bifurcation, montaldi2021singularities}. Moreover, especially from the point of view of our ecosystem-type applications there are two relevant co-dimension 3 scenarios to be unraveled. 
As was already discussed in \Cref{rem:backgroundFig1}, the fate of a complex system that is triggered by a subcritical Turing bifurcation (i.e., $\beta < 0$) is more subtle that the plot in Fig. \ref{fig:Intro-TTT} suggests: the Turing bifurcation may also enable the system to evade tipping in this case, by triggering a transition to a spatial pattern with $\mathcal{O}(1)$ amplitude (see also \cite{van2025vegetation}). This is directly related to the local shape of the Busse balloon and thus it is crucial to study the transition from a subcritical ($\beta < 0$) to a supercritical ($\beta > 0$) Turing-fold bifurcation, i.e. the co-dimension 3 transition between tipping and tipping evasion -- see also \Cref{rem:degenerate}. Furthermore, many $n$-component ecosystem models are conceptually reduced to two-component models --see \cite{meron2015nonlinear, rietkerk2021evasion} and the references therein -- and these models cannot satisfy our generic non-degeneracy condition $k^*>0$ -- see \Cref{rem:k*=0}. From that point of view it thus is also important to study the co-dimension 3 degeneration at which $k^c \to k^* = 0$ as the Turing-fold point is approached. 

Second, as was already noticed in \Cref{rem:validity}, our analysis needs to be rigorously underpinned by justification results akin to the classical results on the validity of the Ginzburg-Landau equation of \cite{mielke2002ginzburg, schneider2017nonlinear} and the references therein. Especially in the simplest setting of higher-order phase-field models, work in progress indicates that a finite-time validity result can be obtained for the AB-system by following this literature. It may therefore, for instance, be more challenging to (also) study the validity of the various intermediate extended AB-systems that together govern the above co-dimension 3 scenario of the transition between tipping and tipping evasion (\Cref{rem:degenerate}). Another intriguing subject of further study in this context could be the relation between the characters of (small amplitude) spatio-temporal chaos in the underlying system and that of its approximating (system of) modulation equation(s) -- see \Cref{subsubsec:chaos:pbab:s} and especially \Cref{fig:chaos in Arjen system}. Such a study could start with a careful numerical exploration of this same relation in the context of a complex Ginzburg-Landau equation and a well-chosen (simple) underlying system. The complex Ginzburg-Landau equation may also exhibit `turbulent' dynamics  \cite{chate1994spatiotemporal}, but we are not aware of a thorough study in which the nature of these dynamics is compared to that of an underlying model: is there a similar distinction as in \Cref{fig:chaos in Arjen system}, i.e., does the underlying system exhibit dynamics at finer scales than those of the modulation equation?
If so, can this distinction be captured by a (higher order) validity analysis? 

Third, given the above claimed universality of the AB-system as normal form describing the dynamics of small amplitude patterns near a Turing-fold bifurcation in systems of evolutionary PDEs, it is natural to study the dynamics of patterns generated by the AB-system beyond the present most simple case of the existence and stability of plane wave patterns (\Cref{sec:coupled system}). Such a study could start with a local analysis of the Turing-type bifurcation that seems to be driving the richer types of dynamics -- quasi-periodic, chaos-like -- observed in the simulations shown in Sections \ref{sec:introduction} and \ref{sec:simulations}. Another approach could be to set up a dynamical systems type search for stationary patterns and/or traveling waves and combine that with a spectral stability analysis. (Both aspects are the subject of work in progress.)

{\bf Acknowledgment.} This research is supported by the ERC Synergy project RESILIENCE (101071417).
\bibliographystyle{plain}
\bibliography{bib}

\end{document}